\documentclass[10pt]{amsart}

\usepackage{xypic,amscd,amssymb,latexsym}
\usepackage{amsmath}	
\usepackage{graphicx}
\usepackage{pdfsync}

\usepackage{color}
\usepackage{cancel}

\usepackage{mathrsfs}

\usepackage{url}

\usepackage{changepage}

\begin{document}
\pagenumbering{arabic}

\newtheorem{theorem}{Theorem}[section]
\newtheorem{proposition}[theorem]{Proposition}
\newtheorem{lemma}[theorem]{Lemma}
\newtheorem{corollary}[theorem]{Corollary}
\newtheorem{remark}[theorem]{Remark}
\newtheorem{definition}[theorem]{Definition}
\newtheorem{question}[theorem]{Question}
\newtheorem{claim}[theorem]{Claim}
\newtheorem{conjecture}[theorem]{Conjecture}
\newtheorem{Prop}[theorem]{Definition and Proposition}
\newtheorem{example}[theorem]{Example}
\newtheorem{deflem}[theorem]{Definition and Lemma}

\def\qed{{\quad \vrule height 8pt width 8pt depth 0pt}}

\newcommand{\vs}[0]{\vspace{2mm}}

\newcommand{\til}[1]{\widetilde{#1}}

\newcommand{\mcal}[1]{\mathcal{#1}}

\newcommand{\ul}[1]{\underline{#1}}

\newcommand{\ol}[1]{\overline{#1}}

\newcommand{\wh}[1]{\widehat{#1}}

\newcommand{\ssum}[0]{{\textstyle \sum}}

\author{Hyun Kyu Kim}
\email{hyunkyukim@ewha.ac.kr, hyunkyu87@gmail.com}

\address{Department of Mathematics, Ewha Womans University, 52 Ewhayeodae-gil, Seodaemun-gu, Seoul 03760, Republic of Korea}

\numberwithin{equation}{section}

\title[Irreducible representations of quantum Teichm\"uller space and the phase constants]{Irreducible self-adjoint representations of quantum Teichm\"uller space and the phase constants}

\begin{abstract}
Quantization of the Teichm\"uller space of a non-compact Riemann surface has emerged in 1980's as an approach to three dimensional quantum gravity. For any choice of an ideal triangulation of the surface, Thurston's shear coordinate functions on the edges form a coordinate system for the Teichm\"uller space, and they should be replaced by suitable self-adjoint operators on a Hilbert space. Upon a change of triangulations, one must construct a unitary operator between the Hilbert spaces intertwining the quantum coordinate operators and satisfying the composition identities up to multiplicative phase constants. In the well-known construction by Chekhov, Fock and Goncharov, the quantum coordinate operators form a family of reducible representations, and the phase constants are all trivial. In the present paper, we employ the harmonic-analytic theory of the Shale-Weil intertwiners for the Schr\"odinger representations, as well as Faddeev-Kashaev's quantum dilogarithm function, to construct a family of irreducible representations of the quantum shear coordinate functions and the corresponding intertwiners for the changes of triangulations. The phase constants are explicitly computed and described by the Maslov indices of the Lagrangian subspaces of a symplectic vector space, and by the pentagon relation of the flips of triangulations. The present work may generalize to the cluster $\mathscr{X}$-varieties.
\end{abstract}

\maketitle

\vspace{-8mm}

\tableofcontents

\vspace{-5mm}

\section{Introduction}

Let $S$ be an oriented non-compact surface of finite type, with negative Euler characteristic: that is, a compact surface of genus $g$ minus $s$ punctures, with $s\ge 1$ and $2-2g-s<0$. The {\em Teichm\"uller space} of $S$ is the space of all complete hyperbolic metrics on $S$ up to isotopy. It is a universal cover of Riemann's classical moduli space, and has been an important object of study since early 20th century. It possesses nice geometric structures on itself, one of them being the Weil-Petersson Poisson structure. Quantization of the Teichm\"uller space has appeared as an approach to three dimensional quantum gravity since the relationship between the Teichm\"uller space and the theory of 3d gravity was pointed out in physics in 1980's, and fundamental constructions were first obtained in 1990's by Kashaev \cite{Kash} and independently by Chekhov-Fock \cite{CF} \cite{Fo}. Later in 2000's, the Chekhov-Fock quantum Teichm\"uller space was generalized by Fock-Goncharov \cite{FG09} to quantum cluster variety. 

\vs

One of the main results of the quantum Teichm\"uller theories is a family of projective unitary representations of the {\em mapping class group} ${\rm MCG}(S)$ on a Hilbert space, for each punctured surface $S$. Here, ${\rm MCG}(S) = {\rm Homeo}^+(S)/{\rm Homeo}(S)_0$ is the group of isotopy classes of all orientation-preserving homeomorphisms $S\to S$. One obtains a Hilbert space $\mathscr{H}_S$ and a map
$$
\pi^\hbar_S : {\rm MCG}(S) \to {\rm U}(\mathscr{H}_S)
$$
that is an `almost-group-homomorphism', where $\hbar$ is the real parameter for quantization (the Planck constant); to each $g\in {\rm MCG}(S)$, one assigns a unitary operator $\pi^\hbar_S(g)$ on $\mathscr{H}_S$ satisfying certain desired `quantization' properties, such that
$$
\pi^\hbar_S(g_1 g_2) = c^\hbar_{g_1,g_2} \, \pi^\hbar_S(g_1) \, \pi^\hbar_S(g_2)
$$
holds for every pair of elements $g_1,g_2\in {\rm MCG}(S)$, for some constants $c^\hbar_{g_1,g_2} \in \mathbb{C}^\times$. Due to the unitarity, these constants must live in ${\rm U}(1)$, and therefore are referred to as {\em phase constants} by the author in previous papers. The map
$$
{\rm MCG}(S) \times {\rm MCG}(S) \to {\rm U}(1) ~ : ~ (g_1,g_2) \mapsto c^\hbar_{g_1,g_2}
$$
is a group $2$-cocycle, hence represents a class in the second cohomology group $H^2({\rm MCG}(S);{\rm U}(1))$. Although this cohomology class did not gain much interest from Chekhov, Fock, and Goncharov, it turned out to be quite an interesting topic, and was an object of study of the papers by Funar-Sergiescu \cite{FS}, Funar-Kashaev \cite{FuKa}, the author \cite{K16a} \cite{K16c}, and Xu \cite{Xu}. The class depends on whether the relevant quantum Teichm\"uller theory is of the Chekhov-Fock(-Goncharov) type or of the Kashaev type, and also on the construction of the representations $(\mathscr{H}_S,\pi^\hbar_S)$, for such a construction is not unique.

\vs

The works \cite{FS} and \cite{Xu} were based on the assumption that the above $2$-cocycle coming from the (Chekhov-)Fock-Goncharov quantum Teichm\"uller theory is of a certain form, in particular cohomologically non-trivial. However, in \cite{K16c} the author proved that this Fock-Goncharov representation of quantum Teichm\"uller space yields trivial constants $c^\hbar_{g_1,g_2} \equiv 1$. One of the motivations of the present paper is a proof of the existence of a Chekhov-Fock-Goncharov type representation of the quantum Teichm\"uller space for which the corresponding cohomology class in $H^2({\rm MCG}(S);{\rm U}(1))$ is non-trivial, possibly rescuing the results of \cite{FS} and \cite{Xu}.

\vs

We now go into a little bit more detail, which would also reveal other significant aspects of the present work. First, we note that the precise version of the Teichm\"uller space relevant to the Chekhov-Fock-Goncharov quantization is the {\em enhanced Teichm\"uller space} $\mathscr{X}^+(S)$, which is a certain $2^s$-to-$1$ branched cover of the Teichm\"uller space \cite{T}. Crucial for the quantization is a special coordinate system on $\mathscr{X}^+(S)$, and it requires the choice of an extra data called an {\em ideal triangulation} $T$ on $S$, which is a collection of {\em edges}, i.e. isotopy classes of simple paths running between punctures, that divide $S$ into {\em triangles}, i.e. regions bounded by three edges. Given an ideal triangulation $T$, for each edge $e\in T$, there is a so-called {\em shear coordinate} function 
$$
x_e = x_{e;T} : \mathscr{X}^+(S) \to \mathbb{R}
$$
studied by Thurston and Penner in 1980's \cite{T} \cite{P87} \cite{P12}. Indeed, the map $\mathscr{X}^+(S) \to \mathbb{R}^T$, $(h,O) \mapsto (x_e(h,O))_{e\in T}$, is a real analytic bijection, providing a global coordinate system for $\mathscr{X}^+(S)$. The Weil-Petersson Poisson structure on $\mathscr{X}^+(S)$ is particularly simple in terms of these coordinate functions; the Poisson brackets among them are constant functions
\begin{align}
\nonumber
\{x_e, x_f\} = \varepsilon_{ef}, \qquad \forall e,f\in T,
\end{align}
where the constants $\varepsilon_{ef} \in \{-2,-1,0,1,2\}$ encode the combinatorics of $T$ as follows:
\begin{align*}
\varepsilon_{ef} & = a_{ef} - a_{fe},  \\
a_{ef} & = \mbox{the number of corners of triangles of $T$ delimited by $e$ on the left and $f$ on the right.}
\end{align*}

\vs

A {\em deformation quantization} of $\mathscr{X}^+(S)$ can be formulated as a map that assigns to each smooth function $F$ on $\mathscr{X}^+(S)$ a self-adjoint operator $\wh{F}^\hbar$ on a Hilbert space $\mathscr{H}$ such that
$$
\mbox{1) $F\mapsto \wh{F}^\hbar$ is $\mathbb{R}$-linear, \qquad 2) $\wh{1} = {\rm Id}_\mathscr{H}$, \qquad 3) 
$
[\wh{F}^\hbar, \wh{G}^\hbar] = 2\pi {\rm i} \hbar \, \wh{\{F,G\}} + o(\hbar), 
\mbox{ as $\hbar \to 0$}.
$}
$$
Constructing such a map for all smooth functions is too difficult, so one begins with easy functions like the shear coordinate functions for $T$. The problem of coming up with a set of self-adjoint operators $\wh{x}^\hbar_e = \wh{x}^\hbar_{e;T}$, $e\in T$, on a Hilbert space $\mathscr{H}_T$ satisfying the commutation relations
\begin{align}
\label{eq:basic_commutation_relations}
[\wh{x}^\hbar_e, \wh{x}^\hbar_f] = 2\pi {\rm i} \hbar \, \varepsilon_{ef} \, \cdot {\rm Id}, \qquad \forall e,f\in T
\end{align}
has a well-known solution on the Hilbert space $L^2(\mathbb{R}^r, dt_1 \cdots dt_r)$ for suitable $r$, by $\mathbb{R}$-linear combinations of the operators $t_1,\ldots,t_r$ and ${\bf i}\frac{\partial}{\partial t_1},\ldots,{\bf i}\frac{\partial}{\partial t_r}$. Finding suitable linear combinations is just a simple linear algebra problem, and one canonical solution is found by Fock-Goncharov \cite{FG09}. 

\vs

However, there is no preferred choice of an ideal triangulation $T$ of the surface $S$, hence it is important to make sure that the answers for the deformation quantization problem for different ideal triangulations $T$ are compatible with each other. In fact, this compatibility under the change of triangulations is the key aspect of the entire quantum Teichm\"uller theory. Let $T,T'$ be any pair of ideal triangulations. Then the shear coordinate functions $x_{e';T'}$ ($e'\in T'$) for $T'$ are related to those $x_{e;T}$ ($e\in T$) for $T$ by certain coordinate change formulas. Thus the corresponding self-adjoint operators $\wh{x}^\hbar_{e';T'}$ on $\mathscr{H}_{T'}$ should be related to the operators $\wh{x}^\hbar_{e;T}$ on $\mathscr{H}_T$ by a quantum version of the classical coordinate change formulas. Hence one looks for a unitary operator
$$
{\bf K}^\hbar_{T,T'} ~:~ \mathscr{H}_{T'} \to \mathscr{H}_T
$$
that {\em intertwines} the two sets of quantum coordinate operators, in the sense that
$$
{\bf K}^\hbar_{T,T'} \, \wh{x}^\hbar_{e';T'} \, ({\bf K}^\hbar_{T,T'})^{-1} = \mbox{quantum coordinate change expression in terms of $\{ \wh{x}^\hbar_{e;T} \,|\, e\in T\} $}
$$
holds for all $e' \in T'$. In particular, one must first come up with the quantum coordinate change expression for the right hand side; making rigorous sense of the above intertwining equations in fact involves notions from the theory of the quantum cluster algebras as well as functional analysis \cite{FG09}. Moreover, there has to be a consistency among these unitary intertwining operators, namely
$$
{\bf K}^\hbar_{T,T'} \, {\bf K}^\hbar_{T',T''} = c^\hbar_{T,T',T''} \, {\bf K}^\hbar_{T,T''}
$$
should hold for each triple of triangulations $T,T',T''$, for some constants $c^\hbar_{T,T',T''} \in {\rm U}(1)$. A solution to this seemingly quite difficult problem of constructing such intertwiners ${\bf K}^\hbar_{T,T'}$ was given by Fock-Goncharov \cite{FG09} (building on the first version by Chekhov-Fock \cite{CF}), with a key ingredient being the quantum dilogarithm function of Faddeev-Kashaev \cite{FaKa} \cite{Fa} (reviewed in \S\ref{subsec:QD} of the present paper), which is a meromorphic function on the complex plane with a parameter $\hbar>0$, defined on the strip $|{\rm Im}z|<\pi(1+\hbar)$ as the contour integral formula (going back to Barnes 100 years ago \cite{B})
$$
\Phi^\hbar(z) = \exp\left( - \frac{1}{4} \int_\Omega \frac{e^{-{\rm i} pz}}{\sinh(\pi p) \sinh(\pi \hbar p)} \frac{dp}{p} \right)
$$
where $\Omega$ is a contour along the real line that avoids the origin by a small half-circle above the origin. As mentioned above, the author proved in \cite{K16c} that the phase constants $c^\hbar_{T,T',T''}$ are identically $1$.

\vs

In fact, what Fock-Goncharov constructed in \cite{FG09} is a family of  projective unitary representations of a quantum version of the so-called cluster $\mathscr{D}$-variety, which is a `symplectic double' of the relevant cluster $\mathscr{X}$-variety, while the enhanced Teichm\"uller space $\mathscr{X}^+(S)$ is a prominent example of a cluster $\mathscr{X}$-variety. On the Hilbert space $\mathscr{H}_T = L^2(\mathbb{R}^T, \wedge_{e\in T} \, da_e)$, they constructed self-adjoint operators $\wh{x}^\hbar_e$ on $\mathscr{H}_T$ for each $e\in T$, as well as their `opposite-copy' operators  $\wh{\til{x}}^\hbar_e$, $e\in T$, together with some more operators; that is, there were more than double the number of coordinate functions to quantize than in the case of $\mathscr{X}^+(S)$. These operators are obtained as linear combinations of the basic operators $\{ a_e, {\rm i} \frac{\partial}{\partial a_e} \, | \, e\in T\}$, with the coefficients canonically determined by the matrix $\varepsilon=(\varepsilon_{ef})_{e,f\in T}$ for $T$. Moreover, these self-adjoint operators together form an irreducible representation on $\mathscr{H}_T$. 

\vs

However, the subset of Fock-Goncharov operators $\{ \wh{x}^\hbar_e \, | \, e\in T\}$ corresponding to just one copy of cluster $\mathscr{X}$-variety, or $\mathscr{X}^+(S)$, does {\em not} form an irreducible representation on $\mathscr{H}_T$. The existence of irreducible representations of the quantum cluster $\mathscr{X}$-variety or the quantum Teichm\"uller space, i.e. irreducible systems of self-adjoint operators satisfying the relations in eq.\eqref{eq:basic_commutation_relations}, as well as the uniqueness up to unitary equivalence, is a well-known fact in functional analysis under the name of generalized Stone-von Neumann theorem \cite{vN} \cite[Thm.14.8]{Hall}. As mentioned in  \cite{FG09}, one way of obtaining such irreducible representations, enumerated by `central characters' $\lambda$ of the quantum cluster $\mathscr{X}$-variety, is to apply the spectral decomposition of $\mathscr{H}_T$ into the direct integral of smaller Hilbert space slices $\mathscr{H}_{T;\lambda}$, over which one has a control for almost every $\lambda$.  Nevertheless, an explicit construction of these irreducible representations $\mathscr{H}_{T;\lambda}$ has not been established in the quantum Teichm\"uller theory literature, and this is one of the tasks that the present paper is undertaking. For each given $T$ and $\lambda$, constructing one irreducible representation $\mathscr{H}_{T;\lambda}$ is not difficult; the only nuisance is that there seems to be no preferred canonical description of it as some $L^2(\mathbb{R}^r)$. In fact, what is not so straightforward to obtain, and hence is one of the key points of the present paper, is an explicit construction of the intertwining operators ${\bf K}^\hbar_{T,T' ; \lambda}$ between these representation spaces $\mathscr{H}_{T';\lambda}$ and $\mathscr{H}_{T;\lambda}$ per each change of ideal triangulations $T\leadsto T'$. 

\vs

A main ingredient we utilized to accomplish the above task is the construction of (Lion-)(Segal-)Shale-Weil projective representations of the symplectic groups \cite{Shale} \cite{Segal} \cite{Weil} \cite{Lion77} \cite{LV}, also known as the representations of the metaplectic groups. Chekhov-Fock already remarked in \cite{CF} that their intertwiners can be viewed as generalizations of these Shale-Weil representations; in the present paper we work this out concretely, and also compute the phase constants precisely. Here we briefly review the Shale-Weil construction; see \S\ref{sec:Schrodinger_representations_and_Weil_intertwiners} for a more detailed review, and the development of some more necessary statements and notions. From the skew-symmetric integer matrix $\varepsilon = (\varepsilon_{ef})_{e,f\in T}$ for an ideal triangulation $T$ or for a cluster $\mathscr{X}$-seed, or more generally from any skew-symmetric real matrix, consider a real vector space $V_T$ with a formal set of basis $\{{\bf x}_e \, | \, e\in T\}$ enumerated by $T$, equipped with a skew-symmetric bilinear form $B_T$ defined on the basis as
$$
B_T({\bf x}_e,{\bf x}_f) = \varepsilon_{ef}, \quad \forall e,f\in T.
$$
This space $(V_T,B_T)$ corresponds to the space of edge weights on a train track made from $T$, equipped with the Thurston intersection form, studied by Bonahon \cite{B96}. Let $\mathfrak{h}_T$ be the real Lie algebra on the vector space $\mathfrak{h}_T = V_T \oplus \mathbb{R}c$ with the Lie bracket given by
$$
[{\bf x}, {\bf y}] = B_T({\bf x}, {\bf y}) \cdot c, \quad \forall {\bf x}, {\bf y} \in V_T, \qquad
[c, {\bf x}]=0, \quad \forall {\bf x}\in V_T.
$$
This Lie algebra is referred to as a {\em generalized Heisenberg algebra} (\S\ref{subsec:V_B}), where usually the form $B_T$ is assumed to be symplectic, i.e. non-degenerate; in our case, it has a non-zero radical $V_T^\perp$ of dimension $s=$ the number of punctures. One must choose a central character $\lambda$, i.e. a linear map $V_T^\perp \to \mathbb{R}$; for the case of the enhanced Teichm\"uller space, it is equivalent to choosing a real number for each puncture. Then, for each choice of a {\em Lagrangian subspace} $\ell$ of $V_T$, i.e. a maximal subspace of $V_T$ s.t. $B_T({\bf x},{\bf y})=0$ for all ${\bf x},{\bf y} \in \ell$, a representation of $\mathfrak{h}_T$ on a Hilbert space $\mathscr{H}_\ell = L^2(\mathbb{R}^r)$ called the Schr\"odinger representation was constructed (see \S\ref{subsec:Schrodinger_representation_of_Heisenberg_algebra} and \S\ref{subsec:irreducible_quantum_representation}), while for each pair of Lagrangian subspaces $\ell,\ell'$ of $V_T$, a canonical unitary map
$$
{\bf F}_{\ell,\ell'} : \mathscr{H}_{\ell} \to \mathscr{H}_{\ell'}
$$
that intertwines the representations of $\mathfrak{h}_T$, called the {\em (Shale-)Weil intertwiner}, was constructed (see \S\ref{subsec:LSSW_intertwiner}), using natural arguments of harmonic analysis, generalizing the Fourier transforms. The analog of the Fourier inversion formula ${\bf F}_{\ell,\ell'} \, {\bf F}_{\ell',\ell} = {\rm Id}$ was established (see Prop.\ref{prop:Weil_intertwiner_involutivity}), and then for every triple of Lagrangian subspaces $\ell,\ell',\ell''$, the composition identity
\begin{align}
\label{eq:F_composition_identity_intro}
{\bf F}_{\ell'',\ell'''} \, {\bf F}_{\ell',\ell''} \, {\bf F}_{\ell,\ell'} = e^{ \frac{{\rm i} \pi }{4} \tau(\ell,\ell',\ell'')} \cdot {\rm Id}
\end{align}
was proved (see Prop.\ref{prop:F_composition_identity}), where $\tau(\ell,\ell',\ell'')$ is the signature of the so-called {\em Maslov index} $Q_{\ell,\ell',\ell''}$ constructed by Kashiwara \cite{Kashiwara} \cite{LV}, defined as the quadratic form on the vector space $\ell \oplus \ell' \oplus \ell''$ as
$$
Q_{\ell,\ell',\ell''}(x+x'+x'') := B_T(x,x')+B_T(x',x'')+B_T(x'',x), \qquad \forall x+x'+x'' \in \ell\oplus \ell' \oplus \ell''
$$
These constructions also yield a family of projective representations of the automorphism group of $(V_T,B_T)$ (see Prop.\ref{prop:projective_representation_of_Aut_V_B}), namely the symplectic group, and its resulting cohomology class was one of the major objects of study in this literature (see e.g. \cite{Weil} \cite{LV} \cite{Rao}).

\vs

Let us concentrate on our setting. Fix a central character $\lambda$, i.e. the choice of a real number for each puncture of $S$. For each ideal triangulation $T$ of $S$, and each Lagrangian subspace $\ell$ for $(V_T,B_T)$, we give an explicit construction of an irreducible representation of the quantum Teichm\"uller space, i.e. a set of self-adjoint operators $\{\wh{x}^\hbar_{e;\lambda;T} \, | \, e\in T\}$ on a Hilbert space $\mathscr{H}_{T;\lambda;\ell}$, acting irreducibly (\S\ref{subsec:irreducible_quantum_representation}). Then, for any different choice of ideal triangulation $T'$ and a Lagrangian subspace $\ell'$ for $(V_{T'},B_{T'})$, we formulate the problem of finding a unitary intertwiner (\S\ref{subsec:quantum_mutation_algebraic}--\S\ref{subsec:mutation_intertwiner_formulation_of_the_problem})
$$
{\bf K}^\hbar_{T,T';\lambda;\ell,\ell'} ~ : ~ \mathscr{H}_{T';\lambda;\ell'} \to \mathscr{H}_{T;\lambda;\ell},
$$
and construct an answer as a mixture of the Weil intertwiner ${\bf F}$ and the Fock-Goncharov intertwiner ${\bf K}$ (\S\ref{subsec:QD}, \S\ref{subsec:our_answer_for_mutation_intertwiner}, \S\ref{subsec:verification_of_intertwining_equations}). As mentioned, these representations and intertwiners had only been hinted but not constructed in the literature. Finally, by carefully keeping track of the interplay between the Weil intertwiners and the Fock-Goncharov intertwiners, we prove the consistency relations (\S\ref{subsec:equality_of_two_signed_decompositions}, \S\ref{sec:phase_constants_for_composition})
\begin{align}
\label{eq:intro_consistency_relations}
{\bf K}^\hbar_{T,T';\lambda;\ell,\ell'} \, {\bf K}^\hbar_{T,T';\lambda;\ell',\ell''} = c^\hbar_{T,T',T'';\lambda;\ell,\ell',\ell''} \, {\bf K}^\hbar_{T,T';\lambda;\ell,\ell''}
\end{align}
We explicitly compute the phase constants $c^\hbar_{T,T',T'';\lambda;\ell,\ell',\ell''}\in {\rm U}(1)$ (Thm.\ref{thm:main}: the main theorem) by using eq.\eqref{eq:F_composition_identity_intro} and the famous pentagon identity of the quantum dilogarithm function $\Phi^\hbar(z)$ (Prop.\ref{prop:Phi_hbar}(4)), as well as by proving and using basic lemmas about the Weil intertwiner ${\bf F}$ (developed in \S\ref{sec:Schrodinger_representations_and_Weil_intertwiners}). So these phase constants are related to both the Maslov indices of the Lagrangian subspaces of degenerate symplectic vector spaces and also the pentagon relations of the flips of ideal triangulations. Using the explicit description of the phase constants, one can deduce that they are cohomologically nontrivial, as assumed by \cite{FS} and \cite{Xu}. One may summarize the results of the present paper as follows:
\begin{theorem}[summary of the results]
Let $S$ be a punctured surface, and $\lambda$ be the choice of a real number for each puncture of $S$. For any ideal triangulation $T$ of $S$, and each Lagrangian subspace $\ell$ for $(V_T,B_T)$, the set of self-adjoint operators $\{\wh{x}^\hbar_{e;\lambda;T} \,|\,e\in T\}$ on a Hilbert space $\mathscr{H}_{T;\lambda;\ell}$ forms an irreducible representation of the quantum Teichm\"uller space of $S$. Per different choice of data $(T,\ell)$ and $(T',\ell')$, there exists a unitary intertwiner ${\bf K}^\hbar_{T,T';\lambda;\ell,\ell'} ~ : ~ \mathscr{H}_{T';\lambda;\ell'} \to \mathscr{H}_{T;\lambda;\ell}$ for these representations. The phase constants $c^\hbar_{T,T',T'';\lambda;\ell,\ell',\ell''} \in {\rm U}(1)$ appearing in the consistency relations eq.\eqref{eq:intro_consistency_relations} for these intertwiners are explicitly computed and are found to be nontrivial, in particular yielding a family of genuinely projective representations of the mapping class group of $S$ whose associated second cohomology class is nontrivial.
\end{theorem}

\vs

The present work can be viewed as a cross-breeding of the Shale-Weil theory of 1960's--1970's with the quantum Teichm\"uller theory of 1990's--2000's. We also note that our results on the quantum Teichm\"uller spaces may generalize to the quantum cluster $\mathscr{X}$-varieties (see the last section).

\vs

{\bf Acknowledgments.} This research was supported by Basic Science Research Program through the National Research Foundation of Korea(NRF) funded by the Ministry of Education(grant number 2017R1D1A1B03030230). This work was supported by the National Research Foundation of Korea(NRF) grant funded by the Korea government(MSIT) (No. 2020R1C1C1A01011151).

\section{The Schr\"odinger representations and the Weil intertwiners}
\label{sec:Schrodinger_representations_and_Weil_intertwiners}

This section is for reviewing the relevant literature and for establishing necessary notations and lemmas. Original papers are \cite{Shale} \cite{Segal} \cite{Weil}, but we mostly follow \cite{Lion77} and \cite{LV} suitably adapted for our purposes. Whenever we add a terminology or provide a proof that is not written in these sources, we indicate so; however, it is possible that these may be found in some of the vast number of works citing these original papers.

\subsection{A real vector space $V$ with a skew-symmetric bilinear form $B$}
\label{subsec:V_B}

Let $V$ be a finite dimensional real vector space, equipped with a skew-symmetric bilinear form $B$, which is not necessarily non-degenerate; so, the radical 
$$
V^\perp = \{v\in V ~|~ B(v,w) =0, ~ \forall w\in V\}
$$
is not necessarily the zero space. 

\begin{remark}
Often in the literature, only the case when $B$ is non-degenerate is dealt with; then $(V,B)$ is called a symplectic vector space. Good thing about \cite{Lion77} is that it deals with the degenerate case. The paper \cite{LV} provides more details overall, but it concerns only the non-degenerate case. Hence we go back and forth between these two references.
\end{remark}

To avoid the triviality, we will assume that $V \neq V^\perp$. To such a data $(V,B)$, we associate a real Lie algebra $\mathfrak{n} = \mathfrak{n}(V,B)$ as follows. The underlying vector space is
$$
\mathfrak{n} = V + \mathbb{R} c, 
$$
which is a vector space direct sum; there is a natural inclusion of $V$ into $\mathfrak{n}$, so the elements of $V$ will be viewed as elements of $\mathfrak{n}$ whenever appropriate. The Lie bracket of $\mathfrak{n}$ is given by the formulas
$$
[v,w] = B(v,w) \cdot c, \quad \forall v,w\in V, \quad\qquad [v,c] = 0, \quad \forall v\in V.
$$

Let $N = N(V,B)$ be a simply-connected real Lie group whose Lie algebra is isomorphic to $\mathfrak{n}$; it is well known that such $N$ exists and is unique up to isomorphism (see e.g. \cite[Thm.20.21]{Lee}). As usual, denote by
$$
\exp : \mathfrak{n} \to N, \qquad x \mapsto \exp(x)= e^x,
$$
the exponential map, which is a diffeomorphism in this case. 

\vs

In particular, $\mathfrak{n}$ is a two-step nilpotent real Lie algebra, and is a generalized version of the so-called \ul{\bf Heisenberg algebra}, while $N$ is a generalized version of the so-called \ul{\bf Heisenberg group}.

\begin{definition}
A \ul{\bf polarization} of this Lie algebra $\mathfrak{n}$ is a maximal abelian Lie subalgebra of $\mathfrak{n}$.
\end{definition}

\begin{definition}
A \ul{\bf Lagrangian} of the vector space $(V,B)$ equipped with the skew-symmetric bilinear form $B$ is a maximal totally isotropic subspace of $V$ with respect to $B$, i.e. a subspace $\ell$ of $V$ s.t. $\ell^\perp = \{v\in V ~|~ B(v,w)=0,\forall w\in \ell\}$ equals $\ell$.
\end{definition}

It is a simple observation that every Lagrangian of $(V,B)$ contains the radical $V^\perp$. For the purposes of our paper, it is convenient to deal with the subspace of $\ell$ complementary to the radical $V^\perp$.

\begin{definition}
An \ul{\bf essential Lagrangian} of the vector space $(V,B)$ is a subspace $\ell$ of some Lagrangian $\ell'$ of $(V,B)$ such that $\ell' = \ell + V^\perp$ and $\ell \cap V^\perp = 0$.
\end{definition}
Among the terms used in the present section, the `essential Lagrangian' is one of few non-standard ones. Note that, an essential Lagrangian determines a unique Lagrangian containing it, but, if $V^\perp \neq 0$, a Lagrangian does not determine a unique essential Lagrangian contained in it.

\vs

For any essential Lagrangian $\ell$ of $V$, consider the vector subspace 
$$
\mathfrak{h} = \mathfrak{h}_\ell := V^\perp + \ell + \mathbb{R}c
$$
of $\mathfrak{n}$. Then $\mathfrak{h}$ is a polarization of $\mathfrak{n}$, and every polarization arises this way. Let $\mathfrak{h} = \mathfrak{h}_\ell$ be a polarization of $\mathfrak{n}$, associated to an essential Lagrangian $\ell$ of $(V,B)$. Denote by $H = H(\mathfrak{h}) = H(\ell)$ be the Lie subgroup of $N$ corresponding to $\mathfrak{h}$, i.e. 
$$
H = \exp(\mathfrak{h}).
$$

\subsection{The Schr\"odinger representation of the Heisenberg group $N$, associated to an essential Lagrangian $\ell$}
\label{subsec:Schrodinger_representation_for_all}

Now, choose any $\mathbb{R}$-linear map
$$
f : V^\perp \to \mathbb{R},
$$
and choose an essential Lagrangian $\ell$ of $(V,B)$. Extend $f$ to a linear map $\til{f}$ on $\mathfrak{h} = V^\perp + \ell + \mathbb{R} c$
$$
\til{f} : \mathfrak{h} \to \mathbb{R}
$$
by setting its value on $\ell$ to be zero and
$$
\til{f}(c) := -1.
$$
This requirement $\til{f}(c) = -1$ might look somewhat artificial; in general $\til{f}(c)$ can be set to be any real number (for example, it is set to be $2\pi$ in \cite[\S1.2.2]{LV}), but $-1$ suits best for our purposes. Then define a character
$$
\chi = \chi_f = \chi_{f,\ell}: H \to {\rm U}(1)
$$
of $H = H(\ell)$ as 
$$
\chi(\exp(x)) = e^{{\rm i} \til{f}(x)}, \quad \forall x \in \mathfrak{h}.
$$
There is a freedom also in the normalization of this character $\chi$, but we choose to work with the above.

\begin{remark}
Lion \cite[p160]{Lion77} starts from a more general linear map $f:\mathfrak{n} \to \mathbb{R}$ from the beginning, instead of $f:V^\perp \to \mathbb{R}$, hence an essential Lagrangian is not necessary to define the character $\chi$ of $H$.
\end{remark}

We shall construct a Hilbert space $\mathscr{H}_\ell$ as follows. First, consider continuous functions $k : N \to \mathbb{C}$ satisfying the covariance relation
\begin{align}
\label{eq:k_covariance}
k(nh) = \chi(h)^{-1} k(n), \quad \forall g\in N, ~ \forall h \in H.
\end{align}
Since $|\chi(h)|=1$ for each $h\in H$, we get a well-defined function $|k| : N/H \to \mathbb{R}_{\ge 0}$, defined as $nH \mapsto |k(nh)|$, for any $n \in N$ and $h\in H$. We denote by $dm_N$ and $dm_{N/H}$ the Radon measures on $N$ and $N/H$ that are invariant under the action of $N$ from the left respectively, i.e. the Haar measures. It is well known that each of them exists and is unique up to multiplication by a positive real constant; we can see these more clearly and concretely once we choose a basis of $\mathfrak{n}$ and $\mathfrak{n}/\mathfrak{h}$, in which case the situation is like the measures on the Euclidean spaces invariant under the translations, and which we shall indeed do shortly.

\vs

One can now ask whether this function $|k|$ on $N/H$ is square-integrable with respect to the measure $dm_{N/H}$. Denote by ${\rm Fun}_H^N$ the set of all continuous functions $k:N\to \mathbb{C}$ that satisfy eq.\eqref{eq:k_covariance} and are square-integrable on $N/H$ with respect to the measure $dm_{N/H}$. Then we have $L^2(N/H,dm_{N/H})$-norm on ${\rm Fun}^N_H$, and define $\mathscr{H}_\ell$ to be the completion of this normed vector space. In the end, we have
$$
\mathscr{H}_\ell = \{ k : N \to \mathbb{C} ~ | ~ \mbox{$k$ is measurable, $k$ satisfies eq.\eqref{eq:k_covariance}, $|k|$ is square-integrable on $N/H$} \}/\sim
$$
where $\sim$ means that an element of $\mathscr{H}_\ell$ is defined up to a set of measure zero. Notice that the space $\mathscr{H}_\ell$ is independent of the choice of a multiplicative constant coming from the ambiguity of the measure $dm_{N/H}$.

\vs
 
If we choose a set $N_0 \subseteq N$ of representatives for all the left $H$-cosets, so that we have a bijection $N_0 \to N/H$, $n \mapsto nH$, then we obtain a natural identification map
\begin{align}
\label{eq:Hilbert_space_correspondence1}
\mathscr{H}_\ell \longleftrightarrow L^2(N/H, dm_{N/H})
\end{align}
Namely, let $k : N \to \mathbb{C}$ be a function representing an element of $\mathscr{H}_\ell$. Define a function $\ol{k} : N/H \to \mathbb{C}$ by
$$
\ol{k}(nH) := k(n), \quad \forall n\in N_0.
$$
Then $\ol{k}$ is well-defined up to a set of measure zero, and one can verify that it represents an element of $L^2(N/H, dm_{N/H})$. Conversely, let $k_0 : N/H \to \mathbb{C}$ be a function representing an element of $L^2(N/H, dm_{N/H})$. Define a function $\til{k}_0 : N \to \mathbb{C}$ by
$$
\til{k}_0 (n h) := \chi(h)^{-1} k_0(n), \qquad \forall n\in N_0, \quad \forall h \in H.
$$
One can verify that $\til{k}_0$ represents a well-defined element of $\mathscr{H}_\ell$, and that the above two correspondences between $\mathscr{H}_\ell$ and $L^2(N/H, dm_{N/H})$ are inverses to each other.

\vs

We now define a unitary representation of the group $N$ on $\mathscr{H}_\ell$, naturally induced by the left action of $N$ on $N$ and on $N/H$. Namely, to each $n \in N$ define an operator 
$$
\pi_\ell(n)  : \mathscr{H}_\ell \to \mathscr{H}_\ell
$$
by the formula
\begin{align}
\label{eq:pi_ell_action}
(\pi_\ell(n) k)(n') := k(n^{-1} n'), \qquad \forall k \in \mathscr{H}_\ell, \quad \forall n' \in N;
\end{align}
one can easily see that $\pi_\ell(n) k \in \mathscr{H}_\ell$ holds for all $k\in \mathscr{H}_\ell$, and can show that $\pi_\ell(n)$ is unitary (using the invariance of the measure $dm_{N/H}$). This representation $(\mathscr{H}_\ell, \pi_\ell)$ is called the \ul{\bf Schr\"odinger representation} of the generalized Heisenberg group $N = N(V,B)$, \ul{\bf associated to} the essential Lagrangian $\ell$ of $(V,B)$. We note that there is no need to choose a basis of $\ell$ so far.

\subsection{The Schr\"odinger representation of $N$, for an essential Lagrangian $\ell$ and an isotropic supplementary basis $\mathcal{B}$}
\label{subsec:Schrodinger_representation_of_N_for_Lagrangian_decomposition}

One could realize the Schr\"odinger representation $\mathscr{H}_\ell$ by the model $L^2(N/H, dm_{N/H})$ via the correspondence in eq.\eqref{eq:Hilbert_space_correspondence1}. Even more concretely, we will now realize this Hilbert space as $L^2(\mathbb{R}^r)$, by choosing suitable bases of vector spaces in $\mathfrak{n}$. We follow the idea of Lion \cite[\S A]{Lion77} of using `B.S.A.' (adapted supplementary basis) and restrict ourselves to more special choices.
\begin{definition}
\label{def:symplectic_decomposition}
For an essential Lagrangian $\ell$ of the vector space $(V,B)$ equipped with a skew-symmetric bilinear form, an \ul{\bf isotropic supplementary basis} to $\ell$ in $V$ is an ordered set $\mathcal{B}$ of elements $v_1,\ldots,v_r$ of $V$ such that the their span ${\rm span}_\mathbb{R} \mathcal{B}$ is an essential Lagrangian of $(V,B)$ and satisfies $\ell \cap {\rm span}_\mathbb{R} \mathcal{B} = 0$.

We say that the pair of such data $(\ell,\mathcal{B})$ is a \ul{\bf symplectic decomposition} of $(V,B)$.
\end{definition}

\begin{lemma}
One has:
\begin{enumerate}
\item[\rm (1)] For each essential Lagrangian $\ell$ of $(V,B)$, an isotropic supplementary basis $\mathcal{B} = \{v_1,\ldots,v_r\}$ to $\ell$ in $V$ exists. 
\item[\rm (2)] $V = V^\perp + \ell + {\rm span}_\mathbb{R} \mathcal{B}$. 

\item[\rm (3)] We have $r = \dim \ell = \frac{1}{2}(\dim V - \dim V^\perp)$. \qed
\end{enumerate}
\end{lemma}

In fact, having a symplectic decomposition is equivalent to having a `symplectic basis'.
\begin{lemma}
\label{lem:symplectic_decomposition_is_equivalent_to_symplectic_basis}
Let $(\ell,\mathcal{B})$ be a symplectic decomposition of $(V,B)$, with $\mathcal{B} = \{v_1,\ldots,v_r\}$. Then there exists a unique ordered basis $\{w_1,\ldots,w_r\}$ of $\ell$ such that $B(w_i, v_j) = \delta_{ij}$, $\forall i,j \in \{1,\ldots,r\}$.
\end{lemma} 
{\it Proof.} Straightforward exercise in linear algebra (we will not really use this lemma). \qed

\vs

Choose an essential Lagrangian $\ell$ of $(V,B)$, and consider the subgroup $H = H(\ell) = H(\mathfrak{h}_\ell)= \exp(V^\perp +\mathbb{R}c + \ell)$ of $N$. Let's now choose any isotropic supplementary basis 
$$
\mathcal{B} = \{v_1,\ldots,v_r\}
$$
to $\ell$ in $V$. Then we shall use the set
$$
N_0 = \{ \exp( \ssum_{i=1}^r t_i v_i ) ~|~ t_i \in \mathbb{R}, ~ \forall i=1,\ldots,r\} \subset N
$$
as the set of representatives of the left $H$-cosets in $N$. So, each element $n$ of $N$ can be written as
$$
n = \exp( \ssum_{i=1}^r t_i v_i ) \cdot h
$$
for unique $t_1,\ldots,t_r \in \mathbb{R}$ and unique $h \in H$. 

\vs

Consider now the following {\em Hilbert space associated to the symplectic decomposition} $(\ell,\mathcal{B})$ of $(V,B)$:
\begin{align}
\label{eq:H_ell_B}
\mathscr{H}_{\ell, \mathcal{B}} := L^2(\mathbb{R}^r, \, dt_1 \cdots dt_r),
\end{align}
where $t_1,\ldots,t_r$ are used as the real coordinate variables for $\mathbb{R}^r$, and $dt_1 \cdots dt_r$ denotes the usual Lebesgue measure on this $\mathbb{R}^r$. By the identification of $N_0$ and $N/H$, one has the natural correspondence
$$
L^2(N/H, dm_{N/H}) \longleftrightarrow \mathscr{H}_{\ell,\mathcal{B}}.
$$
So, by composition, we obtain a canonical identification map
\begin{align}
\label{eq:I}
{\bf I}_{\ell,\mathcal{B}} : \mathscr{H}_\ell \to \mathscr{H}_{\ell,\mathcal{B}}
\end{align}
of Hilbert spaces. Let's describe ${\bf I}_{\ell,\mathcal{B}}$ more concretely. Let $\varphi = \varphi(t_1,\ldots,t_r) \in \mathscr{H}_{\ell,\mathcal{B}}$. Then the corresponding element ${\bf I}_{\ell,\mathcal{B}}^{-1} \varphi \in \mathscr{H}_\ell$ is given by 
\begin{align}
\label{eq:I_ell_B_formula}
({\bf I}_{\ell,\mathcal{B}}^{-1} \varphi)(\exp(\ssum_{j=1}^r t_j v_j)\cdot h) = \chi(h)^{-1} \cdot \varphi(t_1,\ldots,t_r), \qquad \forall t_1,\ldots,t_r \in \mathbb{R}, ~ \forall h \in H.
\end{align}

\vs

Now we define $\pi_{\ell,\mathcal{B}}$ to be the unitary representation of $N$ on the Hilbert space $\mathscr{H}_{\ell,\mathcal{B}}$, induced by $\pi_{\ell}$ via the above map ${\bf I}_{\ell,\mathcal{B}}$, that is,
\begin{align}
\label{eq:pi_ell_B_definition}
\pi_{\ell,\mathcal{B}}(n) \circ {\bf I}_{\ell,\mathcal{B}} = {\bf I}_{\ell,\mathcal{B}} \circ \pi_\ell(n), \quad \forall n \in N.
\end{align}
This representation $(\mathscr{H}_{\ell,\mathcal{B}},\pi_{\ell,\mathcal{B}})$ is called the \ul{\bf Schr\"odinger representation} of the generalized Heisenberg group $N=N(V,B)$, \ul{\bf associated to} the symplectic decomposition $(\ell,\mathcal{B})$ of $(V,B)$.

\vs

What happens if we fix an essential Lagrangian $\ell$ and choose two different isotropic supplementary bases $\mathcal{B}$ and $\mathcal{B}'$ to $\ell$ in $(V,B)$? What is the relationship between the representations $(\mathscr{H}_{\ell,\mathcal{B}}, \pi_{\ell,\mathcal{B}})$ and $(\mathscr{H}_{\ell,\mathcal{B}'}, \pi_{\ell,\mathcal{B}'})$? Is there a natural map $\mathscr{H}_{\ell,\mathcal{B}} \to \mathscr{H}_{\ell,\mathcal{B}'}$? Indeed, by the construction of these representations, there is a canonical map
$$
{\bf I}_{\ell,\mathcal{B}'} \circ {\bf I}_{\ell,\mathcal{B}}^{-1} ~ : ~ \mathscr{H}_{\ell,\mathcal{B}} \longrightarrow \mathscr{H}_{\ell,\mathcal{B}'},
$$
which is unitary and intertwines the actions $\pi_{\ell,\mathcal{B}}$ and $\pi_{\ell,\mathcal{B}'}$ of $N$. Although the construction of this identification map looks natural and trivial, when we write the domain and the codomain as $L^2(\mathbb{R}^r)$, it is in general {\em not} the identity map on $L^2(\mathbb{R}^r)$. In fact, the naive map $\mathscr{H}_{\ell,\mathcal{B}} \to \mathscr{H}_{\ell,\mathcal{B}'}$ representing the identity map on $L^2(\mathbb{R}^r)$ also plays an important role in the present paper, and will be dealt with in \S\ref{subsec:projective_representation_of_the_symplectic_group}.

\subsection{The Schr\"odinger representation of the Heisenberg algebra $\mathfrak{n}$}
\label{subsec:Schrodinger_representation_of_Heisenberg_algebra}

One can also `differentiate' the representations $\pi_\ell$ and $\pi_{\ell,\mathcal{B}}$ of the preceding subsection to obtain representations of the generalized Heisenberg Lie algebra $\mathfrak{n}$, which we again denote by the same symbols $\pi_\ell$ and $\pi_{\ell,\mathcal{B}}$, defined by the following natural relations
$$
\pi_\ell(\exp(x))  = e^{\pi_\ell(x)} \quad\mbox{and}\quad \pi_{\ell,\mathcal{B}}(\exp(x)) = e^{\pi_{\ell,\mathcal{B}}(x)}, \qquad \forall x\in \mathfrak{n}.
$$
Some words must be put in order. The representation $\pi_\ell$ of the real Lie algebra $\mathfrak{n}$ here assigns to each element $x\in \mathfrak{n}$ an operator $\pi_\ell(x)$ on the Hilbert space $\mathscr{H}_\ell$ that is \ul{\bf skew self-adjoint}, i.e. ${\rm i}$ times a self-adjoint operator. In particular, $\pi_\ell(x)$ is in general a densely defined operator, and not defined on the whole space $\mathscr{H}_\ell$. Via the functional calculus for the self-adjoint operator ${\rm i}^{-1} \cdot \pi_\ell(x)$, the unitary operator $e^{\pi_\ell(x)}$ on $\mathscr{H}_\ell$ is well-defined. Likewise for $\pi_{\ell,\mathcal{B}}$ on $\mathscr{H}_{\ell,\mathcal{B}}$. These representations $\pi_\ell$ and $\pi_{\ell,\mathcal{B}}$ of the generalized Heisenberg Lie algebra $\mathfrak{n}$ are also called the \ul{\bf Schr\"odinger representations}. Here we are assuming some basic knowledge on functional analysis; see e.g. \cite[\S14]{Hall}.

\vs

Let us compute the operator $\pi_{\ell,\mathcal{B}}(x)$ for some simple examples of $x\in \mathfrak{n}$.
\begin{example}
\label{ex:operator_in_ell}
Let's compute $\pi_{\ell,\mathcal{B}}(w)$ for an element
$$
w \in \ell \subset V \subseteq \mathfrak{n}.
$$
Define the real numbers $\alpha_i \in \mathbb{R}$, $i=1,\ldots,r$, as
$$
\alpha_i = B(v_i,w), \quad \forall i=1,\ldots,r,
$$
where $\mathscr{B} = \{v_1,\ldots,v_r\}$. Let $\varphi = \varphi(t_1,\ldots,t_r) \in \mathscr{H}_{\ell,\mathcal{B}}$, and let $\til{\varphi} = {\bf I}_{\ell,\mathcal{B}}^{-1} \varphi \in \mathscr{H}_\ell$; see eq.\eqref{eq:I_ell_B_formula}. For each $t\in \mathbb{R}$, let us consider $\exp(t w) \in N$, and its action under $\pi_\ell$ first. We use the relation
\begin{align}
\label{eq:BCH}
\exp(a) \exp(b) = \exp(a+b+{\textstyle \frac{1}{2}}[a,b]) = \exp([a,b]) \exp(b) \exp(a), \qquad \forall a,b\in \mathfrak{n},
\end{align}
which follows e.g. from the Baker-Campbell-Hausdorff formula. In particular, we have
$$
\exp(a) \exp(d) = \exp(a+d) = \exp(d) \exp(a), \quad \forall a\in \mathfrak{n}, \quad \forall d\in V^\perp + \mathbb{R}c,
$$
so that the factor $\exp([a,b])$ in the former equation line commutes with the other factors ($\because [a,b] \in \mathbb{R}c$, $\forall a,b\in \mathfrak{n}$ in our setting). Note
\begin{align*}
( \pi_\ell(\exp(t w)) \til{\varphi})(\exp(\ssum_{j=1}^r t_j v_j)\cdot h) & = 
\til{\varphi} ( \exp(-tw) \exp(\ssum_{j=1}^r t_j v_j) \cdot h) \\
& = \til{\varphi}\left(\exp(\ssum_{j=1}^r t_j v_j) \exp(-tw) \exp([-tw,\ssum_{j=1}^r t_j v_j]) \cdot h\right) \\
& = \til{\varphi}\left(\exp(\ssum_{j=1}^r t_jv_j) \exp(-tw) \exp(t \ssum_{j=1}^r t_j \alpha_j \cdot c) \cdot h\right) \\
& = \chi(h)^{-1} \chi(\exp(t\ssum_{j=1}^r t_j \alpha_j \cdot c))^{-1} \chi(e^{-tw})^{-1} \til{\varphi}(\exp(\ssum_{j=1}^r t_jv_j) )  \\
& = \chi(h)^{-1} \, e^{{\rm i} \cdot (t \sum_{j=1}^r t_j \alpha_j)} \, \til{\varphi}(\exp(\ssum_{j=1}^r t_j v_j)),
\end{align*}
which means, in view of eq.\eqref{eq:I_ell_B_formula}--\eqref{eq:pi_ell_B_definition}, that
$$
(\pi_{\ell,\mathcal{B}} (\exp(tw)) \varphi)(t_1,\ldots,t_r) = e^{{\rm i} t \sum_{j=1}^r t_j \alpha_j} \cdot \varphi(t_1,\ldots,t_r)
$$
thus $\pi_{\ell,\mathcal{B}} (\exp(tw)) = e^{{\rm i} t \sum_{j=1}^r t_j \alpha_j}$, hence by `differentiating' we get
$$
\pi_{\ell,\mathcal{B}}(w) = {\rm i} \ssum_{j=1}^r t_j \alpha_j \qquad (\mbox{on }\mathscr{H}_{\ell,\mathcal{B}} = L^2(\mathbb{R}^r,dt_1\cdots dt_r)).
$$
Here, the operator $\ssum_{j=1}^r t_j \alpha_j$ is given as the multiplication by $\ssum_{j=1}^n t_j \alpha_j$ on a dense subspace of $\mathscr{H}_{\ell,\mathcal{B}}$ such as the Schwartz space, and then extended uniquely to its maximal domain of self-adjointness.
\end{example}

Let's now compute another example of the opposite extreme. 
\begin{example}
\label{ex:operator_in_span_B}
Let's compute $\pi_{\ell,\mathcal{B}}(v)$ for an element
$$
v = \ssum_{i=1}^r a_i v_i \in {\rm span}_\mathbb{R} \mathcal{B} \subseteq V \subseteq \mathfrak{n}
$$
of the span of the chosen isotropic supplementary basis $\mathcal{B}$ to $\ell$ in $V$; so $a_1,\ldots, a_r \in \mathbb{R}$. For any $t\in  \mathbb{R}$, consider the element $\exp(t v) \in N$, and its action under $\pi_\ell$. Let $\varphi \in \mathscr{H}_{\ell,\mathcal{B}}$, and let $\til{\varphi} = {\bf I}_{\ell,\mathcal{B}}^{-1} \varphi \in \mathscr{H}_\ell$. Note
\begin{align*}
( \pi_\ell(\exp(tv)) \til{\varphi})(\exp(\ssum_{j=1}^r t_j v_j)\cdot h) & = 
\til{\varphi} ( \exp(-tv) \exp(\ssum_{j=1}^r t_j v_j) \cdot h) \qquad (\because \mbox{eq.\eqref{eq:pi_ell_action}}) \\
& = \til{\varphi} (\ssum_{j=1}^r (t_j - t a_j) v_j ) \cdot h) \\
& = \chi(h)^{-1} \cdot \varphi(t_1 - ta_1,\ldots,t_r - ta_r)
\end{align*}
which means, in view of eq.\eqref{eq:I_ell_B_formula}, that
$$
(\pi_{\ell,\mathcal{B}}(\exp(tv)) \, \varphi)(t_1,\ldots,t_r) = \varphi(t_1- t a_1,\ldots,t_r-ta_r).
$$
So $\pi_{\ell,\mathcal{B}}(\exp(tv_i)) = e^{ - t \ssum_{j=1}^r a_j \frac{\partial}{\partial t_j}}$, hence by `differentiating' we have
$$
\textstyle \pi_{\ell,\mathcal{B}}(v) = \pi_{\ell,\mathcal{B}}(\ssum_{j=1}^r a_j v_j) = - \ssum_{j=1}^r a_j \, \frac{\partial}{\partial t_j} = {\rm i} \ssum_{j=1}^r a_j \left( {\rm i} \frac{\partial}{\partial t_j} \right) \qquad (\mbox{on }\mathscr{H}_{\ell,\mathcal{B}} = L^2(\mathbb{R}^r,dt_1\cdots dt_r)).
$$
Here, the operator $\ssum_{j=1}^r a_j \, {\rm i} \frac{\partial}{\partial t_j}$ is defined by the formula $\varphi \mapsto \ssum_{j=1}^r a_j {\rm i} \frac{\partial}{\partial t_i} \varphi$ for $\varphi$ living in a dense subspace like the Schwartz space, and then extended to a domain of self-adjointness.

\end{example}

\vs

The following two results can be deduced in a similar manner.
\begin{example}
\label{ex:operator_in_center}
One has
\begin{align*}
\pi_{\ell,\mathcal{B}}(v) & = {\rm i} \, f(v) \cdot {\rm Id}, \qquad \forall v \in V^\perp \subset \mathfrak{n}, \\
\pi_{\ell,\mathcal{B}}(c) & = - {\rm i} \cdot {\rm Id}.
\end{align*}
\end{example}
We note that Examples \ref{ex:operator_in_ell}, \ref{ex:operator_in_span_B} and \ref{ex:operator_in_center} can be partially found in \cite[\S1.2.5]{LV}. One consequence of the above computational results is that, via the (generalized) Stone-von Neumann theorem \cite{vN} \cite[Thm.14.8]{Hall}, we can deduce that the unitary (strongly continuous) representation $\pi_{\ell,\mathcal{B}}$ of the group $N$ on $\mathscr{H}_{\ell,\mathcal{B}} = L^2(\mathbb{R}^r)$ is irreducible in the sense that there is no closed invariant proper subspace, and is uniquely determined by $f:V^\perp \to \mathbb{R}$ up to unitary equivalence.

\subsection{The Weil intertwiner for a change of choice of $\ell$ and $\mathcal{B}$}
\label{subsec:LSSW_intertwiner}
 
Let $\ell_1,\ell_2$ be two essential Lagrangians of $(V,B)$. We will briefly review the construction of a canonical unitary operator
$$
{\bf F}_{\ell_1,\ell_2} : \mathscr{H}_{\ell_1} \to \mathscr{H}_{\ell_2}
$$
that intertwines the representations $\pi_{\ell_1}$ and $\pi_{\ell_2}$ of the Heisenberg Lie group $N$. This operator was established by Segal \cite{Segal}, Shale \cite{Shale}, and Weil \cite{Weil}, for the case when $V^\perp =0$, and for nilpotent Lie algebras $\mathfrak{n}$ by Lion \cite{Lion77}, which concerns our case. See those references (also \cite{LV}) for more details.

\vs

Let $H_i = H(\ell_i)$ and write $\chi_i = \chi_{f,\ell_i} : H_i \to \mathbb{C}$, for $i=1,2$. Formally, the sought-for intertwining operator is defined as the following formula: for $k \in \mathscr{H}_{\ell_1}$, i.e. for measurable functions $k : N \to \mathbb{C}$, satisfying the covariance relation $k(nh) = \chi_1(h)^{-1} k(n)$ for all $n\in N$ and $h \in H_1$, and square-integrable on $N/H_1$, define a new function ${\bf F}_{\ell_1,\ell_2} k : N \to \mathbb{C}$ as
\begin{align}
\label{eq:canonical_intertwiner_formula}
({\bf F}_{\ell_1,\ell_2} k)(n) = \int_{H_2/H_1\cap H_2} k(n h) \, \chi_2 (h) \, dm_{H_2/H_1\cap H_2}(h).
\end{align}
On the right hand side, $h$ ranges over the representatives of the left $(H_1\cap H_2)$-cosets in $H_2$; note that the integrand function $H_2/H_1\cap H_2 \to \mathbb{C} : h(H_1\cap H_2) \mapsto k(nh) \chi_2(h)$ is well-defined, due to the covariance relation of $k$ with respect to the elements of $H_1$, and due to the fact that $\chi_1$ and $\chi_2$ coincide on $H_1 \cap H_2$. The integral is taken with respect to a left Haar measure $dm_{H_2/H_1\cap H_2}$ on $H_2/H_1\cap H_2$, i.e. a Radon measure invariant under the left action of $H_2$; such a measure exists and is unique up to multiplication by a positive real constant. Another issue is on the convergence of the integral. It converges for $k$ living in a dense subspace of $\mathscr{H}_{\ell_1}$, say the Schwartz space, which is more concretely seen in the model $\mathscr{H}_{\ell_1, \mathcal{B}_1} = L^2(\mathbb{R}^r)$, for any chosen isotropic supplementary basis to $\ell_1$ in $V$. In e.g. \cite[Thm.1]{Lion77}, it is shown that ${\bf F}_{\ell_1,\ell_2}$ is a homeomorphism between the Schwartz spaces of $\mathscr{H}_{\ell_1}$ and $\mathscr{H}_{\ell_2}$, with respect to the topologies given by certain families of semi-norms, and is unitary when suitably normalized, hence can be extended to the whole $\mathscr{H}_{\ell_1}$. To remove the ambiguity coming from the choice of the measure $dm_{H_2/H_1\cap H_2}$, consider the operator $\frac{ {\bf F}_{\ell_1,\ell_2} }{||{\bf F}_{\ell_1,\ell_2}||}$ rescaled by its operator norm. This rescaled operator defines a unitary map $\mathscr{H}_{\ell_1} \to \mathscr{H}_{\ell_2}$, and is independent of the choice of the measure $dm_{H_1/H_1\cap H_2}$. From now on, we denote this canonically normalized unitary operator by the same symbol ${\bf F}_{\ell_1,\ell_2}$. It is also proved e.g. in \cite[Thm.1]{Lion77} that this unitary map intertwines the representations, in the sense that
$$
{\bf F}_{\ell_1,\ell_2} \circ \pi_{\ell_1} (n)  = \pi_{\ell_2}(n) \circ {\bf F}_{\ell_1,\ell_2}, \quad \forall n \in N;
$$
at least formally, it is straightforward to verify this from eq.\eqref{eq:canonical_intertwiner_formula}.

\vs

The formula eq.\eqref{eq:canonical_intertwiner_formula} immediately reminds us of the Fourier transform. Indeed, for suitable choices of $\ell_1$ and $\ell_2$, the above ${\bf F}_{\ell_1,\ell_2}$ does coincide with the ordinary Fourier transform; so this operator ${\bf F}_{\ell_1,\ell_2}$ is often referred to as a generalization of the Fourier transform. Of course, such can be explicitly checked using the models $\mathscr{H}_{\ell_1,\mathcal{B}_1}$ and $\mathscr{H}_{\ell_2,\mathcal{B}_2}$, which are realized as $L^2(\mathbb{R}^r)$. 

\vs

Let us be more precise. Let $\mathcal{B}_i$ be an isotropic supplementary basis to $\ell_i$ in $V$, for $i=1,2$. Define
$$
{\bf F}_{(\ell_1,\mathcal{B}_1), (\ell_2,\mathcal{B}_2)} : \mathscr{H}_{\ell_1,\mathcal{B}_1} \to \mathscr{H}_{\ell_2,\mathcal{B}_2}
$$
to be the unique operator making the following diagram to commute:
\begin{align}
\label{eq:Weil_intertwiner_for_ell_B}
\xymatrix@C+17mm{
\mathscr{H}_{\ell_1} \ar[r]^-{{\bf F}_{\ell_1,\ell_2} } \ar[d]_{{\bf I}_{\ell_1,\mathcal{B}_1}} & \mathscr{H}_{\ell_2} \ar[d]^{{\bf I}_{\ell_2,\mathcal{B}_2} } \\
\mathscr{H}_{\ell_1,\mathcal{B}_1} \ar[r]^-{{\bf F}_{(\ell_1,\mathcal{B}_1), (\ell_2,\mathcal{B}_2)}} & \mathscr{H}_{\ell_2,\mathcal{B}_2}
}
\end{align}
Then it is easily seen to be unitary, and to be intertwining the representations $\pi_{\ell_1,\mathcal{B}_1}$ and $\pi_{\ell_2,\mathcal{B}_2}$.

\vs

We refer to the above constructed operators ${\bf F}_{\ell_1,\ell_2}$ and ${\bf F}_{(\ell_1,\mathcal{B}_1), (\ell_2,\mathcal{B}_2)}$ as the \ul{\bf Weil interwiners} or the \ul{\bf Lion-Segal-Shale-Weil intertwiners}. The subscripts of a Weil intertwiner will be often omitted when they are clear from the domain and the codomain of the intertwiner. So we will sometimes write
$$
{\bf F} : \mathscr{H}_{\ell_1} \to \mathscr{H}_{\ell_2} \qquad\mbox{and}\qquad {\bf F} : \mathscr{H}_{\ell_1,\mathcal{B}_1} \to \mathscr{H}_{\ell_2,\mathcal{B}_2}
$$
to mean the Weil intertwiners ${\bf F}_{\ell_1,\ell_2}$ and ${\bf F}_{(\ell_1,\mathcal{B}_1),(\ell_2,\mathcal{B}_2)}$ respectively.

\subsection{The phase constants for a composition of the Weil intertwiners via the Maslov indices}
\label{subsec:phase_constants_for_the_Weil_intertwiners}

Now, for any two essential Lagrangians $\ell_1$ and $\ell_2$ of $(V,B)$, we have an intertwiner ${\bf F}_{\ell_1,\ell_2}$ from the representation $(\mathscr{H}_{\ell_1}, \pi_{\ell_1})$ of $N$ to the representation $(\mathscr{H}_{\ell_2}, \pi_{\ell_2})$.  Of our principal interest is a composition of these intertwiners, starting and ending at a same representation $(\mathscr{H}_\ell,\pi_\ell)$. By the (generalized) Stone-von Neumann theorem \cite[Thm.14.8]{Hall}, any unitary intertwiner from $(\mathscr{H}_\ell,\pi_\ell)$ to itself must be the identity operator times a complex scalar of modulus $1$, and the question is to exactly pin down this scalar, which we call a \ul{\bf phase constant}.

\vs

The first observation is the following trivial one:
\begin{align}
\label{eq:F_ell_ell}
{\bf F}_{\ell,\ell} = {\rm Id} : \mathscr{H}_\ell \to \mathscr{H}_\ell.
\end{align}
The second one is already quite non-trivial, whose proof involves certain case-by-case investigation.
\begin{proposition}[the involutivity of the Weil intertwiners; follows from {\cite[end of p173]{Lion77}}]
\label{prop:Weil_intertwiner_involutivity}
For any two essential Lagrangians $\ell_1$ and $\ell_2$ of $(V,B)$, one has
$$
{\bf F}_{\ell_2,\ell_1} \circ {\bf F}_{\ell_1,\ell_2} = {\rm Id} : \mathscr{H}_{\ell_1} \to \mathscr{H}_{\ell_1}.
$$
\end{proposition}
See also \cite{Weil} and \cite[Prop.1.4.7]{LV} for the cases $V^\perp = 0$. This result corresponds to the Fourier inversion formula.

\vs

The next step is then to determine the constant for ${\bf F}_{\ell_3,\ell_1} \circ {\bf F}_{\ell_2,\ell_3} \circ {\bf F}_{\ell_1,\ell_2}$. To describe the result, we first have to recall the following definition. 

\begin{definition}[Kashiwara \cite{Kashiwara} {\cite[p174]{Lion77}} {\cite[\S1.5]{LV}}]
\label{def:Kashiwara_index}
For any three Lagrangians $\ell_1,\ell_2,\ell_3$ of $(V,B)$, the \ul{\bf Maslov index} is the quadratic form $Q_{\ell_1,\ell_2,\ell_3}$ on the vector space $\ell_1\oplus \ell_2 \oplus \ell_3$ defined by
\begin{align}
\label{eq:Maslov_index}
Q_{\ell_1,\ell_2,\ell_3}( x_1 + x_2 + x_3) := B(x_1,x_2)+ B(x_2,x_3)+B(x_3,x_1), \qquad \forall x_1+x_2 + x_3 \in \ell_1\oplus \ell_2 \oplus \ell_3,
\end{align}
and the \ul{\bf Kashiwara index} $\tau(\ell_1,\ell_2,\ell_3) \in \mathbb{Z}$ is defined as
\begin{align}
\label{eq:tau}
\tau(\ell_1,\ell_2,\ell_3) = \mbox{the signature of the Maslov index $Q_{\ell_1,\ell_2,\ell_3}$;}
\end{align}
that is, the number of $+$ signs minus the number of $-$ signs in a diagonalized expression of the quadratic form $Q_{\ell_1,\ell_2,\ell_3}$. 

\end{definition}
Note that this definition is about the Lagrangians, instead of the essential Lagrangians. We adapt it to our setting as follows.
\begin{definition}
\label{def:Kashiwara_index2}
For any three essential Lagrangians $\ell_1,\ell_2,\ell_3$ of $V(B)$, define their Maslov index as that of the corresponding Lagrangians 
$$
Q_{\ell_1,\ell_2,\ell_3} := Q_{V^\perp + \ell_1, V^\perp + \ell_2, V^\perp + \ell_3}
$$
and the Kashiwara index $\tau(\ell_1,\ell_2,\ell_3)=$ the signature of $Q_{\ell_1,\ell_2,\ell_3}$ accordingly.
\end{definition}
In fact, the same formulas in eq.\eqref{eq:Maslov_index}--\eqref{eq:tau} just work for the essential Lagrangians.

\begin{proposition}[follows from {\cite[Thm.2]{Lion77}}] 
\label{prop:F_composition_identity}
Let $\ell_1,\ell_2,\ell_3$ be any three essential Lagrangians of $(V,B)$. Then
\begin{align}
{\bf F}_{\ell_3,\ell_1} \circ {\bf F}_{\ell_2,\ell_3} \circ {\bf F}_{\ell_1,\ell_2} = e^{\frac{{\rm i}\pi}{4} \tau(\ell_1,\ell_2,\ell_3)} \cdot {\rm Id} : \mathscr{H}_{\ell_1} \to \mathscr{H}_{\ell_1}.
\end{align}
\end{proposition}
See \cite[Thm.1.6.1]{LV} for the cases $V^\perp=0$.

\vs

In the present paper, to simplify the situation, we choose not to pay too much attention to the powers of $e^{\frac{{\rm i}\pi}{4}}$, because there will be another phase constant which will play a more significant role, and because it may be possible to get rid of all these powers $e^{\frac{{\rm i}\pi}{4}}$; see the very last section.

\begin{definition}[the equality of the unitary operators up to an integer multiple of $e^{{\rm i}\pi /4}$]
If two unitary operators ${\bf A}$ and ${\bf B}$ satisfy ${\bf A} = e^{ \frac{ {\rm i} \pi n }{4}} {\bf B}$ for some integer $n$, then we write
$$
{\bf A} \sim {\bf B}.
$$
\end{definition}

\begin{corollary}[the relationship between compositions of ${\bf F}$]
\label{cor:relationship_between_two_compositions_of_F}
Let $(\ell_1,\mathcal{B}_1), \ldots, (\ell_n,\mathcal{B}_n)$ and $(\ell_1',\mathcal{B}_1'),\ldots,(\ell_m',\mathcal{B}_m')$ be any symplectic decompositions of $(V,B)$, such that $(\ell_1,\mathcal{B}_1) = (\ell'_1,\mathcal{B}_1')$ and $(\ell_n, \mathcal{B}_n) = (\ell'_m, \mathcal{B}_m')$. Then
$$
{\bf F}_{(\ell_{n-1},\mathcal{B}_{n-1}),(\ell_{n},\mathcal{B}_{n})} \circ \cdots \circ
{\bf F}_{(\ell_1,\mathcal{B}_1),(\ell_2,\mathcal{B}_2)}
~\sim~ 
{\bf F}_{(\ell'_{m-1},\mathcal{B}_{m-1}'),(\ell'_{m},\mathcal{B}'_{m})} \circ \cdots \circ
{\bf F}_{(\ell'_1,\mathcal{B}'_1),(\ell'_2,\mathcal{B}'_2)}.
$$
Here, if $n=1$ (resp. $m=1$), the left hand side (resp. the right hand side) is set to be ${\rm Id}$.
\end{corollary}

{\it Proof.} By definition of ${\bf F}_{(\ell,\mathcal{B}),(\ell',\mathcal{B}')}$ in \S\ref{subsec:LSSW_intertwiner}, it suffices to show ${\bf F}_{\ell_{n-1},\ell_{n}} \cdots {\bf F}_{\ell_1,\ell_2} \sim {\bf F}_{\ell'_{m-1},\ell'_{m}} \cdots {\bf F}_{\ell'_1,\ell'_2}$. By applying Prop.\ref{prop:Weil_intertwiner_involutivity} several times, one can boil down the situation to proving
$$
{\bf F}_{\ell_{n-1},\ell_{n}} \cdots 
{\bf F}_{\ell_1,\ell_2}
~\sim~ {\rm Id}
$$
in the case when $\ell_n = \ell_1$. The case $n=1$ is observed in eq.\eqref{eq:F_ell_ell}, and the case $n=2$ in Prop.\ref{prop:Weil_intertwiner_involutivity}, so it remains to show this when $n\ge 3$. By Prop.\ref{prop:F_composition_identity} and Prop.\ref{prop:Weil_intertwiner_involutivity} one obtains
\begin{align}
\label{eq:F_composition_combine}
{\bf F}_{\ell_2,\ell_3} \circ {\bf F}_{\ell_1,\ell_2} \sim {\bf F}_{\ell_1,\ell_3}.
\end{align}
By using eq.\eqref{eq:F_composition_combine} repeatedly, one can reduce the length of the composition ${\bf F}_{\ell_{n-1},\ell_{n}} \cdots 
{\bf F}_{\ell_1,\ell_2}$ to $2$, finishing the proof (by Prop.\ref{prop:Weil_intertwiner_involutivity});
\begin{align*}
\underbrace{ {\bf F}_{\ell_{n-1},\ell_{n}} {\bf F}_{\ell_{n-2},\ell_{n-1}} }_{\mbox{\tiny apply eq.\eqref{eq:F_composition_combine}}} \, \cdots \, {\bf F}_{\ell_1,\ell_2}
~\sim~ \underbrace{ {\bf F}_{\ell_{n-2},\ell_n} {\bf F}_{\ell_{n-3},\ell_{n-2}} }_{\mbox{\tiny apply eq.\eqref{eq:F_composition_combine}}} \, \cdots \, {\bf F}_{\ell_1,\ell_2}
~\sim~ \ldots ~\sim~ {\bf F}_{\ell_2,\ell_n} \, {\bf F}_{\ell_1,\ell_2} ~ \sim ~ {\rm Id}. \qed
\end{align*}

\subsection{The projective representations of the symplectic group}
\label{subsec:projective_representation_of_the_symplectic_group}

We now review how the Weil intertwiners lead to a family of projective representations of the symplectic group. This result itself is not directly used in the present paper, but an intermediate lemma needed for this construction will be used later in our proof, so it is worth reviewing it. Besides, the main result of the present paper can be viewed as a generalization of this construction.

\vs

Consider the automorphism group of the vector space $(V,B)$ equipped with a skew-symmetric bilinear form $B$, defined as:
$$
{\rm Aut}(V,B) := \{ C: V \to V ~|~ \mbox{$C$ is linear, invertible, preserves $B$, and $C|_{V^\perp} = {\rm Id}_{V^\perp}$} \}.
$$
The last condition $C|_{V^\perp} = {\rm Id}_{V^\perp}$ is necessary for the construction to work uniformly for any choice of the (central chacacter) function $f:V^\perp \to \mathbb{R}$. It is easy to observe that this group is isomorphic to the usual symplectic group ${\rm Sp}(2r)$ (or denoted sometimes by ${\rm Sp}(r)$).

\vs

This group ${\rm Aut}(V,B)$ acts naturally on the set of all symplectic decompositions of $(V,B)$; if $(\ell,\mathcal{B})$ is a symplectic decomposition, it is easy to see that
$$
C(\ell,\mathcal{B}) := (C(\ell),C(\mathcal{B}))
$$
is also a symplectic decomposition.

\vs

Choose any symplectic decomposition of $(\ell,\mathcal{B})$ of $(V,B)$. Choose any linear map $f : V^\perp \to \mathbb{R}$, so that the constructions of the present section apply. In particular, we have the Schr\"odinger representation $(\mathscr{H}_{\ell,\mathcal{B}}, \pi_{\ell,\mathcal{B}})$ of the Heisenberg group $N$ and the Heisenberg algebra $\mathfrak{n}$. We will now construct a family of unitary representations $\rho$ of ${\rm Aut}(V,B)$ on $\mathscr{H}_{\ell,\mathcal{B}}$, i.e. a homomorphism
$$
\rho = \rho_{\ell,\mathcal{B}} : {\rm Aut}(V,B) \longrightarrow {\rm U}(\mathscr{H}_{\ell,\mathcal{B}}) = \{\mbox{all unitary maps $\mathscr{H}_{\ell,\mathcal{B}} \to \mathscr{H}_{\ell,\mathcal{B}}$}\}
$$

\vs

Let $C\in {\rm Aut}(V,B)$. Then $C(\ell,\mathcal{B})$ is a symplectic decomposition, hence we can consider the Schr\"odinger representation $(\mathscr{H}_{C(\ell,\mathcal{B})}, \pi_{C(\ell,\mathcal{B})})$, and the Weil intertwiner ${\bf F} : \mathscr{H}_{C(\ell,\mathcal{B})} \to \mathscr{H}_{\ell,\mathcal{B}}$. We now describe how to canonically identify the domain $\mathscr{H}_{C(\ell,\mathcal{B})}$ with $\mathscr{H}_{\ell,\mathcal{B}}$. For this and for later use, we consider a slightly more general situation.
\begin{definition}[the re-labling map ${\bf R}$ for a symplectic transformation]
\label{def:R}
Let $(V,B)$ and $(V',B')$ be vector spaces with skew-symmetric bilinear forms, with the fixed choice of functions $f : V^\perp \to \mathbb{R}$ and $f' : {V'}^\perp \to \mathbb{R}$ which are used to build the Schr\"odinger representations. Let $C : V' \to V$ be an invertible linear transformation respecting $B'$ and $B$, and also $f'$ and $f$ in the sense that $f' = f \circ C$.

Let $(\ell,\mathcal{B})$ and $(\ell',\mathcal{B}')$ be symplectic decompositions of $(V,B)$ and $(V',B')$ respectively. Suppose that
$$
(\ell,\mathcal{B}) = C(\ell',\mathcal{B}')
$$
holds, where $C(\ell',\mathcal{B}') := (C(\ell'),C(\mathcal{B}'))$. In this case, define
$$
{\bf R} = {\bf R}_{(\ell,\mathcal{B}),(\ell',\mathcal{B}')} = {\bf R}_{C(\ell',\mathcal{B}'),(\ell',\mathcal{B}')} ~ : ~ \mathscr{H}_{\ell,\mathcal{B}} = \mathscr{H}_{C(\ell',\mathcal{B}')} \longrightarrow \mathscr{H}_{\ell',\mathcal{B}'}
$$
as the map representing the identity map $L^2(\mathbb{R}^r) \to L^2(\mathbb{R}^r)$ when the domain $\mathscr{H}_{\ell,\mathcal{B}}$ and the codomain $\mathscr{H}_{\ell',\mathcal{B}'}$ are realized as $L^2(\mathbb{R}^r)$ as in their very constructions in eq.\eqref{eq:H_ell_B}.
\end{definition}
One remark that is implicitly used in the above definition is that $C : V'\to V$ naturally induces the identification $C : {V'}^\perp \to V^\perp$, which is easy to check. 

\vs

The following lemma is obvious: 
\begin{lemma}
\label{lem:R_composition_identity}
The map ${\bf R}$ is unitary, and satisfies the composition identity
\begin{align}
\label{eq:R_composition_identity}
{\bf R}_{(\ell',\mathcal{B}'),(\ell',\mathcal{B}'')} \circ {\bf R}_{(\ell,\mathcal{B}),(\ell',\mathcal{B}')} = {\bf R}_{(\ell,\mathcal{B}),(\ell'',\mathcal{B}'')}
\end{align}
or ${\bf R} \circ {\bf R} = {\bf R}$ in short, whenever these three ${\bf R}$'s are defined.
\end{lemma}

We omit the subscripts of ${\bf R}$ when they are clear. This map ${\bf R}$ can be viewed as a simple re-labling map, and can be defined in a more general situation; namely, we do not really need a transformation $C$, and $(V,B)$, $(V',B')$ can be any two vector spaces with skew-symmetric bilinear forms such that $\dim V - \dim V^\perp = \dim V' - \dim{V'}^\perp \ge 1$. 

\vs

However, for the case as in Def.\ref{def:R}, the map ${\bf R}$ has more natural meaning. 
\begin{lemma}[viewing ${\bf R}$ as a pullback map]
\label{lem:R_as_pullback_map}
When we view $\mathscr{H}_{\ell,\mathcal{B}}$ and $\mathscr{H}_{\ell',\mathcal{B}'}$ as $L^2({\rm span}_\mathbb{R} \mathcal{B})$ and $L^2({\rm span}_\mathbb{R} \mathcal{B}')$ via the identifications ${\rm span}_\mathbb{R} \mathcal{B} \cong \mathbb{R}^r$ and ${\rm span}_\mathbb{R} \mathcal{B}' \cong \mathbb{R}^r$ induced by their respective bases $\mathcal{B}$ and $\mathcal{B}'$, the map ${\bf R} : \mathscr{H}_{\ell,\mathcal{B}} \to \mathscr{H}_{\ell',\mathcal{B}'}$ can be understood as the pullback map along $C$:
$$
{\bf R} = C^* ~:~ \mathscr{H}_{\ell,\mathcal{B}} = L^2(\mathbb{R}^r) =  L^2({\rm span}_\mathbb{R} \mathcal{B}) \longrightarrow L^2({\rm span}_\mathbb{R} \mathcal{B}') = L^2(\mathbb{R}^r) = \mathscr{H}_{\ell',\mathcal{B}'}.
$$
\end{lemma}
{\it Proof.} Indeed, $C : V' \to V$ restricts to an isomorphism $C : {\rm span}_\mathbb{R} \mathcal{B}' \to {\rm span}_\mathbb{R} \mathcal{B}$, which, when viewed as an isomorphism $C : \mathbb{R}^r \cong {\rm span}_\mathbb{R} \mathcal{B}' \to {\rm span}_\mathbb{R} \mathcal{B} \cong \mathbb{R}^r$, equals the identity map ${\rm Id} : \mathbb{R}^r \to \mathbb{R}^r$; note that $C$ sends the $i$-th basis vector $v_i'$ of ${\rm span}_\mathbb{R} \mathcal{B}'$ to the $i$-th basis vector $C(v_i')$ of ${\rm span}_\mathbb{R} C(\mathcal{B}') = {\rm span}_\mathbb{R} \mathcal{B}$.  \qed

\vs

Sometimes we want $(\ell,\mathcal{B})$ and $(\ell', \mathcal{B}')$ to completely determine the relating map $C$, so we establish the following uniqueness lemma.
\begin{lemma}[the uniqueness of the isomorphism $C$]
\label{lem:uniqueness_of_C}
Let $(V,B)$, $(V',B')$, $f$, $f'$ be as in Def.\ref{def:R}, and let $(\ell,\mathcal{B})$ and $(\ell',\mathcal{B}')$ be symplectic decompositions of $(V,B)$ and $(V',B')$. Let $C_1,C_2 : V' \to V$ be invertible linear transformations respecting $B'$ and $B$, and also $f'$ and $f$. Suppose further that
$$
C_1(\ell',\mathcal{B}') = (\ell,\mathcal{B}) = C_2(\ell',\mathcal{B}'), \quad \mbox{and} \quad C_1|_{{V'}^\perp} = C_2|_{{V'}^\perp}.
$$
Then $C_1 = C_2$.
\end{lemma}

This is comparable to Lem.\ref{lem:symplectic_decomposition_is_equivalent_to_symplectic_basis}. We prove this one, as we will implicitly use it throughout the paper.

\vs

{\it Proof.} Note $V' = {V'}^\perp + \ell' + {\rm span}_\mathbb{R} \mathcal{B}'$, and $C_1$ coincides with $C_2$ on ${V'}^\perp + {\rm span}_\mathbb{R} \mathcal{B}'$ already. Now let $v' \in \ell'$. Then $C_1(v') \in \ell$ and $C_2(v') \in \ell$. Denote $\mathcal{B}' = \{v_1',\ldots,v_r'\}$, so that $\mathcal{B} = \{C_j(v_1'),\ldots,C_j(v_r')\}$ holds for $j=1,2$, as the equality of ordered sets. For each $i=1,\ldots,r$, observe $B(C_1(v'),C_1(v_i')) = B'(v',v_i') = B(C_2(v'),C_2(v_i'))$ holds. Hence $B(C_1(v')- C_2(v'),x) = 0$ for all $x\in V^\perp + {\rm span}_\mathbb{R}\mathcal{B}$. Meanwhile, since $V^\perp + {\rm span}_\mathbb{R} \mathcal{B}$ is a maximal totally isotropic subspace of $V$, it follows that $C_1(v') - C_2(v')=0$, so $C_1(v')=C_2(v')$. \qed

\vs

Here is another natural viewpoint on ${\bf R}$, involving the ${\bf I}$ map of eq.\eqref{eq:I}. 
\begin{lemma}[a version $\til{\bf R}$ of ${\bf R}$ just for the Lagrangians]
\label{lem:til_R}
In the setting of Def.\ref{def:R}, define the map
$$
\til{{\bf R}} = \til{{\bf R}}_{\ell,\ell'} = \til{{\bf R}}_{C(\ell'),\ell'} ~:~ \mathscr{H}_\ell = \mathscr{H}_{C(\ell')} \longrightarrow \mathscr{H}_{\ell'}
$$
as the unique map making the following diagram to commute:
\begin{align}
\label{eq:til_R_diagram}
\xymatrix@C+2mm{
\mathscr{H}_\ell \ar[d]_{ {\bf I}_{\ell,\mathcal{B}} } \ar[r]^-{\til{\bf R}} & \mathscr{H}_{\ell'} \ar[d]^{ {\bf I}_{\ell',\mathcal{B}'} } \\
\mathscr{H}_{\ell,\mathcal{B}} \ar[r]^-{{\bf R}} & \mathscr{H}_{\ell',\mathcal{B}'}
}
\end{align}
i.e. define $\til{\bf R}$ by the equation:
$$
{\bf R}_{C(\ell',\mathcal{B}'),(\ell',\mathcal{B}')} \, \, {\bf I}_{C(\ell',\mathcal{B}')} = {\bf I}_{\ell',\mathcal{B}'} \,\, \til{{\bf R}}_{C(\ell'),\ell'}
$$
Then $\til{\bf R}$ is given by the formula
\begin{align}
\label{eq:til_R_formula}
(\til{\bf R}\, k)(n') = k( e^C ( n') ), \qquad \forall k \in \mathscr{H}_\ell, \quad \forall n' \in N',
\end{align}
where $e^C : N' \to N$ is the map defined as
\begin{align}
\label{eq:e_C}
e^C(e^{v' + \alpha c'}) = e^{C(v') + \alpha c}, \qquad \forall v' \in V', \quad \forall \alpha \in \mathbb{R}.
\end{align}
\end{lemma}

We note that, in case $(V,B) = (V',B')$ and $V^\perp = 0$, the map $\til{\bf R}$ as defined in eq.\eqref{eq:til_R_formula} is written as the symbol $A(g)$ (for an element $g$ of the symplectic group) in \cite[\S1.6]{LV}. As the maps ${\bf R}, {\bf I}, \til{\bf R}$ are used crucially in the present paper, we try to present complete proofs of their basic properties, for not all such proofs can be easily found in the literature.

\vs

{\it Proof of Lem.\ref{lem:til_R}.} We will show that the map $\til{\bf R}$ defined by eq.\eqref{eq:til_R_formula} indeed makes the diagram in eq.\eqref{eq:til_R_diagram} to commute. First, it is easy to observe that the map $\mathfrak{n}' \to \mathfrak{n} : v'+\alpha c' \mapsto C(v')+\alpha c$ is a Lie algebra isomorphism; we see that $e^C: N' \to N$ is the corresponding Lie group isomorphism. One can easily see that $C({V'}^\perp) = V^\perp$, hence $C(\mathfrak{h}') = C({V'}^\perp + \ell'+\mathbb{R}c') = C({V'}^\perp) + C(\ell') + C(\mathbb{R}c')= V^\perp + \ell + \mathbb{R}c = \mathfrak{h}$, and therefore $e^C(H')=H$. So, in view of the definitions of $\chi =\chi_f : H \to {\rm U}(1)$ and $\chi'=\chi'_{f'} : H' \to {\rm U}(1)$, and because $f' = f \circ C$, it is easy to observe that
$$
\chi' = \chi \circ e^C : H' \to {\rm U}(1).
$$

Write $\mathcal{B}' = \{v_1',\ldots,v_r'\}$, so that $\mathcal{B} = C(\mathcal{B}') = \{C(v_1'),\ldots,C(v_r')\}$. Let $\varphi = \varphi(t_1,\ldots,t_r) \in \mathscr{H}_{\ell,\mathcal{B}} = L^2(\mathbb{R}^r,dt_1\cdots dt_r)$. Write $\mathscr{H}_{\ell',\mathcal{B}'} = L^2(\mathbb{R}^r,dt_1'\cdots dt_r')$. By the definition of the map ${\bf R}$ (Def.\ref{def:R}), the element ${\bf R} \varphi \in \mathscr{H}_{\ell',\mathcal{B}'}$ is given by
$$
({\bf R} \varphi)(t_1',\ldots,t_r') = \varphi(t_1',\ldots,t_r'), \qquad \forall t_1',\ldots,t_r' \in \mathbb{R}.
$$
In view of eq.\eqref{eq:I_ell_B_formula} for ${\bf I}_{\ell',\mathcal{B}'}$, the element ${\bf I}_{\ell',\mathcal{B}'}^{-1} {\bf R} \varphi \in \mathscr{H}_{\ell'}$ is then given by
\begin{align*}
({\bf I}_{\ell',\mathcal{B}'}^{-1} {\bf R}\varphi)(\exp(\ssum_{j=1}^r t_j' v_j')\cdot h') & = \chi'(h')^{-1} \cdot ({\bf R}\varphi)(t_1',\ldots,t_r'), \\
& = \chi'(h')^{-1} \cdot \varphi(t_1',\ldots,t_r'), \qquad \forall t_1',\ldots,t_r' \in \mathbb{R}, ~ \forall h' \in H'.
\end{align*}
Meanwhile, the element $\til{\bf R} {\bf I}_{\ell,\mathcal{B}}^{-1} \varphi \in \mathscr{H}_{\ell'}$ (where $\til{\bf R}$ defined by eq.\eqref{eq:til_R_formula}) is given by 
\begin{align*}
(\til{\bf R} {\bf I}_{\ell,\mathcal{B}}^{-1} \varphi )(\exp(\ssum_{j=1}^r t_j' v_j')\cdot h')
& = ({\bf I}_{\ell,\mathcal{B}}^{-1} \varphi)( e^C(\exp(\ssum_{j=1}^r t_j' v_j')\cdot h'))
\end{align*}
We observe that $e^C : N' \to N$ defined by eq.\eqref{eq:e_C} is a group homomorphism. Indeed, for any $v',w' \in V'$ and $\alpha,\beta\in \mathbb{R}$, note
\begin{align*}
& e^C(e^{v'+\alpha c'} e^{w'+ \beta c'}) \underset{\mbox{\tiny eq.\eqref{eq:BCH}}}{=} e^C(e^{v'+w'+(\alpha+\beta)c'+\frac{1}{2}B'(v',w') c'}) 
\underset{\mbox{\tiny eq.\eqref{eq:e_C}}}{=} e^{C(v'+w')+(\alpha+\beta)c + \frac{1}{2}B'(v',w') c}  \\
& \qquad = e^{C(v')+C(w')+(\alpha+\beta)c + \frac{1}{2}B(C(v),C(w)) c} \underset{\mbox{\tiny eq.\eqref{eq:BCH}}}{=} e^{C(v') + \alpha c } e^{C(w')+\beta c} 
\underset{\mbox{\tiny eq.\eqref{eq:e_C}}}{=} e^C(e^{v'+\alpha c}) e^C(e^{w'+\beta c}).
\end{align*}
Hence we have
\begin{align*}
e^C(\exp(\ssum_{j=1}^r t_j' v_j')\cdot h')= \exp(\ssum_{j=1}^r t_j' C(v_j')) \cdot e^C(h').
\end{align*}
Since $e^C(h') \in H$, from eq.\eqref{eq:I_ell_B_formula} for ${\bf I}_{\ell,\mathcal{B}}$ we get
\begin{align*}
({\bf I}_{\ell,\mathcal{B}}^{-1} \varphi)( e^{C}(\exp(\ssum_{j=1}^r t_j' v_j')\cdot h')) & = \chi(e^C(h'))^{-1} \cdot \varphi(t_1',\ldots,t_r') = \chi'(h')^{-1} \cdot \varphi(t_1',\ldots,t_r').
\end{align*}
Hence, we showed that ${\bf I}_{\ell',\mathcal{B}'}^{-1} {\bf R} \varphi = \til{\bf R} {\bf I}_{\ell,\mathcal{B}}^{-1} \varphi$ holds for all $\varphi \in \mathscr{H}_{\ell,\mathcal{B}}$, as desired. \qed

\vs

Using Lem.\ref{lem:til_R}, we see that ${\bf R}$ serves as a pullback map for the Schr\"odinger representations of $N$ and $\mathfrak{n}$.
\begin{corollary}[a pullback of representations]
\label{cor:R_is_pullback}
The unitary operator ${\bf R} = {\bf R}_{C(\ell',\mathcal{B}'), (\ell',\mathcal{B}')}$ satisfies
\begin{align}
\label{eq:conjugation_by_R}
{\bf R} \, \pi_{C(\ell',\mathcal{B}')} (C(x')) \, {\bf R}^{-1} = \pi_{\ell',\mathcal{B}'}(x'), \qquad \forall x' \in \mathfrak{n}',
\end{align}
which hold as equalities of skew self-adjoint operators.
\end{corollary}

{\it Proof.} It suffices to prove the exponentiated version, namely
$$
{\bf R} \, \pi_{C(\ell',\mathcal{B}')}(e^{C(x')}) \, {\bf R}^{-1} = \pi_{\ell',\mathcal{B}'}(e^{x'}), \qquad \forall x' \in \mathfrak{n}',
$$
which are equations of unitary operators. In view of eq.\eqref{eq:til_R_diagram}, it suffices to prove
$$
\til{\bf R} \, \pi_{C(\ell')}(e^{C(x')}) \, \til{\bf R}^{-1} = \pi_{\ell'}(e^{x'}), \qquad \forall x' \in \mathfrak{n}'.
$$
From eq.\eqref{eq:e_C}, we have $e^{C(x')} = e^C(e^{x'}) \in N$. Let's now show the equality $\til{\bf R} \, \pi_{C(\ell')}(e^C(e^{x'})) = \pi_{\ell'}(e^{x'}) \, \til{\bf R}$ of maps $\mathscr{H}_{C(\ell')} \to \mathscr{H}_{\ell'}$. From the definitions of $\til{\bf R}$, $\pi_{C(\ell')}$, and $\pi_{\ell'}$, one observes that, for each $k \in \mathscr{H}_{C(\ell')} = \mathscr{H}_\ell$ and $n' \in N'$, the following holds.
\begin{align*}
( \til{\bf R} \, (\pi_{C(\ell')}(e^C(e^{x'})) \, k) )(n')
& = (\pi_{C(\ell')}(e^C(e^{x'})) \, k) ( e^C(n')) \qquad (\because \mbox{eq.\eqref{eq:til_R_formula}}) \\
& = k( (e^C(e^{x'}) )^{-1} \, e^C(n') ) \qquad\qquad (\because \mbox{eq.\eqref{eq:pi_ell_action}}) \\
& = k( e^C( (e^{x'})^{-1} \, n') ) \qquad\qquad\quad (\because \mbox{$e^C$ is a homomorphism}) \\
& = (\til{\bf R} \, k)( (e^{x'})^{-1} n') \qquad\qquad\quad (\because \mbox{eq.\eqref{eq:til_R_formula}}) \\
& = ( \pi_{\ell'}(e^{x'}) (\til{\bf R} \, k) )(n') \qquad\qquad (\because \mbox{eq.\eqref{eq:pi_ell_action}})  \quad \qed
\end{align*}

\vs

Back to our situation of $C\in {\rm Aut}(V,B)$. Now, $C : V\to V$ is an invertible linear transformation respecting $B$. We chose one symplectic decomposition $(\ell,\mathcal{B})$ of $(V,B)$, considered another symplectic decomposition $C(\ell,\mathcal{B}) = (C(\ell),C(\mathcal{B}))$ of $(V,B)$, and the Weil intertwiner ${\bf F} : \mathscr{H}_{C(\ell,\mathcal{B})} \to \mathscr{H}_{\ell,\mathcal{B}}$. Now we identify the domain $\mathscr{H}_{C(\ell,\mathcal{B})}$ with the codomain $\mathscr{H}_{\ell,\mathcal{B}}$ via the pullback map ${\bf R} : \mathscr{H}_{C(\ell,\mathcal{B})} \to \mathscr{H}_{\ell,\mathcal{B}}$ in Def.\ref{def:R}. That is, we define $\rho(C)$ as the unique map making the following diagram to commute:
$$
\xymatrix@R+0mm@C+3mm{
\mathscr{H}_{C(\ell,\mathcal{B})} \ar[r]^-{{\bf F}} \ar[d]_{{\bf R}} & \mathscr{H}_{\ell,\mathcal{B}} \\ 
\mathscr{H}_{\ell,\mathcal{B}} 
\ar[ur]_{\rho(C)} &
}
$$
We can reverse the vertical arrow ${\bf R}$ and still call it ${\bf R}$, because ${\bf R} :\mathscr{H}_{C(\ell,\mathcal{B})} \to \mathscr{H}_{\ell,\mathcal{B}}$ (for $C$) and ${\bf R} : \mathscr{H}_{\ell,\mathcal{B}} \to \mathscr{H}_{C(\ell,\mathcal{B})}$ (for $C^{-1}$) are inverses to each other. So the definition of $\rho(C)$ can be written as
\begin{align}
\label{eq:rho_C}
\rho(C) = \rho_{\ell,\mathcal{B}}(C) := 
{\bf F}_{C(\ell,\mathcal{B}),(\ell,\mathcal{B})} \circ {\bf R}_{(\ell,\mathcal{B}),C(\ell,\mathcal{B})}
~:~ \mathscr{H}_{\ell,\mathcal{B}} \longrightarrow \mathscr{H}_{\ell,\mathcal{B}}, \quad \forall C \in {\rm Aut}(V,B).
\end{align}
See e.g. eq.(1.6.9) of \cite{LV}. We need to verify that $\rho(C_1) \rho(C_2) = \rho(C_1C_2)$ holds up to constant for all $C_1,C_2\in {\rm Aut}(V,B)$. For this, we first establish the following lemma, which is written more generally than needed at the moment.  
\begin{lemma}[the compatibility between the Weil intertwiners ${\bf F}$ and the pullback maps ${\bf R}$]
\label{lem:compatibility_between_F_and_R}
Let $(V,B)$, $(V',B')$, $C:V' \to V$, $f$, and $f'$ be as in Def.\ref{def:R}. Let $(\ell_1',\mathcal{B}_1')$ and $(\ell_2',\mathcal{B}_2')$ be any symplectic decompositions of $(V',B')$; so, $C(\ell_1',\mathcal{B}_1')$ and $C(\ell_2',\mathcal{B}_2')$ are symplectic decompositions of $(V,B)$. Then the following diagram of the Weil intertwiners and the pullback maps commutes:
$$
\xymatrix@C+2mm@R+0mm{
\mathscr{H}_{C(\ell_1',\mathcal{B}_1')} \ar[r]^-{{\bf F}} \ar[d]_{{\bf R}} & \mathscr{H}_{C(\ell_2',\mathcal{B}_2')} \ar[d]^{{\bf R}} \\
\mathscr{H}_{\ell_1',\mathcal{B}_1'} \ar[r]^-{{\bf F}} & \mathscr{H}_{\ell_2',\mathcal{B}_2'}
}
$$
That is,
$$
{\bf R}_{C(\ell_2',\mathcal{B}_2'),(\ell_2',\mathcal{B}_2')} \, {\bf F}_{C(\ell_1',\mathcal{B}_1'),C(\ell_2',\mathcal{B}_2')} 
= {\bf F}_{(\ell_1',\mathcal{B}_1'),(\ell_2',\mathcal{B}_2')} \, {\bf R}_{C(\ell_1',\mathcal{B}_1'),(\ell_1',\mathcal{B}_1')}.
$$
\end{lemma}
This lemma is asserted in \S1.6.8 of \cite{LV} for the case when $(V,B)=(V',B')$ and $V^\perp=0$ and used in the construction and computation of the symplectic group representation, but not explicitly proved there. Since it will be used crucially throughout the present paper, we present a complete proof.

\vs

{\it Proof of Lem.\ref{lem:compatibility_between_F_and_R}.} In view of the diagrams in eq.\eqref{eq:Weil_intertwiner_for_ell_B} and eq.\eqref{eq:til_R_diagram}, it suffices to show the commutativity of the following diagram
$$
\xymatrix{
\mathscr{H}_{C(\ell_1')} \ar[r]^{{\bf F}} \ar[d]_{\til{\bf R}} & \mathscr{H}_{C(\ell_2')} \ar[d]^{\til{\bf R}} \\
\mathscr{H}_{\ell_1'} \ar[r]^{{\bf F}} & \mathscr{H}_{\ell_2'}.
}
$$
Let $k \in \mathscr{H}_{C(\ell_1')}$. By eq.\eqref{eq:til_R_formula}, $\til{\bf R} k \in \mathscr{H}_{\ell_1'}$ is given by
$$
(\til{\bf R}k)(n') = k(e^C n'), \qquad \forall n \in N',
$$
where $e^C:N' \to N$ is defined as in Lem.\ref{lem:til_R}, i.e. by eq.\eqref{eq:e_C}. By the definition of ${\bf F}$ as written in eq.\eqref{eq:canonical_intertwiner_formula}, ${\bf F} \til{\bf R} k \in \mathscr{H}_{\ell_2'}$ is given by
\begin{align*}
({\bf F} \til{\bf R} k)(n') & = \int_{H_2'/H_1'\cap H_2'} (\til{\bf R} k)(n' h') \, \chi_2' (h') \, dm_{H_2'/H_1'\cap H_2'}(h'), \\
& = \int_{H_2'/H_1'\cap H_2'} k(e^C(n' h')) \, \chi_2' (h') \, dm_{H_2'/H_1'\cap H_2'}(h'), \qquad \forall n' \in N',
\end{align*}
where $\chi_2' : H_2' \to {\rm U}(1)$ is the character for $H_2'$, and the measure $dm_{H_2'/H_1'\cap H_2'}$ is uniquely normalized so that the formula in eq.\eqref{eq:canonical_intertwiner_formula} gives our final normalized unitary map ${\bf F} : \mathscr{H}_{\ell_1'} \to \mathscr{H}_{\ell_2'}$ on the nose.

\vs

Meanwhile, for the essential Lagrangians $C(\ell_1')$ and $C(\ell_2')$ of $(V,B)$, the corresponding polarizations of $\mathfrak{n}$ are $C(\mathfrak{h}_1')$ and $C(\mathfrak{h}_2')$, as seen in the proof of Lem.\ref{lem:til_R}; so their exponentiated subgroups are $e^C(H_1')$ and $e^C(H_2')$. Hence, in view of eq.\eqref{eq:canonical_intertwiner_formula}, ${\bf F} k \in \mathscr{H}_{C(\ell_2')}$ is given by
\begin{align*}
({\bf F} k)(n) = \int_{e^C(H_2')/e^C(H_1')\cap e^C(H_2')} k(n h) \, \chi_2 (h) \, dm_{e^C(H_2')/e^C(H_1')\cap e^C(H_2')}(h), \quad \forall n \in N,
\end{align*}
where $\chi_2 : e^C(H_2') \to {\rm U}(1)$ is the character for $C(H_2')$, and the measure $dm_{e^C(H_2')/e^C(H_1')\cap e^C(H_2')}$ is suitably uniquely normalized. As seen in the proof of Lem.\ref{lem:til_R}, we have $\chi_2' = \chi_2 \circ e^C$. So we have, for each $n' \in N'$, that
\begin{align*}
({\bf F} \til{\bf R} k)(n') & = \int_{H_2'/H_1'\cap H_2'} k(e^C(n) e^C(h')) \, \chi_2 (e^C(h')) \, dm_{H_2/H_1\cap H_2}(h'), \\
(\til{\bf R} {\bf F} k)(n') & = \int_{e^C(H_2')/e^C(H_1')\cap e^C(H_2')} k(e^C(n') h) \, \chi_2 (h) \, dm_{e^C(H_2)/e^C(H_1)\cap e^C(H_2)}(h).
\end{align*}
As $e^C$ gives a natural identification from $H_2' / H_1' \cap H_2'$ to $e^C(H_2') / e^C(H_1') \cap e^C(H_2')$, we can observe that the above two integrals coincide via the change of variables $h = e^C(h')$, maybe up to a positive real constant, due to the choice of the normalizations of the measures. So ${\bf F} \til{\bf R} = \til{\bf R} {\bf F}$ holds up to a positive real constant; since ${\bf F}$'s and $\til{\bf R}$'s are unitary, it follows that this positive constant must be $1$. \qed

\vs

Now we prove that $C \mapsto \rho(C) = \rho_{\ell,\mathcal{B}}(C)$ is a projective unitary representation of ${\rm Aut}(V,B)$ on $\mathscr{H}_{\ell,\mathcal{B}}$, and that the projective phase constants are given by the Maslov indices.
\begin{proposition}
\label{prop:projective_representation_of_Aut_V_B}
For any $C_1,C_2 \in {\rm Aut}(V,B)$, the assignment $\rho = \rho_{\ell,\mathcal{B}}$ satisfies
$$
\rho(C_1) \, \rho(C_2) = e^{-\frac{{\rm i} \pi}{4} \tau(\ell, C_1(\ell),C_1(C_2(\ell)))} \, \rho(C_1C_2),
$$
where $\tau(\cdot,\cdot,\cdot)$ is as defined in Def.\ref{def:Kashiwara_index}--\ref{def:Kashiwara_index2}.
\end{proposition}
See \cite[Thm.1.6.11]{LV} for the cases $V^\perp=0$.

\vs

{\it Proof.} Consider the diagram
$$
\xymatrix@C-1mm{
& & \mathscr{H}_{C_1(C_2(\ell,\mathcal{B}))} \ar[dr]^{{\bf F}} \ar@/^4.0pc/[ddrr]^-{{\bf F}} & & \\
& \mathscr{H}_{C_2(\ell,\mathcal{B})} \ar[dr]^{{\bf F}} \ar[ur]^{{\bf R}} & \mbox{square} & \mathscr{H}_{C_1(\ell,\mathcal{B})} \ar[dr]^{{\bf F}} & \\
\mathscr{H}_{\ell,\mathcal{B}} \ar[rr]_-{\rho_{\ell,\mathcal{B}}(C_2)} \ar[ur]^-{{\bf R}} \ar@/^4.0pc/[uurr]^{{\bf R}} & & \mathscr{H}_{\ell,\mathcal{B}} \ar[rr]_-{\rho_{\ell,\mathcal{B}}(C_1)} \ar[ur]^{{\bf R}} & & \mathscr{H}_{\ell,\mathcal{B}} \\
}
$$
consisting of one square and four triangles. The bottom two triangles are commutative by definition, and the square is commutative by Lem.\ref{lem:compatibility_between_F_and_R} with $C_1^{-1}$ playing the role of $C$. The upper-left triangle is commutative by Lem.\ref{lem:R_composition_identity} (or eq.\eqref{eq:R_composition_identity}). The upper-right triangle is commutative up to a constant; by Prop.\ref{prop:Weil_intertwiner_involutivity} and Prop.\ref{prop:F_composition_identity}, we see
$$
{\bf F}_{C_1(\ell,\mathcal{B}),(\ell,\mathcal{B})} \circ {\bf F}_{C_1(C_2(\ell,\mathcal{B})),C_1(\ell,\mathcal{B})} = e^{\frac{{\rm i}\pi}{4} \tau(\ell,C_1(C_2(\ell)),C_1(\ell))} \, {\bf F}_{C_1(C_2(\ell,\mathcal{B})),(\ell,\mathcal{B})}
$$
Finally, note that the composition of the two curved arrows, ${\bf R}$ and ${\bf F}$, equals $\rho_{\ell,\mathcal{B}}(C_1C_2)$ by the definition of $\rho_{\ell,\mathcal{B}}$, and use an easy observation $\tau(\ell_1,\ell_2,\ell_3) = - \tau(\ell_1,\ell_3,\ell_2)$. \qed

\section{The irreducible representations of the quantum Teichm\"uller space}

\subsection{The irreducible quantum representations, associated to an ideal triangulation $T$}
\label{subsec:irreducible_quantum_representation}

Let $S$ be a compact oriented surface of genus $g$ minus a finite set $\mathcal{P}$ of \ul{\bf punctures} of size $|\mathcal{P}|=s$. Let $T$ be an \ul{\bf ideal triangulation} of $S$, i.e. a `triangulation of $S$ with vertices in $\mathcal{P}$'. That is, $T$ is a collection of mutually non-intersecting non-homotopic simple paths running between punctures, called the \ul{\bf edges} of $T$, that divide $S$ into regions bounded by three edges, called the \ul{\bf triangles}; see \cite[Def.2.6]{FST} for a precise definition. We consider $T$ up to isotopy that respects the above conditions. It is easy to see that $T$ has $6g-6+3s$ edges.  Define a real vector space $V_T$ freely generated by a formal set $\{ x_e : e\in T\}$ enumerated by $T$, equipped with the skew-symmetric $\mathbb{R}$-bilinear form $B_T$ given on the basis vectors by
\begin{align}
\label{eq:B_T}
& \hspace{55mm} B_T(x_e, x_f) = \varepsilon_{ef}, \qquad \mbox{where} \\
\nonumber
& \begin{array}{l}
\varepsilon_{ef} = a_{ef} - a_{fe}, \\
a_{ef} = \mbox{the number of corners of triangles of $T$ delimited by $e$ from the left and by $f$ from the right}.
\end{array}
\end{align}
The matrix $\varepsilon_T = \varepsilon = (\varepsilon_{ef})_{e,f\in T}$ is called the \ul{\bf exchange matrix} for $T$. We now apply the results of \S\ref{sec:Schrodinger_representations_and_Weil_intertwiners} to $(V_T,B_T)$. Denote the corresponding Heisenberg algebra by $\mathfrak{n}_T$, and the Heisenberg group by $N_T$.

\begin{remark}
The notion of the left and the right appearing in the definition of $a_{ef}$ can be consistently defined, using the orientation of $S$. See e.g. \cite[\S4]{FST} for a more precise treatment of these numbers and $\varepsilon_{ef}$. 
\end{remark}

For each puncture $p\in \mathcal{P}$ and an edge $e\in T$, denote by $\sigma^e_p$ be the number of incidences of $e$ at $p$. That is, if $p$ is not an endpoint of $e$, then $\sigma^e_p=0$. If $e$ has two distinct endpoints and one of them is $p$, then $\sigma^e_p=1$. If the two endpoints $e$ both coincide with $p$, then $\sigma^e_p=2$.
\begin{lemma}[\cite{BL}; the puncture elements]
\label{lem:puncture_elements}
The radical $V_T^\perp$ of $(V_T,B_T)$ has a basis $\{x_p : p \in \mathcal{P}\}$, where
$$
x_p := {\textstyle \underset{e\in T}{\sum}} \sigma^e_p x_e, \quad \forall p \in \mathcal{P}.
$$
\end{lemma}

{\it Proof.} Use \cite[Lem.8]{BL} and an easy observation $\sum_{e\in T} x_e = 2 \sum_{p \in \mathcal{P}} x_p$ (also in \cite{BL}). \qed

\vs

For each $T$, we shall use this as a preferred basis of $V^\perp_T$. Choose and fix any function
\begin{align}
\label{eq:lambda}
\lambda : \mathcal{P} \to \mathbb{R}, \quad p \mapsto \lambda(p).
\end{align}
Note that this choice of $\lambda$ does not involve $T$. All representations from now on will depend on this fixed choice of $\lambda$; that is, the function 
$$
f = f_T = f_{T;\lambda} :V_T^\perp \to \mathbb{R}
$$
which is a basic ingredient of the construction in \S\ref{subsec:Schrodinger_representation_for_all} is set as
$$
f_T(x_p) = \lambda(p), \qquad \forall p\in \mathcal{P}.
$$

\vs

Now, choose any symplectic decomposition $(\ell_T,\mathcal{B}_T)$ of  $(V_T,B_T)$; when clear, we shall sometimes say that $(\ell_T,\mathcal{B}_T)$ is a symplectic decomposition \ul{\bf for} the triangulation $T$. Then we have a representation $\pi_{\ell_T}$ on the Hilbert space $\mathscr{H}_{\ell_T}$ of the Heisenberg algebra $\mathfrak{n}_T$ (and of the Heisenberg group $N_T$ too), and also a representation $\pi_{\ell_T,\mathcal{B}_T}$ on the Hilbert space (see eq.\eqref{eq:H_ell_B}).
$$
\mathscr{H}_{\ell_T,\mathcal{B}_T} = L^2(\mathbb{R}^{3g-3+s}, dt_1 \cdots t_{3g-3+s})
$$ 
with the identification map (see eq.\eqref{eq:I_ell_B_formula})
$$
{\bf I}_{\ell_T,\mathcal{B}_T} : \mathscr{H}_{\ell_T} \to \mathscr{H}_{\ell_T,\mathcal{B}_T}.
$$
In particular, since $\pi_\ell$ and $\pi_{\ell,\mathcal{B}}$ are Lie algebra representations,  the skew self-adoint operators $\pi_\ell(x_e)$ on $\mathscr{H}_{\ell_T}$ and $\pi_{\ell_T,\mathcal{B}_T}(x_e)$ on $\mathscr{H}_{\ell_T,\mathcal{B}_T}$ associated to the vectors $x_e \in V_T \subset \mathfrak{n}_T$ satisfy
$$
[ \pi_\ell(x_e), \, \pi_\ell(x_f) ] = - {\rm i} \, \varepsilon_{ef} \cdot {\rm Id}, \qquad \forall e,f \in T
$$
(see Example \ref{ex:operator_in_center}), which makes sense e.g. on some dense subspace, and likewise for $\pi_{\ell_T,\mathcal{B}_T}$. 

\vs

Certain normalization would be suitable for us here. Note that the element $x_e$ of $V_T$ is to symbolize Thurston's shear coordinate function on the enhanced Teichm\"uller space of $S$, associated to the edge $e$ of the ideal triangulation $T$. In particular, it is a real valued function, and hence its quantization $\wh{x}_e$ should be a self-adjoint operator, and the usual requirement adopted in the literature is the following commutation relations
\begin{align}
\label{eq:usual_commutation_of_hat_x_e}
[\wh{x}_e,\wh{x}_f] = 2\pi {\rm i} \hbar \, \varepsilon_{ef} \cdot {\rm Id}, \qquad \forall e,f\in T,
\end{align}
where $\hbar \in \mathbb{R}$ is the real parameter for the quantization (the Planck constant). From now on, we fix this parameter
$$
\hbar>0.
$$
\begin{definition}[the re-scaled self-adjoint representations of the Heisenberg algebra]
We re-scale the representations $\pi_{\ell_T}$ and $\pi_{\ell_T,\mathcal{B}_T}$ of the Heisenberg algebra $\mathfrak{n}_T$ by ${\rm i} \sqrt{2\pi \hbar}$, to define the representations $\pi_{\ell_T}^\hbar$ and $\pi_{\ell_T, \mathcal{B}_T}^\hbar$ of $\mathfrak{n}_T$: 
$$
\pi^\hbar_{\ell_T}(x) = {\rm i} \sqrt{2\pi \hbar} \cdot \pi_{\ell_T}(x), \qquad
\pi^\hbar_{\ell_T,\mathcal{B}_T}(x) = {\rm i} \sqrt{2\pi \hbar} \cdot \pi_{\ell_T,\mathcal{B}_T}(x), \qquad \forall x\in \mathfrak{n}_T.
$$
Now, for each $e\in T$ we denote
$$
\wh{x}^\hbar_{e;\ell_T} = \pi^\hbar_{\ell_T}(x_e) \quad \mbox{and}\quad \wh{x}^\hbar_{e;\ell_T,\mathcal{B}_T} = \pi^\hbar_{\ell_T,\mathcal{B}_T}(x_e), \qquad \forall e\in T,
$$
and denote either of them by
$$
\wh{x}_e^\hbar = \wh{x}^\hbar_{e;T}
$$
by abuse of notation.
\end{definition}
Then, for each $x\in \mathfrak{n}_T$, the operators $\pi^\hbar_{\ell_T}(x)$ and $\pi^\hbar_{\ell_T,\mathcal{B}_T}$ are self-adjoint. And the operators $\wh{x}^\hbar_e$ defined this way are self-adjoint operators and satisfy the sought-for usual commutation relations in eq.\eqref{eq:usual_commutation_of_hat_x_e}. 

\vs

We thus obtained a family of {\em irreducible} self-adjoint representations of the Heisenberg algebra $\mathfrak{n}_T$ generated by the family of quantum shear coordinate functions for an ideal triangulation $T$. Such a family of irreducible representations has not been dealt with in the quantum Teichm\"uller theory literature so concretely, systematically, and generally as we just did here. Its existence had been widely accepted and mentioned already in Chekhov-Fock \cite{CF}, but not really constructed explicitly. In order to exactly compute the multiplicative constants coming from the composition of intertwiners between these representations, as we are trying to accomplish in the present paper, one must have an explicit construction at hand. We also note that Fock-Goncharov \cite{FG09} considered a more canonical but not irreducible family of representations, which are irreducible only as representations of a bigger algebra, namely the symplectic double.

\subsection{The quantum mutation: the algebraic quantum coordinate change for a flip along an edge}
\label{subsec:quantum_mutation_algebraic}

Suppose now that two ideal triangulations $T$ and $T'$ are related by the \ul{\bf flip at the edge $k$}, i.e. they differ only at one edge, namely $k$; the edges of $T$ are in natural bijection with those of $T'$, hence we identify these two as sets. Choose symplectic decompositions $(\ell_T,\mathcal{B}_T)$ and $(\ell_{T'},\mathcal{B}_{T'})$ of $(V_T,B_T)$ and $(V_{T'},B_{T'})$ respectively. Our goal in the present subsection is to construct an intertwiner
\begin{align}
\label{eq:K_hbar_k}
{\bf K}^\hbar_k = {\bf K}^\hbar_{TT'} = {\bf K}^\hbar_{\ell_T, \ell_{T'}} : \mathscr{H}_{\ell_{T'}} \to \mathscr{H}_{\ell_T},
\end{align}
i.e. a unitary operator that intertwines the action $\pi^\hbar_{\ell_{T'}}$ of $\mathfrak{n}_{T'}$ and the action $\pi^\hbar_{\ell_T}$ of $\mathfrak{n}_T$. To distinguish from the Weil intertwiner, we call this sought-for operator ${\bf K}_k^\hbar$ a \ul{\bf mutation intertwiner}, for a flip is an example of the notion of the {\em mutation} in the theory of cluster varieties.

\vs

To even just formulate this problem properly, we first need to describe the relationship between the two Lie algebras $\mathfrak{n}_T$ and $\mathfrak{n}_{T'}$. In this subsection, we briefly review the setting from the literature. At the classical level, upon the flip along the edge $k$, the shear coordinate functions for the edges of $T$ on the enhanced Teichm\"uller space are related to the shear coordinate functions for the edges of $T'$ by certain formulas: if we denote by $x_e$ the shear coordinate for the edge $e$ of $T$ and its exponential by $X_e$
$$
X_e = \exp(x_e),
$$
then the coordinate change formula is given by
\begin{align}
\label{eq:mutation_formula_for_X}
X_e' = \left\{
\begin{array}{ll}
X_k^{-1}, & \mbox{if $e=k$,} \\
X_e(1+X_k^{{\rm sgn}(-\varepsilon_{ek})})^{-\varepsilon_{ek}}, & \mbox{if $e\neq k$}.
\end{array}
\right.
\end{align}

\begin{remark}
\label{rem:mutation_only_for_regular_flips}
The transformation formula in eq.\eqref{eq:mutation_formula_for_X} works for flips when the triangles adjacent to the flipped edge $k$ are not self-folded; in the present paper, we restrict our attention to these flips, for convenience. When a self-folded triangle is involved, one must be careful; see e.g. \cite{FST} \cite{BL} \cite{Liu}.
\end{remark}

For quantum version, for each $T$ we must first construct a non-commutative algebra $\mathcal{X}^q_T$, which can be regarded as a generalized quantum torus algebra.
\begin{definition}[\cite{BL} \cite{Liu}]
\label{def:CF_algebra}
For an ideal triangle $T$, the \ul{\bf Chekhov-Fock algebra} $\mathcal{X}^q_T$ is  the associative $\mathbb{C}$-algebra generated by $\wh{X}_e$'s (where $e$ ranges in $T$) and their inverses, mod out by the relations
$$
\wh{X}_e \wh{X}_f = q^{2\varepsilon_{ef}} \wh{X}_f \wh{X}_e, \quad \forall e,f \in T.
$$
Here, $q$ may be regarded either as a complex number of modulus $1$ defined as
$$
q = e^{{\rm i} \pi \hbar},
$$
or as a formal symbol so that $\mathcal{X}^q_T$ is a $\mathbb{C}[q,q^{-1}]$-algebra defined in terms of generators and relations as above.
\end{definition}
This algebra $\mathcal{X}^q_T$ satisfies the so-called Ore condition, hence its skew-field of fractions ${\rm Frac}(\mathcal{X}^q_T)$ can be considered. It is in fact a $*$-algebra, with the $*$-structure given on the generators by
$$
*\wh{X}_e = \wh{X}_e, \quad \forall e\in T.
$$

\vs

For a flip $T\leadsto T'$ along the edge $k$ one must construct a \ul{\bf quantum coordinate change map}
$$
\Phi^q_{TT'} : {\rm Frac}(\mathcal{X}^q_{T'}) \to {\rm Frac}(\mathcal{X}^q_T)
$$
that is an isomorphism recovering the classical coordinate change map as $q\to 1$, and that satisfies the consistency relations corresponding to the ones satisfied by their classical counterparts. The consistency relations guarantee that one can construct a well-defined quantum coordinate change map between any two ideal triangulations, not just for flips. A first version of such a map $\Phi^q_{TT'}$ was found by Chekhov-Fock \cite{CF} and also by Kashaev \cite{Kash} in a slightly different setting, then established by Bonahon and Liu \cite{BL} \cite{Liu} for all possible flips, and later generalized to cluster varieties by Fock-Goncharov \cite{FG09} and developed by others, including Kashaev-Nakanishi \cite{KN}.  As an explicit formula for this map $\Phi^q_{TT'}$ shall follow from our sought-for intertwiner which can be viewed as an operator version of $\Phi^q_{TT'}$ and which we will construct, we do not describe it here at the moment; a reader can consult e.g. \cite{CF} \cite{BL} \cite{Liu} \cite{FG09} \cite{K16c} \cite{KN} for a detailed algebraic treatment.

\vs

For each $T$, an element of ${\rm Frac}(\mathcal{X}^q_T)$ is said to be \ul{\bf Laurent for $T$} if it belongs to $\mathcal{X}^q_T$. As we shall soon see, we will construct a family of representations of $\mathcal{X}^q_T$, but not of ${\rm Frac}(\mathcal{X}^q_T)$, hence the intertwining equations for a flip $T\leadsto T'$ can be asked only for  those elements that are Laurent for both $T$ and $T'$. In the end, we shall ask the intertwining equations for the elements that are Laurent for {\em every} ideal triangulation $T$, i.e. the \ul{\bf universally Laurent} elements \cite{FG09}; define this algebra as
$$
{\bf L}^q_T := \bigcap_{T'} \Phi^q_{TT'} (\mathcal{X}^q_{T'}) \quad \subset \quad  \mathcal{X}^q_T \quad \subset \quad {\rm Frac}(\mathcal{X}^q_T)
$$
where the intersection is taken over all possible ideal triangulations $T'$; then one can observe that ${\bf L}^q_T$ for different $T$'s are all canonically identified via the (restrictions of the) maps $\Phi^q_{TT'}$:
$$
\Phi^q_{TT'} : {\bf L}^q_{T'} \overset{\cong}{\longrightarrow} {\bf L}^q_T
$$

\subsection{The mutation intertwiner: the formulation of the problem}
\label{subsec:mutation_intertwiner_formulation_of_the_problem}

Coming back to the intertwiner problem, we should consider a family of representations on a Hilbert space of each Chekhov-Fock algebra $\mathcal{X}^q_T$, or really, of the universally Laurent $*$-subalgebra ${\bf L}^q_T$. We use the representations $(\mathscr{H}_{\ell_T}, \pi^\hbar_{\ell_T})$ and $(\mathscr{H}_{\ell_T,\mathcal{B}_T}, \pi^\hbar_{\ell_T,\mathcal{B}_T})$ of the Heisenberg algebra $\mathfrak{n}_T$. Recall the self-adjoint operator $\wh{x}^\hbar_{e;T} = \wh{x}^\hbar_{e;\ell_T}$ or $\wh{x}^\hbar_{e;\ell_T,\mathcal{B}_T}$ defined in \S\ref{subsec:irreducible_quantum_representation}. Then to each generator $\wh{X}_e^{\pm 1}$ of $\mathcal{X}^q_T$ we associate the self-adjoint operator $\exp(\pm \wh{x}^\hbar_{e;T})$, defined via the functional calculus for $\wh{x}^\hbar_{e;T}$. This operator for $\wh{X}_e^{\pm 1}$ is restricted to a certain nice dense subspace $D_{\ell_T}$ of $\mathscr{H}_{\ell_T}$, or the corresponding subspace $D_{\ell_T,\mathcal{B}_T}$ of $\mathscr{H}_{\ell_T,\mathcal{B}_T} = L^2(\mathbb{R}^r, dt_1 \cdots dt_r)$, defined as (\cite{K16c} \cite{FG09})
$$
D_{\ell_T,\mathcal{B}_T} := {\rm span}_\mathbb{C} \left\{ 
e^{\vec{t} M \vec{t}^{\,{\rm t}} + \vec{a} \cdot \vec{t}} \, P(\vec{t}) \left| 
\begin{array}{l}
\mbox{$M\in {\rm Mat}_{r\times r}(\mathbb{C})$ with a negative-definite real part ${\rm Re}(M)$,} \\
\mbox{$\vec{a} = (a_1\cdots a_r) \in \mathbb{C}^n$, and $P$ a polynomial in $t_1,\ldots t_r$ over $\mathbb{C}$}
\end{array}
\right.
\right\}.
$$
In particular, the operator for $\wh{X}_e^{\pm 1}$ preserves this space. Then to each element $u$ of ${\bf L}^q_T$, which is a polynomial in $\wh{X}_e^{\pm 1}$, we associate the corresponding polynomial in $\exp(\pm \wh{x}^\hbar_{e;T})$, denoted by $\wh{u}$, which makes sense as an operator on $D_{\ell_T,\mathcal{B}_T}$ to itself. As it turns out, this subspace $D_{\ell_T,\mathcal{B}_T}$ is too small for the intertwining problem to make sense, and we need to consider the common maximal domain of all $\wh{u}$'s, in the following way:
$$
\mathscr{S}_{\ell_T,\mathcal{B}_T}^q := \bigcap_{u \in {\bf L}^q_T} {\rm Dom}(\wh{u}^*) \quad \subset \quad \mathscr{H}_{\ell_T,\mathcal{B}_T},
$$
where $\wh{u}^*$ is the operator adjoint of $\wh{u}$, namely ${\rm Dom}(\wh{u}^*)$ is the set of all $\xi \in \mathscr{H}_{\ell_T,\mathcal{B}_T}$ such that the map $D_{\ell_T,\mathcal{B}_T} \to \mathbb{C}$, $w \mapsto \langle \wh{u}(w), \xi\rangle$, is a continuous functional, where $\langle \cdot,\cdot\rangle$ denotes the Hilbert space inner product on $\mathscr{H}_{\ell_T,\mathcal{B}_T}$. This space $\mathscr{S}^q_{\ell_T,\mathcal{B}_T}$ is a direct analog of what Fock-Goncharov refers to as the \ul{\bf Schwartz space} in \cite{FG09}. Indeed, it generalizes the role of the classical Schwartz space; so, we might call it the \ul{\bf Fock-Goncharov Schwartz space}. Finally, on this space, we can define the representation $\pi^q_{\ell_T,\mathcal{B}_T}$ of the algebra ${\bf L}^q_T$ as
$$
\pi^q_{\ell_T,\mathcal{B}_T}(u) := \wh{*u}^* \restriction \mathscr{S}^q_{\ell_T,\mathcal{B}_T}, \quad \forall u \in {\bf L}^q_T.
$$
Namely, consider the image $*u \in {\bf L}^q_T$ of $u$ under the $*$-map, the associated operator $\wh{*u}$ on the nice subspace $D_{\ell_T,\mathcal{B}_T}$, take its operator adjoint, and then restrict down to the Fock-Goncharov Schwartz space $\mathscr{S}^q_{\ell_T,\mathcal{B}_T}$. Analogously to the classical Schwartz space, this Fock-Goncharov Schwartz space $\mathscr{S}^q_{\ell_T,\mathcal{B}_T}$ is then given the Frech\'et topology defined by the family of semi-norms $\rho_u$, defined for each $u\in {\bf L}^q_T$ as $\rho_u(\xi) := || \pi^q_{\ell_T,\mathcal{B}_T}(u) \xi ||$, where $||\cdot ||$ is the Hilbert space norm. It can be shown that each operator $\pi^q_{\ell_T,\mathcal{B}_T}(u)$ preserves $\mathscr{S}^q_{\ell_T,\mathcal{B}_T}$. We note that the Fock-Goncharov Schwartz space considered here is not precisely same as the one in \cite{FG09}; for example, in \cite{FG09}, the algebra ${\bf L}^q_T$ is replaced by a bigger algebra which incorporates the concepts of the `symplectic double' and the `modular double'. 

\vs

We can finally rigorously formulate the intertwiner problem for each flip $T \leadsto T'$. Choose symplectic decompositions $(\ell_T,\mathcal{B}_T)$and $(\ell_{T'},\mathcal{B}_{T'})$ for $T$ and $T'$, as before. We would like to find a unitary map ${\bf K}^\hbar_k$ as in eq.\eqref{eq:K_hbar_k} that induces a bijection between the Schwartz spaces $\mathscr{S}^q_{\ell_{T'},\mathcal{B}_{T'}}$ and $\mathscr{S}^q_{\ell_T, \mathcal{B}_T}$, preferably a homeomorphism with respect to the Frech\'et topologies, that satisfies the \ul{\bf intertwining equations}
\begin{align}
\label{eq:intertwining_equations}
{\bf K}^\hbar_k \circ \pi^q_{\ell_{T'},\mathcal{B}_{T'}}(u) = \pi^q_{\ell_T,\mathcal{B}_T}( \Phi^q_{TT'}(u)) \circ {\bf K}^\hbar_k ~ : ~ \mathscr{S}^q_{T'} \to \mathscr{S}^q_T, \qquad \forall u \in {\bf L}^q_{T'},
\end{align}
and satisfies the consistency relations up to multiplicative constants of modulus $1$.

\vs

A couple of comments would be appropriate here. Firstly, Fock-Goncharov \cite{FG09} constructed an answer for such an intertwiner ${\bf K}^\hbar_k$ and proved desired properties, but that was an intertwiner for certain canonical but {\em reducible} representations of ${\bf L}^q_T$'s. As noted in \cite{FG09}, by general theory of the spectral decomposition of the operators representing the generators of the center of ${\bf L}^q_T$ enumerated by the punctures $\mathcal{P}$, one can easily expect the existence of an intertwiner ${\bf K}^\hbar_k$ for the irreducible representations for almost every choice of $\lambda : \mathcal{P} \to \mathbb{R}$ (eq.\eqref{eq:lambda}), solving the above problem. But, an explicit construction of such an intertwiner ${\bf K}^\hbar_k$ for the irreducible representations of ${\bf L}^q_T$'s (i.e. for one copy of the quantum Teichm\"uller space, instead of Fock-Goncharov's quantum symplectic double of the Teichm\"uller space) does not follow from this almost-everywhere existence, nor had it been systematically written down in the literature; after all, the irreducible representations themselves were not really dealt with. For our purposes, it is important to have an explicit operator ${\bf K}^\hbar_k$, instead of just an abstract existence statement, because we want to compose several of them and compute the resulting constants. Secondly, since the major result of the present paper is on the multiplicative constants appearing in the consistency relations satisfied by the intertwiners, we must describe these consistency relations in more detail, which we will do in \S\ref{subsect:consistency_relations}.

\subsection{The non-compact quantum dilogarithm function $\Phi^\hbar$}
\label{subsec:QD}

Here we recall an important ingredient for the mutation intertwiner, namely the quantum dilogarithm function of Faddeev-Kashaev \cite{FaKa}. We just give a definition and necessary properties, and refer the interested readers to \cite{FaKa}, \cite{FG09}, \cite{KN}, \cite{K16c}, and references therein.

\begin{proposition}[the non-compact quantum dilogarithm function \cite{Fa} \cite{Kash00}]
\label{prop:Phi_hbar}
Let $\hbar>0$. Define a complex function $\Phi^\hbar(z)$ as
$$
\Phi^\hbar(z) = \exp\left( - \frac{1}{4} \int_\Omega \frac{e^{-{\rm i} pz}}{\sinh(\pi p) \sinh(\pi \hbar p)} \frac{dp}{p} \right)
$$
for $z$ living in the strip $|{\rm Im}(z)|<\pi(1 + \hbar)$, where $\Omega$ is a contour along the real line that avoids the origin along a small half-circle above the origin.
\begin{enumerate}
\item[\rm (1)] Then $\Phi^\hbar(z)$ is a non-vanishing analytic function on this strip, satisfying the difference relations
$$
\Phi^\hbar(z+2\pi {\rm i}\hbar) = (1+e^{{\rm i} \pi\hbar} e^z) \Phi^\hbar(z), \qquad \Phi^\hbar(z+2\pi {\rm i}) = (1+ e^{{\rm i}\pi/\hbar} e^{z/\hbar}) \Phi^\hbar(z)
$$
in this strip. With the help of these relations, $\Phi^\hbar(z)$ analytically continues to a meromorphic function on the whole complex plane $\mathbb{C}$, called the \ul{\bf non-compact quantum dilogarithm}.

\item[\rm (2)] (the unitarity)
$$
|\Phi^\hbar(x)| = 1, \qquad \forall x\in \mathbb{R}.
$$

\item[\rm (3)] (the reflexivity)
$$
\Phi^\hbar(z) \Phi^\hbar(-z) = \alpha_\hbar^{-2} \cdot e^{z^2/(4\pi {\rm i} \hbar)}
$$
where
\begin{align}
\label{eq:alpha_hbar}
\alpha_\hbar = e^{\frac{{\rm i}\pi }{24} (\hbar + \hbar^{-1})} \in {\rm U}(1).
\end{align}

\item[\rm (4)] (the pentagon identity \cite{FKV} \cite{W} \cite{G}) If ${\bf P}$ and ${\bf Q}$ are self-adjoint operators on a separable Hilbert space $\mathscr{H}$ and satisfies the Weyl-relation-version of the commutation relation
$$
[{\bf P}, {\bf Q}] = 2\pi {\rm i} \hbar \cdot {\rm Id},
$$
i.e. if $e^{{\rm i} a {\bf P}} e^{{\rm i} b {\bf Q}} = e^{-2\pi {\rm i} ab \hbar} e^{{\rm i} b {\bf Q}} e^{{\rm i} a{\bf P}}$ holds for all $a,b\in \mathbb{R}$, then the following equality of unitary operators hold:
$$
\Phi^\hbar({\bf P}) \, \Phi^\hbar({\bf Q}) = \Phi^\hbar({\bf Q}) \, \Phi^\hbar({\bf P} + {\bf Q}) \, \Phi^\hbar({\bf P}).
$$
\end{enumerate}
\end{proposition}
Another version of the quantum dilogarithm is the `compact' version, defined for a complex parameter ${\bf q}$ with $|{\bf q}|<1$ as the meromorphic function
\begin{align}
\label{eq:Psi_q}
\Psi^{{\bf q}}(z) := \prod_{i=1}^\infty (1+ {\bf q}^{2i-1} z)^{-1};
\end{align}
see e.g. \cite{FaKa} \cite{FG09} \cite{KN} \cite{K16c} for its relationship with $\Phi^\hbar$.

\subsection{Our answer for the mutation intertwiner}
\label{subsec:our_answer_for_mutation_intertwiner}

We now describe our answer for the mutation intertwiner ${\bf K}^\hbar_k$, associated to the flip $T\leadsto T'$ along the edge $k$, for any choice of symplectic decompositions $(\ell_T,\mathcal{B}_T)$ for $T$ and $(\ell_{T'},\mathcal{B}_{T'})$ for $T'$. Following e.g. \cite{KN} and \cite{K16c}, we present two descriptions, depending on the choice of a `tropical' sign $\epsilon \in \{+,-\}$. We present the answer as the composition
$$
{\bf K}^{\hbar(\epsilon)}_k = \alpha_\hbar^\epsilon \, {\bf K}^{\sharp \hbar (\epsilon)}_k \circ {{\bf K}'}^{(\epsilon)}_{\hspace{-1mm}k} : \mathscr{H}_{\ell_{T'},\mathcal{B}_{T'}} \to \mathscr{H}_{\ell_T,\mathcal{B}_T}, \qquad \epsilon \in \{+,-\},
$$
in the style of \cite{FG09} and \cite{KN}, with
$$
{\bf K}^{\sharp \hbar (\epsilon)}_k : \mathscr{H}_{\ell_T,\mathcal{B}_T} \to \mathscr{H}_{\ell_T,\mathcal{B}_T}
$$
being called the \ul{\bf automorphism part},
$$
{{\bf K}'}^{(\epsilon)}_{\hspace{-1mm}k} : \mathscr{H}_{\ell_{T'}, \mathcal{B}_{T'}} \to \mathscr{H}_{\ell_T, \mathcal{B}_T}
$$
being called the \ul{\bf monomial transformation part}, where $\alpha_\hbar$ is as defined in eq.\eqref{eq:alpha_hbar}. The symbols $(\epsilon)$ appearing in ${\bf K}_k^{\hbar(\epsilon)}$, ${\bf K}_k^{\sharp \hbar (\epsilon)}$, and ${{\bf K}'}^{(\epsilon)}_{\hspace{-1mm}k}$ are just labels but not actual numbers, while $\epsilon \in \{+,-\}$ appearing in $\alpha_\hbar^\epsilon$ stands for the number $+1,-1$ respectively. Note that the operators ${\bf K}_k^{\hbar(\epsilon)}$, ${\bf K}_k^{\sharp \hbar (\epsilon)}$, and ${{\bf K}'}^{(\epsilon)}_{\hspace{-1mm}k}$ all depend on the choices $(\ell_T,\mathcal{B}_T)$ and $(\ell_{T'},\mathcal{B}_{T'})$. To emphasize this dependence, one might put these symplectic decompositions in the subscripts of these operators, like
$$
{\bf K}_k^{\hbar(\epsilon)} = {\bf K}_{(\ell_T,\mathcal{B}_T),(\ell_{T'},\mathcal{B}_{T'})}^{\hbar(\epsilon)}, \quad
{{\bf K}'}^{(\epsilon)}_{\hspace{-1mm}k}  = {{\bf K}'}^{(\epsilon)}_{\hspace{-1mm}(\ell_T,\mathcal{B}_T),(\ell_{T'},\mathcal{B}_{T'})}, \quad 
{\bf K}^{\sharp \hbar (\epsilon)}_k = {\bf K}^{\sharp \hbar (\epsilon)}_{(\ell_T,\mathcal{B}_T),(\ell_{T'},\mathcal{B}_{T'})}.
$$

\vs

The {\em construction of the automorphism part} is via the functional calculus for the self-adjoint operator $\epsilon \, \wh{x}^\hbar_k = \epsilon \, \wh{x}^\hbar_{k;\ell_T,\mathcal{B}_T}$ applied to the non-compact quantum dilogarithm function:
\begin{align}
\label{eq:automorphism_part}
{\bf K}^{\sharp \hbar (\epsilon)}_k := \Phi^\hbar(\epsilon\, \wh{x}^\hbar_k)^\epsilon : \mathscr{H}_{\ell_T,\mathcal{B}_T} \to \mathscr{H}_{\ell_T,\mathcal{B}_T},
\end{align}
where $\epsilon=+,-$ in the right hand side means $\epsilon =+1,-1$ respectively. From the property Prop.\ref{prop:Phi_hbar}(2) of the quantum dilogarithm function, we see that ${\bf K}^{\sharp \hbar (\epsilon)}_k$ is a unitary operator on $\mathscr{H}_{\ell_T,\mathcal{B}_T}$.

\vs

For the monomial transformation part, we consider the two `cluster (tropical)' linear maps
\begin{align}
\nonumber
& C^{(\epsilon)}_k : V_{T'} \to V_T, \\
\label{eq:C_k_epsilon}
& C^{(\epsilon)}_k (x_e') = \left\{
\begin{array}{ll}
- x_k & \mbox{for $e\in T'$ with $e=k$,} \\
x_e + [\epsilon \, \varepsilon_{ek}]_+ \, x_k & \mbox{for $e\in T'$ with $e\neq k$},
\end{array}
\right.
\end{align}
where $T$ and $T'$ are naturally identified as sets, $(\varepsilon_{ef})_{e,f \in T},(\varepsilon'_{ef})_{e,f \in T'}$ are the exchange matrices for $T,T'$ respectively, and the symbol $[a]_+$ denotes the positive part of a real number $a$:
$$
[a]_+ = \left\{
\begin{array}{ll}
a & \mbox{if $a \ge 0$,} \\
0 & \mbox{if $a<0$,}
\end{array}
\right.
$$
which can also be defined as $[a]_+ = \frac{a+|a|}{2}$. A simple but important observation is that these linear maps $C_k^{(\epsilon)}$ depend only on $T$ and $T'$, but not on the choice of symplectic decompositions $(\ell_T,\mathcal{B}_T)$ and $(\ell_{T'},\mathcal{B}_{T'})$. We also observe the following basic properties.
\begin{lemma}
\label{lem:C_k_epsilon_is_symplectic}
One has:
\begin{enumerate}
\item[\rm (1)] (see Lem.2.7 of \cite{FG09}) $C^{(\epsilon)}_k$ is invertible, and respects the forms $B_{T'}$ and $B_T$.

\item[\rm (2)] (from \cite[Lem.24]{BL} \cite[Prop.14]{Liu}) $C_k^{(\epsilon)}$ sends $x'_p \in V_{T'}^\perp$ to $x_p \in V_T^\perp$, $\forall p \in \mathcal{P}$, where $x'_p$ and $x_p$ are the puncture elements defined in Lem.\ref{lem:puncture_elements}. In particualr the functions $f_T : V_T^\perp \to \mathbb{R}$ and $f_{T'} : V_{T'}^\perp \to \mathbb{R}$ are related by $C_k^{(\epsilon)}$.
\end{enumerate}
\end{lemma}

{\it Proof.} (1) It is an easy check that the inverse map $(C_k^{(\epsilon)})^{-1} : V_T \to V_{T'}$ is given by the formula
\begin{align}
\label{eq:C_k_epsilon_inverse}
( C_k^{(\epsilon)} )^{-1} (x_k) = \left\{
\begin{array}{ll}
- x_k' & \mbox{if $e=k$,} \\
x_e' + [\epsilon \, \varepsilon_{ek}]_+ \, x_k' & \mbox{if $e\neq k$,}
\end{array}
\right.
\end{align}
Indeed, for example, $(C_k^{(\epsilon)}) (C_k^{(\epsilon)})^{-1}(x_k) = C_k^{(\epsilon)}(-x_k') = x_k$, and $(C_k^{(\epsilon)}) (C_k^{(\epsilon)})^{-1}(x_e) = C_k^{(\epsilon)}(x_e' + [\epsilon \, \varepsilon_{ek}]_+ \, x_k') = (x_e + [\epsilon \, \varepsilon_{ek}]_+ \, x_k) + [\epsilon\,\varepsilon_{ek}]_+ \, (-x_k) = x_e$ for each $e\neq k$. 

\vs

We recall the transformation formula for the exchange matrices (see e.g. \cite{FG06})
\begin{align}
\label{eq:varepsilon_prime}
\varepsilon'_{ef} = \left\{
\begin{array}{ll}
- \varepsilon_{ef} & \mbox{if $k\in \{e,f\},$} \\
\varepsilon_{ef} + \frac{ |\varepsilon_{ek}| \, \varepsilon_{kf} + \varepsilon_{ek} \, |\varepsilon_{kf}|}{2} & \mbox{if $k\notin \{e,f\}$}.
\end{array}
\right.
\end{align}
Now let's check whether
$$
B_{T'}( x_e', x_f' ) = B_T( C_k^{(\epsilon)}(x_e') ,\, C_k^{(\epsilon)}(x_f') )
$$
holds for all $e,f\in T'$; the left hand side is $\varepsilon'_{ef}$. In case $f=k$, the left hand side is $\varepsilon'_{ek} = - \varepsilon_{ek}$. If also $e=k$, then both sides are equal, due to the skew-symmetry. If $e\neq k$, then the right hand side is $B_T(x_e + [\epsilon \, \varepsilon_{ek}]_+ \, x_k, \, -x_k) = B_T(x_e,\, - x_k) = - \varepsilon_{ek}$, coinciding with the left hand side. Now suppose $f\neq k$. The case $e=k$ is dealt with by the skew-symmetry, so let $e\neq k$. Then the right hand side is
\begin{align*}
& B_T(x_e+[\epsilon \, \varepsilon_{ek}]_+\, x_k, \, x_f + [\epsilon \, \varepsilon_{fk}]_+ \, x_k) = \varepsilon_{ef} + [\epsilon \, \varepsilon_{fk}]_+ \, \varepsilon_{ek} + [\epsilon\,\varepsilon_{ek}]_+ \, \varepsilon_{kf} \\
& = \varepsilon_{ef} + \frac{ \epsilon \, \varepsilon_{fk} + |\epsilon \, \varepsilon_{fk}|}{2} \, \varepsilon_{ek} + \frac{ \epsilon \, \varepsilon_{ek} + |\epsilon \, \varepsilon_{ek}|}{2} \, \varepsilon_{kf} = \varepsilon_{ef} + \frac{|\varepsilon_{kf}| \varepsilon_{ek} + |\varepsilon_{ek}| \varepsilon_{kf}}{2} = \varepsilon'_{ef},
\end{align*}
coinciding with the left hand side. 

\vs

(2) Note $x'_p = \sum_{e\in T'} {\sigma'}_{\hspace{-1mm}p}^e x_e' \in V_{T'}^\perp$ and $x_p = \sum_{e\in T} \sigma_p^e x_e \in V_T$. Note
$$
C_k^{(\epsilon)}(x_p') = \ssum_{e\in T'} {\sigma'}^e_{\hspace{-1mm}p} C_k^{(\epsilon)}(x_e')
= - {\sigma'}_{\hspace{-2mm}p}^k x_k + \ssum_{e\neq k} {\sigma'}^e_{\hspace{-2mm}p} (x_e + [\epsilon\, \varepsilon_{ek}]_+ x_k).
$$
It is easy to see that ${\sigma'}^e_{\hspace{-2mm}p} = \sigma^e_p$ holds for all $e\neq k$, hence it remains to show that $\sigma_p^k = -{\sigma'}^k_{\hspace{-2mm}p} + \sum_{e\neq k} {\sigma'}^e_{\hspace{-2mm}p} [\epsilon \, \varepsilon_{ek}]_+$ holds. We leave it as a straightforward exercise on the combinatorics of ideal triangulations; the result itself has been already known by indirect arguments, e.g. in \cite[Lem.24]{BL}, \cite[Prop.14]{Liu}, \cite{FG06}, and \cite[Lem.3.3]{AK17}. \qed

\begin{remark}
The same comment as in Rem.\ref{rem:mutation_only_for_regular_flips} applies to the transformation formula eq.\eqref{eq:varepsilon_prime}.
\end{remark}

As a result of this lemma, each of the pairs 
\begin{align}
\label{eq:pulllback_Lagrangian_decompositions}
C_k^{(\epsilon)} (\ell_{T'},\mathcal{B}_{T'}) = ( C_k^{(\epsilon)} (\ell_{T'}), C_k^{(\epsilon)} (\mathcal{B}_{T'})), \qquad \mbox{for} \quad \epsilon \in \{+,-\},
\end{align}
forms a symplectic decomposition for $T$, or more precisely, that of $(V_{T},B_{T})$. Thus we can consider the Weil intertwiner constructed in \S\ref{subsec:LSSW_intertwiner}:
$$
{\bf F} = {\bf F}_{C_k^{(\epsilon)}(\ell_{T'},\mathcal{B}_{T'}),(\ell_T,\mathcal{B}_T)} ~ : ~ \mathscr{H}_{C_k^{(\epsilon)}(\ell_{T'},\mathcal{B}_{T'})} \longrightarrow \mathscr{H}_{\ell_T,\mathcal{B}_T}.
$$
Now consider the pullback map ${\bf R}$ considered in Def.\ref{def:R} (with our $C_k^{(\epsilon)}$ playing the role of $C^{-1}$ there):
$$
{\bf R} = {\bf R}_{(\ell_{T'},\mathcal{B}_{T'}), C_k^{(\epsilon)}(\ell_{T'},\mathcal{B}_{T'})} ~ :~ \mathscr{H}_{\ell_{T'},\mathcal{B}_{T'}} \longrightarrow \mathscr{H}_{C_k^{(\epsilon)}(\ell_{T'},\mathcal{B}_{T'})}
$$
which can be understood either as representing the identity map $L^2(\mathbb{R}^r) \to L^2(\mathbb{R}^r)$, or as the pullback map $((C_k^{(\epsilon)})^{-1})^* : L^2({\rm span}_\mathbb{R} \mathcal{B}_{T'}) \to L^2({\rm span}_\mathbb{R} C_k^{(\epsilon)}(\mathcal{B}_{T'}))$ as in Lem.\ref{lem:R_as_pullback_map}.

\vs

We finally give a {\em construction of the monomial transformation part} as follows:
\begin{align}
\label{eq:monomial_transformation_part}
{{\bf K}'}^{(\epsilon)}_{\hspace{-1mm} k} := {\bf F}_{C_k^{(\epsilon)}(\ell_{T'},\mathcal{B}_{T'}),(\ell_T,\mathcal{B}_T)} \circ {\bf R}_{(\ell_{T'},\mathcal{B}_{T'}), C_k^{(\epsilon)}(\ell_{T'},\mathcal{B}_{T'})} ~ :~ \mathscr{H}_{\ell_{T'},\mathcal{B}_{T'}} \longrightarrow \mathscr{H}_{\ell_T,\mathcal{B}_T}. 
\end{align}
In a diagram format, one can write this map ${{\bf K}'}^{(\epsilon)}_{\hspace{-1mm} k}$ as the following composition
$$
\xymatrix@R+5mm{
\mathscr{H}_{\ell_{T'},\mathcal{B}_{T'}} \ar[r]^-{{\bf R}} & \mathscr{H}_{C_k^{(\epsilon)}(\ell_{T'},\mathcal{B}_{T'})} \ar[r]^-{{\bf F}} & \mathscr{H}_{\ell_T,\mathcal{B}_T}.
}
$$
One can immediately observe that this operator ${{\bf K}'}^{(\epsilon)}_{\hspace{-1mm} k}$ is a generalization of the operator $\rho(C)$ in eq.\eqref{eq:rho_C} which appeared in the representation of the symplectic group, in \S\ref{subsec:projective_representation_of_the_symplectic_group}. This operator written as ${\bf FR}$ can also be written as ${\bf RF}$, by Lem.\ref{lem:compatibility_between_F_and_R}; more precisely, as ${\bf R}_{(C_k^{(\epsilon)})^{-1}(\ell_T,\mathcal{B}_T), (\ell_T,\mathcal{B}_T)} \, {\bf F}_{(\ell_{T'},\mathcal{B}_{T'}),(C_k^{(\epsilon)})^{-1}(\ell_T,\mathcal{B}_T)}$.

\vs

We claim that the operator ${\bf K}_k^{\hbar(\epsilon)}$ constructed in this subsection is indeed an interwining operator.
\begin{proposition}[the mutation intertwiner for the irreducible representations]
\label{prop:intertwining}
For each choice of a sign $\epsilon \in \{+,-\}$, the operator ${\bf K}^{\hbar(\epsilon)}_k = \alpha_\hbar^\epsilon \, {\bf K}^{\sharp \hbar (\epsilon)}_k \circ {{\bf K}'}^{(\epsilon)}_{\hspace{-1mm}k} : \mathscr{H}_{\ell_{T'},\mathcal{B}_{T'}} \to \mathscr{H}_{\ell_T,\mathcal{B}_T}$ constructed above satisfies the intertwining equations, i.e. eq.\eqref{eq:intertwining_equations}.
\end{proposition}

Very crucial for the purposes of the present paper is the equality of these two answers ${\bf K}^{\hbar(+)}_k$ and ${\bf K}^{\hbar(-)}_k$ for the mutation intertwiner problem.
\begin{proposition}[the equality of the two signed decompositions of the mutation intertwiner]
\label{prop:equality_of_two_signed_decompositions}
One has
$$
{\bf K}^{\hbar(+)}_k ~\sim~ {\bf K}^{\hbar(-)}_k
$$
i.e. these two operators coincide up to a multiplicative constant which is a power of $e^{{\rm i}\pi/4}$.
\end{proposition}

Finally, upon the change of the choice of symplectic decompositions for $T$ and $T'$, our intertwiner is compatible with the Weil intertwiners, in the following sense.
\begin{proposition}[the compatibility between the mutation intertwiners and the Weil intertwiners]
\label{prop:compatibility_of_mutation_intertwiner_and_Weil_intertwiners}
Let $\epsilon \in \{+,-\}$. For any two symplectic decompositions $(\ell_T,\mathcal{B}_T)$ and $(\ol{\ell}_T, \ol{\mathcal{B}}_T)$ for $T$, and for any two symplectic decompositions $(\ell_{T'},\mathcal{B}_{T'})$ and $(\ol{\ell}_{T'},\ol{\mathcal{B}}_{T'})$ for $T'$, the following diagram commutes, up to a multiplicative constant that is a power of $e^{{\rm i}\pi/4}$:
\begin{align}
\label{eq:K_and_F_diagram}
\begin{array}{c}
\xymatrix@C+25mm@R+3mm{
\mathscr{H}_{\ell_{T'},\mathcal{B}_{T'}} \ar[r]^-{{\bf K}^{\hbar(\epsilon)}_k =\, {\bf K}^{\hbar(\epsilon)}_{(\ell_T,\mathcal{B}_T),(\ell_{T'},\mathcal{B}_{T'})}} \ar[d]_{{\bf F}} & \mathscr{H}_{\ell_T,\mathcal{B}_T} \ar[d]^{{\bf F}}  \\
\mathscr{H}_{\ol{\ell}_{T'}, \ol{\mathcal{B}}_{T'}} \ar[r]^-{{\bf K}^{\hbar(\epsilon)}_k = \, {\bf K}^{\hbar(\epsilon)}_{(\ol{\ell}_T, \ol{\mathcal{B}}_T), (\ol{\ell}_{T'}, \ol{\mathcal{B}}_{T'})}} & \mathscr{H}_{\ol{\ell}_T,\ol{\mathcal{B}}_T}
}
\end{array}
\end{align}
where the vertical arrows are the Weil intertwiners.
\end{proposition}

{\it Proof of Prop.\ref{prop:compatibility_of_mutation_intertwiner_and_Weil_intertwiners}.} We consider the following factored version diagram, consisting of three squares:
$$
\xymatrix@R-1mm{
\mathscr{H}_{\ell_{T'},\mathcal{B}_{T'}} \ar[r]^-{{\bf R}}  \ar[d]_{{\bf F}} & \mathscr{H}_{C_k^{(\epsilon)}(\ell_{T'},\mathcal{B}_{T'})} \ar[r]^-{{\bf F}}  \ar[d]^{{\bf F}} & \mathscr{H}_{\ell_T,\mathcal{B}_T} \ar[r]^-{{\bf K}_k^{\sharp \hbar (\epsilon)}}  \ar[d]_{{\bf F}} & \mathscr{H}_{\ell_T,\mathcal{B}_T}   \ar[d]^{{\bf F}} \\
\mathscr{H}_{\ol{\ell}_{T'}, \ol{\mathcal{B}}_{T'}} \ar[r]^-{{\bf R}} & \mathscr{H}_{C_k^{(\epsilon)}(\ol{\ell}_{T'},\ol{\mathcal{B}}_{T'})} \ar[r]^-{{\bf F}} & \mathscr{H}_{\ol{\ell}_T, \ol{\mathcal{B}}_T} \ar[r]^-{{\bf K}_k^{\sharp \hbar (\epsilon)}} & \mathscr{H}_{\ol{\ell}_T, \ol{\mathcal{B}}_T} \\
}
$$
We see that the left square is commutative by Lem.\ref{lem:compatibility_between_F_and_R}. The middle square is commutative up to a power of $e^{{\rm i} \pi /4}$, by Cor.\ref{cor:relationship_between_two_compositions_of_F}. For the right square, observe that ${\bf F} : \mathscr{H}_{\ell_T,\mathcal{B}_T} \to \mathscr{H}_{\ol{\ell}_T,\ol{\mathcal{B}}_T}$ is an intertwiner of representations, in the sense that the equalities
\begin{align}
\label{eq:conjugation_by_F}
{\bf F}_{(\ell_T,\mathcal{B}_T), (\ol{\ell}_T,\ol{\mathcal{B}}_T)} \, \pi_{\ell_T,\mathcal{B}_T}(x) \, ({\bf F}_{(\ell_T,\mathcal{B}_T), (\ol{\ell}_T,\ol{\mathcal{B}}_T)})^{-1} = \pi_{\ol{\ell}_T, \ol{\mathcal{B}}_T}(x), \qquad \forall x\in \mathfrak{n}_T,
\end{align}
hold as equalities of self-adjoint operators. Therefore
\begin{align*}
{\bf F}  \, {\bf K}_k^{\sharp\hbar (\epsilon)} \, {\bf F}^{-1} = {\bf F} \, {\bf K}_{(\ell_T,\mathcal{B}_T),(\ell_{T'},\mathcal{B}_{T'})}^{\sharp \hbar(\epsilon)} \, {\bf F}^{-1}
& = {\bf F} \, \Phi^\hbar( \epsilon \, {\rm i} \sqrt{2 \pi \hbar} \cdot \pi_{\ell_T,\mathcal{B}_T}(x_k)  )^\epsilon \, {\bf F}^{-1} \\
& = \Phi^\hbar(  {\bf F}\, (\epsilon \, {\rm i} \sqrt{2 \pi \hbar} \cdot \pi_{\ell_T,\mathcal{B}_T}(x_k) ) \, {\bf F}^{-1} )^\epsilon \\
& = \Phi^\hbar(  \epsilon \, {\rm i} \sqrt{2 \pi \hbar} \cdot \pi_{\ol{\ell}_T,\ol{\mathcal{B}}_T}(x_k) )^\epsilon  \\
& = {\bf K}_{(\ol{\ell}_T,\ol{\mathcal{B}}_T),(\ol{\ell}_{T'},\ol{\mathcal{B}}_{T'})}^{\sharp \hbar (\epsilon)}
= {\bf K}_k^{\sharp \hbar (\epsilon)}.
\end{align*}
Here we used the fact that the conjugation by a unitary operator commutes with functional calculus; that is, if $U$ is a unitary operator, $\wh{x}$ a self-adjoint operator, and $\Phi :\mathbb{R} \to {\rm U}(1)$ a measurable function, then
\begin{align}
\label{eq:conjugation_commutes_with_functional_calculus}
U \, \Phi(\wh{x}) \, U^{-1} = \Phi(U\, \wh{x} \, U^{-1})
\end{align}
holds (see e.g. Lem.47 in the version 2 of \cite{K16c}, i.e. arXiv:1602.00797v2). So the right square of the diagram is commutative. \qed

\vs

The proofs of Prop.\ref{prop:intertwining} and Prop.\ref{prop:equality_of_two_signed_decompositions} will occupy the rest of the present section.

\subsection{The verification of the intertwining equations}
\label{subsec:verification_of_intertwining_equations}

Here we provide a proof of Prop.\ref{prop:intertwining}. Core of the proof has been established by Fock-Goncharov \cite{FG09} for their non-irreducible representations, and we just present how to adapt it to our situation. We inherit the notations from the previous subsection. In particular, we work with two triangulations $T$ and $T'$ related by the flip along the edge $k$. We begin with an easy step, which will be used later again in the present paper.
\begin{lemma}[the conjugation action of the monomial transformation part ${\bf K}'_k$]
\label{lem:conjugation_of_primed_part}
The monomial transformation part operator ${{\bf K}'}^{(\epsilon)}_{\hspace{-1mm}k} : \mathscr{H}_{\ell_{T'},\mathcal{B}_{T'}} \to \mathscr{H}_{\ell_T,\mathcal{B}_T}$ satisfies
$$
{{\bf K}'}^{(\epsilon)}_{\hspace{-1mm}k} \, \pi^\hbar_{\ell_{T'},\mathcal{B}_{T'}}(x') \, ({{\bf K}'}^{(\epsilon)}_{\hspace{-1mm}k})^{-1}
= \pi^\hbar_{\ell_T,\mathcal{B}_T}( C_k^{(\epsilon)}(x') ), \qquad \forall x' \in \mathfrak{n}_{T'},
$$
which are equalities of self-adjoint operators. In particular, we have
\begin{align}
\label{eq:conjugation_of_primed_part_on_generators}
{{\bf K}'}^{(\epsilon)}_{\hspace{-1mm}k} \, {\wh{x}{}'}^{\hbar}_{\hspace{-1mm} e} \, ({{\bf K}'}^{(\epsilon)}_{\hspace{-1mm}k})^{-1}
= \left\{
\begin{array}{ll}
- \wh{x}^\hbar_k, & \mbox{if $e =k$}, \\
\wh{x}^\hbar_e + [\epsilon \, \varepsilon_{ek}]_+ \, \wh{x}^\hbar_k, & \mbox{if $e\neq k$}, 
\end{array}
\right.
\end{align}

\end{lemma}

{\it Proof.} Let $x' \in \mathfrak{n}_{T'}$. In view of the definition in eq.\eqref{eq:monomial_transformation_part} of ${{\bf K}'}^{(\epsilon)}_{\hspace{-1mm}k}$, we first investigate the conjugation action by ${\bf R}$, then that by ${\bf F}$. From Cor.\ref{cor:R_is_pullback}, one gets
\begin{align*}
{\bf R}_{(\ell_{T'},\mathcal{B}_{T'}), C_k^{(\epsilon)}(\ell_{T'},\mathcal{B}_{T'})} \, \pi^\hbar_{\ell_{T'},\mathcal{B}_{T'}}(x') \, ({\bf R}_{(\ell_{T'},\mathcal{B}_{T'}), C_k^{(\epsilon)}(\ell_{T'},\mathcal{B}_{T'})})^{-1} = \pi^\hbar_{C_k^{(\epsilon)}(\ell_{T'},\mathcal{B}_{T'})}( C_k^{(\epsilon)}(x') )
\end{align*}
From eq.\eqref{eq:conjugation_by_F}, one gets
\begin{align*}
{\bf F}_{C_k^{(\epsilon)}(\ell_{T'},\mathcal{B}_{T'}), (\ell_T,\mathcal{B}_T)} \, \pi^\hbar_{C_k^{(\epsilon)}(\ell_{T'},\mathcal{B}_{T'})}( C_k^{(\epsilon)}(x') ) \, ({\bf F}_{C_k^{(\epsilon)}(\ell_{T'},\mathcal{B}_{T'}), (\ell_T,\mathcal{B}_T)})^{-1} = \pi^\hbar_{\ell_T,\mathcal{B}_T}( C_k^{(\epsilon)}(x') )
\end{align*}
finishing the proof. \qed

\vs

\begin{lemma}[the conjugation action of the automorphism part ${\bf K}^{\sharp \hbar (\epsilon)}_k$; \cite{FG09}]
\label{lem:conjugation_action_of_the_automorphism_part}
The automorphism part operator ${\bf K}^{\sharp \hbar (\epsilon)}_k : \mathscr{H}_{\ell_T,\mathcal{B}_T} \to \mathscr{H}_{\ell_T,\mathcal{B}_T}$ satisfies the following: for each $u \in {\bf L}^q_T$ and $\psi \in D_{\ell_T,\mathcal{B}_T}$, 
\begin{align}
\label{eq:conjugation_action_of_automorphism_part}
{\bf K}^{\sharp\hbar(\epsilon)}_k \, \pi^q_{\ell_T,\mathcal{B}_T}(u) \, \psi = \pi^q_{\ell_T,\mathcal{B}_T} ({\rm Ad}_{\Psi^q(\wh{X}_k^\epsilon)^\epsilon}(u)) \, {\bf K}^{\sharp\hbar(\epsilon)}_k \, \psi
\end{align}
where $\Psi^q$ is the compact quantum dilogarithm as defined in eq.\eqref{eq:Psi_q}.
\end{lemma}
Here, $\Psi^q(\wh{X}_k^\epsilon)^\epsilon = (1+q\wh{X}_k^\epsilon)^{-\epsilon} (1+q^3 \wh{X}_k^\epsilon)^{-\epsilon} \cdots$ is only a formal element, but ${\rm Ad}_{\Psi^q(\wh{X}_k^\epsilon)^\epsilon}(u)$ is well-defined. A rigorous proof of Lem.\ref{lem:conjugation_action_of_the_automorphism_part} requires somewhat heavy arguments in analysis. As we will not really use the above intertwining equations in the present paper, we only claim that Fock-Goncharov's proof of \cite[Thm.5.6]{FG09} should work for this lemma almost verbatim. In view of the relationship between $\Phi^\hbar$ and $\Psi^q$ (see \cite{Fa} \cite{FG09} \cite{K16c}), one can observe at least formally that the conjugation by ${\bf K}^{\sharp\hbar(\epsilon)}_k = \Phi^\hbar(\epsilon \, \wh{x}^\hbar_k)^\epsilon$ would yield the same result as the conjugation by $\Psi^q(\wh{X}_k^\epsilon)^\epsilon$; see \cite{FG09}. We also claim that Lem.\ref{lem:conjugation_action_of_the_automorphism_part} holds for all vectors $\psi$ in the Schwartz space $\mathscr{S}^q_T$, via the arguments in \cite[\S5.4]{FG09}. Combining Lem.\ref{lem:conjugation_of_primed_part} and Lem.\ref{lem:conjugation_action_of_the_automorphism_part} would provide a proof of the sought-for Prop.\ref{prop:intertwining}.

\subsection{The equality of the two signed decompositions of the mutation intertwiner}
\label{subsec:equality_of_two_signed_decompositions}

The present subsection is solely devoted to a proof of Prop.\ref{prop:equality_of_two_signed_decompositions}, i.e. the equality ${\bf K}^{\hbar(+)}_k = {\bf K}^{\hbar(-)}_k$ up to $e^{{\rm i}\pi/4}$. So we should show $\alpha_\hbar \, {\bf K}_k^{\sharp \hbar (+)} \, {{\bf K}'}_{\hspace{-1mm}k}^{(+)} \sim \alpha_\hbar^{-1} {\bf K}_k^{\sharp \hbar (-)} \, {{\bf K}'}_{\hspace{-1mm}k}^{(-)}$, or equivalently,
$$
\alpha_\hbar^2 \left({\bf K}_k^{\sharp \hbar (-)} \right)^{-1} \, {\bf K}_k^{\sharp \hbar (+)} ~\sim~ {{\bf K}'}_{\hspace{-1mm}k}^{(-)} \, \left({{\bf K}'}_{\hspace{-1mm}k}^{(+)} \right)^{-1} : \mathscr{H}_{\ell_T,\mathcal{B}_T} \longrightarrow \mathscr{H}_{\ell_T,\mathcal{B}_T}.
$$
In view of the definition in eq.\eqref{eq:automorphism_part} of the automorphism parts, the left hand side is
\begin{align*}
\alpha_\hbar^2 \left({\bf K}_k^{\sharp \hbar (-)} \right)^{-1} \, {\bf K}_k^{\sharp \hbar (+)} = \alpha_\hbar^2 \, \Phi^\hbar(-\wh{x}^\hbar_k) \, \Phi^\hbar(\wh{x}^\hbar_k) = e^{(\wh{x}^\hbar_k)^2/(4\pi {\rm i} \hbar)},
\end{align*}
by the property Prop.\ref{prop:Phi_hbar}(3) of the quantum dilogarithm, where $e^{(\wh{x}^\hbar_k)^2/(4\pi {\rm i} \hbar)}$ is the unitary operator made sense via the functional calculus for the self-adjoint operator $\wh{x}^\hbar_k$. In view of the definition in eq.\eqref{eq:monomial_transformation_part} of the monomial transformation parts, the right hand side is
\begin{align}
\label{eq:equality_of_signed_decompositions_monomial_part_RHS}
\begin{array}{l}
{{\bf K}'}_{\hspace{-1mm}k}^{(-)} \, \left({{\bf K}'}_{\hspace{-1mm}k}^{(+)} \right)^{-1} \\
= 
{\bf F}_{C_k^{(-)}(\ell_{T'}, \mathcal{B}_{T'}), (\ell_T,\mathcal{B}_T)} \, \underbrace{ {\bf R}_{(\ell_{T'},\mathcal{B}_{T'}),C_k^{(-)}(\ell_{T'},\mathcal{B}_{T'})} \, {\bf R}_{(\ell_{T'},\mathcal{B}_{T'}), C_k^{(+)}(\ell_{T'},\mathcal{B}_{T'})}^{-1} }   \,
{\bf F}_{C_k^{(+)}(\ell_{T'},\mathcal{B}_{T'}), (\ell_T,\mathcal{B}_T)}^{-1}
\end{array}
\end{align}
By the composition identity (Lem.\ref{lem:R_composition_identity}, eq.\eqref{eq:R_composition_identity}) for ${\bf R}$'s, the underbraced part equals
$$
{\bf R}_{C_k^{(+)}(\ell_{T'}, \mathcal{B}_{T'}), C_k^{(-)}(\ell_{T'}, \mathcal{B}_{T'})} ~:~ \mathscr{H}_{C_k^{(+)}(\ell_{T'},\mathcal{B}_{T'})} \longrightarrow \mathscr{H}_{C_k^{(-)}(\ell_{T'},\mathcal{B}_{T'})}
$$
which, through Lem.\ref{lem:R_as_pullback_map}, can be understood as the pullback map
$$
C_k^* : L^2({\rm span}_\mathbb{R} C_k^{(+)}(\mathcal{B}_{T'})) \to L^2({\rm span}_\mathbb{R} C_k^{(-)}(\mathcal{B}_{T'}))
$$
where $C_k : {\rm span}_\mathbb{R} C_k^{(-)}(\mathcal{B}_{T'}) \to {\rm span}_\mathbb{R} C_k^{(+)}(\mathcal{B}_{T'})$ is the restriction of the map
$$
C_k := C_k^{(+)} \circ (C_k^{(-)})^{-1} : V_T \to V_T.
$$
Recall that the map ${\bf R}$ by itself is in fact just the identity map when its domain $\mathscr{H}_{C_k^{(+)}(\ell_{T'},\mathcal{B}_{T'})}$ and the codomain $\mathscr{H}_{C_k^{(-)}(\ell_{T'},\mathcal{B}_{T'})}$ are written as $L^2(\mathbb{R}^r)$. So, one observes that the nontriviality of the computation of ${{\bf K}'}_{\hspace{-1mm}k}^{(-)} \, \left({{\bf K}'}_{\hspace{-1mm}k}^{(+)} \right)^{-1}$ is due to the discrepancy between this map  ${\bf R} : \mathscr{H}_{C_k^{(+)}(\ell_{T'},\mathcal{B}_{T'})} \to \mathscr{H}_{C_k^{(-)}(\ell_{T'},\mathcal{B}_{T'})}$ and the Weil intertwiner ${\bf F} : \mathscr{H}_{C_k^{(+)}(\ell_{T'},\mathcal{B}_{T'})} \to \mathscr{H}_{C_k^{(-)}(\ell_{T'},\mathcal{B}_{T'})}$ with the same domain and the codomain.

\vs

We first compute the linear map $C_k : V_T \to V_T$; note again that this map is completely determined by $T$ and $T'$, and is independent of the choice of symplectic decompositions $(\ell_T,\mathcal{B}_T)$ and $(\ell_{T'},\mathcal{B}_{T'})$. From the formulas in eq.\eqref{eq:C_k_epsilon} and eq.\eqref{eq:C_k_epsilon_inverse}, one obtains
\begin{align*}
C_k(x_e) = (C_k^{(+)} \circ (C_k^{(-)})^{-1})(x_e)
& = \left\{
\begin{array}{ll}
x_k & \mbox{if $e=k$}, \\
x_e + ([\varepsilon_{ek}]_+ - [-\varepsilon_{ek}]_+) x_k & \mbox{if $e\neq k$},
\end{array}
\right. \\
& = \left\{
\begin{array}{ll}
x_k & \mbox{if $e=k$}, \\
x_e + \varepsilon_{ek} \, x_k & \mbox{if $e\neq k$};
\end{array}
\right.
\end{align*}
in the last equality, we used the fact $[a]_+ - [-a]_+ = a$ which holds for any real number $a$. Note from eq.\eqref{eq:B_T} that $\varepsilon_{ef} = B_T(x_e,x_f)$, and hence by linearity one can easily check that the above map $C_k : V_T \to V_T$ is given by the following neat formula:
\begin{align}
\label{eq:C_k_neat}
C_k(x) = x + B_T(x,x_k) \, x_k, \qquad \forall x\in V_T.
\end{align}
Before proceeding, let us summarize some properties of $C_k$. As each $C_k^{(+)}$ and $C_k^{(-)}$ is an invertible symplectic transformation (by Lem.\ref{lem:C_k_epsilon_is_symplectic}), so is $C_k$. The two symplectic decompositions $C_k^{(+)}(\ell_{T'},\mathcal{B}_{T'})$ and $C_k^{(-)}(\ell_{T'},\mathcal{B}_{T'})$ are related by the map $C_k$, i.e.
$$
C_k (C_k^{(-)}(\ell_{T'}, \mathcal{B}_{T'})) = C_k^{(+)}(\ell_{T'}, \mathcal{B}_{T'}).
$$

\vs

In eq.\eqref{eq:equality_of_signed_decompositions_monomial_part_RHS}, the rightmost factor ${\bf F}_{C_k^{(+)}(\ell_{T'},\mathcal{B}_{T'}), (\ell_T,\mathcal{B}_T)}^{-1}$ in the right hand side can be written as
$$
{\bf F}_{C_k^{(+)}(\ell_{T'},\mathcal{B}_{T'}), (\ell_T,\mathcal{B}_T)}^{-1} =  {\bf F}_{(\ell_T,\mathcal{B}_T), C_k^{(+)}(\ell_{T'},\mathcal{B}_{T'})} : \mathscr{H}_{\ell_T,\mathcal{B}_T} \to \mathscr{H}_{C_k^{(+)}(\ell_{T'},\mathcal{B}_{T'})},
$$
by Prop.\ref{prop:Weil_intertwiner_involutivity}. Further, one can write $C_k^{(+)}(\ell_{T'},\mathcal{B}_{T'})$ as $C_k(C_k^{(-)}(\ell_{T'},\mathcal{B}_{T'}))$. So one can represent the composition in the right hand side of eq.\eqref{eq:equality_of_signed_decompositions_monomial_part_RHS} by the following diagram, keeping in mind that a Weil intertwiner can be denoted just by the symbol ${\bf F}$ when its domain and the codomain are explicitly specified:
\begin{align}
\label{eq:diagram_one_row}
\xymatrix@C+5mm@R+2mm{
\mathscr{H}_{\ell_T,\mathcal{B}_T} \ar[r]^-{{\bf F}}  & \mathscr{H}_{C_k(C_k^{(-)}(\ell_{T'},\mathcal{B}_{T'}))} \ar[r]^-{{\bf R}} & \mathscr{H}_{C_k^{(-)}(\ell_{T'},\mathcal{B}_{T'})} \ar[r]^-{{\bf F}} & \mathscr{H}_{\ell_T,\mathcal{B}_T} 
}
\end{align}
The composition of the above three arrows equals ${{\bf K}'}_{\hspace{-1mm}k}^{(-)} \, \left({{\bf K}'}_{\hspace{-1mm}k}^{(+)} \right)^{-1} : \mathscr{H}_{\ell_T,\mathcal{B}_T} \to \mathscr{H}_{\ell_T,\mathcal{B}_T}$. Applying Lem.\ref{lem:compatibility_between_F_and_R} to the composition ${\bf F} {\bf R}$ of the latter two maps of eq.\eqref{eq:diagram_one_row} to re-write it as ${\bf R} {\bf F}$, we see that the composition of the diagram in eq.\eqref{eq:diagram_one_row} equals the composition of the following diagram:
\begin{align}
\label{eq:diagram_one_row2}
\xymatrix@C+5mm@R+2mm{
\mathscr{H}_{\ell_T,\mathcal{B}_T} \ar[r]^-{{\bf F}}  & \mathscr{H}_{C_k(C_k^{(-)}(\ell_{T'},\mathcal{B}_{T'}))} 
\ar[r]^-{{\bf F}} & \mathscr{H}_{C_k(\ell_T,\mathcal{B}_T)} \ar[r]^-{{\bf R}} & \mathscr{H}_{\ell_T,\mathcal{B}_T}
}
\end{align}
Then, applying Cor.\ref{cor:relationship_between_two_compositions_of_F} (or eq.\eqref{eq:F_composition_combine}) to the composition ${\bf F}{\bf F}$ of the first two maps of this new diagram eq.\eqref{eq:diagram_one_row2} to re-write it as ${\bf F}$ we see that the composition of the diagram in eq.\eqref{eq:diagram_one_row2} equals the composition of the following diagram, up to a multiplicative constant that is a power of $e^{{\rm i}\pi /4}$:
\begin{align}
\label{eq:diagram_one_row3}
\xymatrix@C+5mm@R+2mm{
\mathscr{H}_{\ell_T,\mathcal{B}_T} \ar[r]^-{{\bf F}} & \mathscr{H}_{C_k(\ell_T,\mathcal{B}_T)} \ar[r]^-{{\bf R}} & \mathscr{H}_{\ell_T,\mathcal{B}_T}
}
\end{align}
To summarize, we proved
$$
{{\bf K}'}_{\hspace{-1mm}k}^{(-)} \, \left({{\bf K}'}_{\hspace{-1mm}k}^{(+)} \right)^{-1} ~\sim~ {\bf R}_{C_k(\ell_T,\mathcal{B}_T), (\ell_T,\mathcal{B}_T)} \circ {\bf F}_{(\ell_T,\mathcal{B}_T), C_k(\ell_T,\mathcal{B}_T)}
$$

\vs

\ul{\it A special case.} Still, a direct computation of the map ${{\bf K}'}_{\hspace{-1mm}k}^{(-)} \, \left({{\bf K}'}_{\hspace{-1mm}k}^{(+)} \right)^{-1} : \mathscr{H}_{\ell_T,\mathcal{B}_T} \to \mathscr{H}_{\ell_T,\mathcal{B}_T}$ is not so straightforward, hence we take a detour. We first deal with an easy case; assume that the condition
$$
x_k \in \ell_T
$$
holds. For this special case, let's compute the composition of the diagram in eq.\eqref{eq:diagram_one_row3}. One observation is that $x_k \in \ell_T$ implies
$$
C_k(\ell_T) = \ell_T.
$$ 
Indeed, for each $x \in \ell_T$, one has $B(x,x_k)=0$ because $\ell_T$ is isotropic, hence $C_k(x) = x + B(x,x_k) x_k = x$; so $C_k$ restricts to the identity map $\ell_T \to \ell_T$. By means of the maps ${\bf I}$, consider the diagram:
$$
\xymatrix@C+5mm@R+0mm{
\mathscr{H}_{\ell_T,\mathcal{B}_T} \ar[r]^-{{\bf F}} & \mathscr{H}_{C_k(\ell_T), C_k(\mathcal{B}_T)} \ar[r]^-{{\bf R}} & \mathscr{H}_{\ell_T,\mathcal{B}_T} \\
\mathscr{H}_{\ell_T} \ar[u]_{{\bf I}} \ar[r]^-{{\bf F}} & \mathscr{H}_{C_k(\ell_T)} \ar[u]_{{\bf I}} \ar[r]^-{\til{\bf R}} & \mathscr{H}_{\ell_T} \ar[u]_{{\bf I}} 
}
$$
The upper row is eq.\eqref{eq:diagram_one_row3}. The left square commutes by the definition of the Weil intertwiner ${\bf F}$ associated to symplectic decompositions. The right square commutes by the definition of the map $\til{\bf R}$ as appeared in Lem.\ref{lem:til_R}. From $C_k(\ell_T) = \ell_T$ and eq.\eqref{eq:F_ell_ell}, it follows that the map ${\bf F} : \mathscr{H}_{\ell_T} \to \mathscr{H}_{C_k(\ell_T)}$ in the lower row is the identity map. 

\vs

Let $\varphi = \varphi(t_1,\ldots,t_r) \in \mathscr{H}_{\ell_T,\mathcal{B}_T} = L^2(\mathbb{R}^r, dt_1 \cdots dt_r)$ (the upper-left corner of the above diagram). Then ${\bf I}^{-1} \varphi \in \mathscr{H}_{\ell_T}$ is given by formula eq.\eqref{eq:I_ell_B_formula}, as a function $N \to \mathbb{C}$. Next, ${\bf F} {\bf I}^{-1} \varphi \in \mathscr{H}_{C_k(\ell_T)} = \mathscr{H}_{\ell_T}$ is identical to ${\bf I}^{-1} \varphi$, as functions $N \to \mathbb{C}$. Now let's compute $\til{\bf R}{\bf F}{\bf I}^{-1} \varphi \in \mathscr{H}_{\ell_T}$; as a function on $N$, it is
\begin{align*}
(\til{\bf R}{\bf F}{\bf I}^{-1} \varphi)(n) ~\underset{{\rm eq}.\eqref{eq:til_R_formula}}{=}~ ({\bf F}{\bf I}^{-1} \varphi)(e^{C_k}(n)) = ({\bf I}^{-1}\varphi)(e^{C_k}(n)), \quad \forall n\in N.
\end{align*}
To compute the final result ${\bf I} \til{\bf R}{\bf F}{\bf I}^{-1} \varphi \in \mathscr{H}_{\ell_T,\mathcal{B}_T}$ with the help of eq.\eqref{eq:I_ell_B_formula}, we put $n = \exp(\ssum_{j=1}^r t_j v_j) \cdot h$ for $t_1,\ldots,t_r \in \mathbb{R}$ and $h\in H$. Note $h\in H = \exp \mathfrak{h} = \exp(V^\perp + \ell_T + \mathbb{R}c)$, and observe that $C_k$ restricts to the identity map on $V^\perp + \ell_T + \mathbb{R}c$, because $B(x,x_k)=0$ for all $x\in V^\perp + \ell_T + \mathbb{R}c$. Hence we get $e^{C_k}(h) = h$, so $e^{C_k}(n)=e^{C_k}( \exp(\ssum_{j=1}^r t_j v_j) \cdot h) = \exp(\ssum_{j=1}^r t_j C_k(v_j) )\cdot h$. Let us put
$$
\alpha_j = B(v_j, x_k) \in \mathbb{R}, \qquad j =1,2,\ldots,r,
$$
so that
$$
C_k(v_j) = v_j + B(v_j,x_k) x_k = v_j + \alpha_j x_k, \qquad \forall j=1,\ldots,r.
$$
Thus
\begin{align}
\nonumber
\exp( \ssum_{j=1}^r t_j C_k(v_j) ) & = \exp(\ssum_{j=1}^r t_j v_j + (\ssum_{i=1}^r t_i \alpha_i) x_k) \\
\nonumber
& = \exp(\ssum_{j=1}^r t_j v_j) \exp( (\ssum_{i=1}^r t_i \alpha_i) x_k) \exp( {\textstyle -\frac{1}{2}} [ \ssum_{j=1}^r t_j v_j, \, (\ssum_{i=1}^r t_i \alpha_i) x_k] ) \\
\nonumber
& = \exp(\ssum_{j=1}^r t_j v_j) \exp( (\ssum_{i=1}^r t_i \alpha_i) x_k) \exp( {\textstyle -\frac{1}{2}} \ssum_{j=1}^r t_j (\ssum_{i=1}^r t_i \alpha_i) B(v_j,x_k) c) \\
\nonumber
& = \exp(\ssum_{j=1}^r t_j v_j) \exp( (\ssum_{i=1}^r t_i \alpha_i) x_k) \exp( {\textstyle -\frac{1}{2}} \ssum_{j=1}^r t_j (\ssum_{i=1}^r t_i \alpha_i) \alpha_j c) \\
\label{eq:square_computation1}
& = \exp(\ssum_{j=1}^r t_j v_j) \underbrace{ \exp( (\ssum_{i=1}^r t_i \alpha_i) x_k) \exp( {\textstyle -\frac{1}{2}} (\ssum_{i=1}^r t_j \alpha_j)^2 c)},
\end{align}
where the underbraced part is in $H = \exp \mathfrak{h}$. Assembling these, for $n = \exp(\ssum_{j=1}^r t_j v_j) \cdot h$ observe
\begin{align*}
& (\til{\bf R}{\bf F}{\bf I}^{-1} \varphi)(\exp(\ssum_{j=1}^r t_j v_j) \cdot h) \\
& = (\til{\bf R}{\bf F}{\bf I}^{-1} \varphi)(n) = ({\bf I}^{-1}\varphi)(e^{C_k}(n)) \\
& = ({\bf I}^{-1}\varphi)(\exp(\ssum_{j=1}^r t_j C_k(v_j) )\cdot h) \\
& = ({\bf I}^{-1}\varphi)( \exp(\ssum_{j=1}^r t_j v_j) \exp( (\ssum_{i=1}^r t_i \alpha_i) x_k) \exp( {\textstyle -\frac{1}{2}} (\ssum_{i=1}^r t_j \alpha_j)^2 c) \cdot h) \qquad (\because \mbox{eq.\eqref{eq:square_computation1}}) \\
& = \chi\left( \exp( (\ssum_{i=1}^r t_i \alpha_i) x_k) \exp( {\textstyle -\frac{1}{2}} (\ssum_{i=1}^r t_j \alpha_j)^2 c) \cdot h\right)^{-1} \cdot \varphi(t_1,\ldots,t_r) \qquad (\because \mbox{eq.\eqref{eq:I_ell_B_formula}}) \\
& = e^{ -{\rm i} ({\textstyle \frac{1}{2}} (\ssum_{i=1}^r t_j \alpha_j)^2 )} \cdot \chi(h)^{-1} \cdot \varphi(t_1,\ldots,t_r) \qquad (\because\mbox{$(\ssum_{i=1}^r t_i \alpha_i) x_k \in \ell_T$, and from definition of $\chi$}).
\end{align*}
So, in view of eq.\eqref{eq:I_ell_B_formula}, the element ${\bf I} \til{\bf R}{\bf F}{\bf I}^{-1} \varphi \in \mathscr{H}_{\ell_T,\mathcal{B}_T} = L^2(\mathbb{R}^r, dt_1\cdots dt_r)$ is given by
$$
({\bf I} \til{\bf R}{\bf F}{\bf I}^{-1} \varphi)(t_1,\ldots,t_r) = e^{ -{\rm i} ({\textstyle \frac{1}{2}} (\ssum_{i=1}^r t_j \alpha_j)^2 )} \cdot \varphi(t_1,\ldots,t_r).
$$
Hence, the unitary operator ${\bf I} \til{\bf R}{\bf F}{\bf I}^{-1}$ is the multiplication operator
$$
{\bf I} \til{\bf R}{\bf F}{\bf I}^{-1} = e^{ -{\rm i} ({\textstyle \frac{1}{2}} (\ssum_{i=1}^r t_j \alpha_j)^2 )}.
$$
Meanwhile, recall from \S\ref{subsec:irreducible_quantum_representation} that the self-adjoint operator $\wh{x}^\hbar_k$ is 
$$
\wh{x}^\hbar_k = \wh{x}^\hbar_{k;\ell_T,\mathcal{B}_T} = \pi^\hbar_{\ell_T,\mathcal{B}_T}(x_k) = {\rm i} \sqrt{2\pi\hbar} \cdot \pi_{\ell_T,\mathcal{B}_T}(x_k).
$$
Since $x_k \in \ell_T$, Example \ref{ex:operator_in_ell} applies, so that $\pi_{\ell_T,\mathcal{B}_T}(x_k)$ is the multiplication by ${\rm i} \ssum_{j=1}^r t_j \alpha_j$. Thus $\wh{x}^\hbar_k$ is the following multiplication operator
$$
\wh{x}^\hbar_k = - \sqrt{2\pi\hbar} \cdot \ssum_{j=1}^r t_j \alpha_j
$$
on its domain of self-adjointness. Therefore, the unitary operator $e^{(\wh{x}^\hbar_k)^2/(4\pi {\rm i} \hbar)}$ defined by means of the functional calculus for $\wh{x}^\hbar_k$ is the multiplication operator
$$
e^{(\wh{x}^\hbar_k)^2/(4\pi {\rm i} \hbar)} = e^{ 2\pi \hbar (\ssum_{j=1}^r t_j\alpha_j)^2 / (4\pi{\rm i}\hbar)}
= e^{ - \frac{\rm i}{2} (\ssum_{j=1}^r t_j\alpha_j)^2}
$$
hence precisely coincides with ${\bf I} \til{\bf R}{\bf F}{\bf I}^{-1}$.

\vs

So, in this special case $x_k \in \ell_T$, we showed that ${{\bf K}'}_{\hspace{-1mm}k}^{(-)} \, \left({{\bf K}'}_{\hspace{-1mm}k}^{(+)} \right)^{-1} = {\bf I} \til{\bf R}{\bf F}{\bf I}^{-1} : \mathscr{H}_{\ell_T,\mathcal{B}_T} \to \mathscr{H}_{\ell_T,\mathcal{B}_T}$ equals $e^{(\wh{x}^\hbar_k)^2/(4\pi {\rm i} \hbar)}$, hence also equals $\alpha_\hbar^2 \, \left({\bf K}_k^{\sharp \hbar (-)} \right)^{-1} \, {\bf K}_k^{\sharp \hbar (+)}$, up to $e^{{\rm i}\pi/4}$, as desired in this subsection.

\vs

\ul{\it A general case.} We now show that a general case boils down to the above special case. Let $(\ell_T,\mathcal{B}_T)$ and $(\ol{\ell}_T,\ol{\mathcal{B}}_T)$ be any two symplectic decompositions for $T$. We compare the composition of operators in eq.\eqref{eq:diagram_one_row3}, for these two choices. Consider the diagram
$$
\xymatrix@C+5mm@R+2mm{
\mathscr{H}_{\ell_T,\mathcal{B}_T} \ar[r]^-{{\bf F}} \ar[d]^{{\bf F}} & \mathscr{H}_{C_k(\ell_T,\mathcal{B}_T)} \ar[r]^-{{\bf R}} \ar[d]^{{\bf F}} & \mathscr{H}_{\ell_T,\mathcal{B}_T} \ar[d]^{{\bf F}} \\
\mathscr{H}_{\ol{\ell}_T, \ol{\mathcal{B}}_T} \ar[r]^-{{\bf F}} & \mathscr{H}_{C_k(\ol{\ell}_{T},\ol{\mathcal{B}}_{T})} \ar[r]^-{{\bf R}}  & \mathscr{H}_{\ol{\ell}_{T}, \ol{\mathcal{B}}_{T}}
}
$$
where the upper row is the diagram in eq.\eqref{eq:diagram_one_row3} for the former choice, the lower row is that for the latter choice, and the vertical arrows are the Weil intertwiners. In this diagram, the left square is commutative up to a power of $e^{{\rm i}\pi/4}$, by Cor.\ref{cor:relationship_between_two_compositions_of_F}. The commutativity of the right square follows from Lem.\ref{lem:compatibility_between_F_and_R}.

\vs

Then we can observe that the sought-for equality
$$
{{\bf K}'}_{\hspace{-1mm}k}^{(-)} \, \left({{\bf K}'}_{\hspace{-1mm}k}^{(+)} \right)^{-1} ~\sim~ e^{(\wh{x}^\hbar_k)^2/(4\pi {\rm i} \hbar)} 
$$
holds for the first choice $(\ell_T,\mathcal{B}_T)$ if and only if it holds for the second choice $(\ol{\ell}_T, \ol{\mathcal{B}}_T)$. Indeed, in view of the above diagram we just looked at, the left hand side ${{\bf K}'}_{\hspace{-1mm}k}^{(-)} \, \left({{\bf K}'}_{\hspace{-1mm}k}^{(+)} \right)^{-1}$ for the first choice equals that for the second choice conjugated by ${\bf F} : \mathscr{H}_{\ell_T,\mathcal{B}_T} \to \mathscr{H}_{\ol{\ell}_T,\ol{\mathcal{B}}_T}$, up to a power of $e^{{\rm i} \pi /4}$. In the meantime, one can easily observe from the definition and the construction of $\wh{x}_k$ that the right hand side $e^{(\wh{x}_k)^2/(4\pi {\rm i} \hbar)}$ for the first batch equals that for the second batch conjugated by the same ${\bf F}$, using the fact that ${\bf F}$ is an intertwiner for representations; a similar argument as in the proof of Prop.\ref{prop:compatibility_of_mutation_intertwiner_and_Weil_intertwiners} would work.

\vs

So, to prove this equality, it suffices to prove it for any special choice of symplectic decomposition $(\ell_T,\mathcal{B}_T)$ for $T$. It is indeed possible to choose $(\ell_T,\mathcal{B}_T)$ for $T$ so that $x_k \in \ell_T$ holds, by a simple linear algebra. Then the situation boils down to the special case which we already settled.

\vs

To summarize, we showed that ${\bf K}_k^{\hbar(+)} = {\bf K}_k^{\hbar(-)}$ holds up to a power of $e^{{\rm i} \pi /4}$, establishing Prop.\ref{prop:equality_of_two_signed_decompositions}. We leave to readers to check whether this equality holds exactly without a multiplicative constant, under a certain assumption if necessary.

\section{The phase constants for the compositions of the mutation intertwiners}
\label{sec:phase_constants_for_composition}

\subsection{The consistency relations for the flips}
\label{subsect:consistency_relations}

A \ul{\bf change of ideal triangulations} of a punctured oriented surface $S$ is a pair $[T,T']$ of ideal triangulations, sometimes denoted $T\leadsto T'$. We think of {\em applying} this change $[T,T']$ to $T$ to obtain $T'$, and write $[T,T'](T) = T'$. They form the following groupoid.
\begin{definition}[\cite{P87} \cite{P12}]
A \ul{\bf labeled (ideal) triangulation} $T$ of a punctured surface $S$ consists of an ideal triangulation of $S$ and the labeling of edges of $T$ by a fixed index set $I$.

The \ul{\bf (labeled) Ptolemy groupoid} of a punctured surface is the groupoid ${\rm Pt}(S)$ whose set of objects is the set of all labeled triangulations, and from any object $T$ to any object $T'$, there is exactly one morphism, denoted by $[T,T']$.
\end{definition}
Two morphisms $[T,T']$ and $[T'',T''']$ can be composed if $T' = T''$, and we write
$$
[T',T'''] \circ [T,T'] = [T,T'''],
$$
or sometimes without the composition symbol $\circ$. The Ptolemy groupoid is an example of the cluster modular groupoids generated by the mutations, in the theory of the cluster varieties and the cluster algebras \cite{FG09}. 
\begin{definition}[\cite{P87} \cite{P12} \cite{K16b}]
A morphism $[T,T']$ of ${\rm Pt}(S)$ is called a \ul{\bf flip at $k$} (for $k\in I$) and represented by the symbol $\mu_k$, if $T$ coincides with $T'$ except at one edge labeled by $k$. 

A morphism $[T,T']$ of ${\rm Pt}(S)$ is called a \ul{\bf label-change by $\sigma$} (for a permutation $\sigma$ of $I$) and represented by the symbol $P_\sigma$, if $T$ and $T'$ have same underlying triangulations, and for each $i \in I$, the edge labeled by $i$ in $T$ is labeled by $\sigma(i)$ in $T'$.

These morphisms are called the \ul{\bf elementary morphisms} of ${\rm Pt}(S)$.
\end{definition}
A flip is an example of a {\em mutation} in the theory of the cluster varieties. Each symbol $\mu_k$ and $P_\sigma$ represents many different changes of triangulations. A classical result is:
\begin{proposition}[\cite{FLP} {\cite[Fact 1.24]{P12}}]
Any two ideal triangulations $T$ and $T'$ are connected by a finite sequence of flips and label-changes; elementary morphisms of ${\rm Pt}(S)$ generate all morphisms of ${\rm Pt}(S)$.
\end{proposition}
By a sequence of flips and label-changes, we mean the corresponding composition of the members of the sequence. Also classically known is the complete list of relations satisfied by the flips and the label-changes. 
\begin{proposition}[{\cite[Thm.1.1]{CP}} {\cite[Thm.3.10]{FST}} {\cite[Prop.2.30]{K16b}}: the presentation of the Ptolemy groupoid]
\label{prop:presentation_of_Pt}
Any algebraic relations satisfied by the flips and the label-changes are consequences of the following:
\begin{align*}
\begin{array}{rl}
\mbox{\rm (the twice-flip identity)} & \mu_k \mu_k = {\rm Id} \quad (\mbox{when applied to any $T$}) \\ 
\mbox{\rm (the square identity)} & \mu_j \mu_i \mu_j \mu_i = {\rm Id} \quad (\mbox{when applied to any $T$ with $\varepsilon_{ij}=0$}) \\
\mbox{\rm (the pentagon identity)} & P_{(ij)} \mu_i \mu_j \mu_i \mu_j \mu_i  =  {\rm Id} \quad (\mbox{when applied to any $T$ with $\varepsilon_{ij}\in \{1,-1\}$}) \\
\mbox{\rm (the permutation identities)} & P_\sigma \, P_\gamma = P_{\sigma \circ \gamma}, \qquad P_\sigma \, \mu_k \, P_{\sigma^{-1}} = \mu_{\sigma(k)}, \qquad P_{\rm Id} = {\rm Id}.  \quad(\mbox{when applied to any $T$})
\end{array}
\end{align*}
\end{proposition}

We shall check that the composition identities for the mutation intertwiners corresponding to the above classical relations hold up to constants, and compute these constants. A convenient way to formulate the problem and the result is to use the notion of a projective representation of a groupoid. For the present paper, we first have to extend the Ptolemy groupoid:
\begin{definition}
The \ul{\bf symplectic Ptolemy groupoid} of a punctured surface $S$ is the groupoid ${\rm SPt}(S)$ defined as the following category. The set of objects is the collection of all triples $(T,\ell_T,\mathcal{B}_T)$, where $T$ is a labeled ideal triangulation of $S$, and $(\ell_T,\mathcal{B}_T)$ is a symplectic decomposition of $(V_T,B_T)$ in the sense as defined in Def.\ref{def:symplectic_decomposition}. From any object $(T,\ell_T,\mathcal{B}_T)$ to any object $(T',\ell_{T'},\mathcal{B}_{T'})$, there is exactly one morphism, denoted by $[(T,\ell_T,\mathcal{B}_T), (T',\ell_{T'},\mathcal{B}_{T'})]$
\end{definition}
There is a forgetful functor ${\rm SPt}(S) \to {\rm Pt}(S)$; the notion of the elementary morphisms of ${\rm Pt}(S)$ carries over, with an additional one.

\begin{definition}
A morphism $[(T,\ell_T,\mathcal{B}_T), (T',\ell_{T'},\mathcal{B}_{T'})]$ of ${\rm SPt}(S)$ is called a \ul{\bf flip at $k \in I$} and denoted by $\mu_k$ if $[T,T']$ is $\mu_k$ of ${\rm Pt}(S)$, and then we write $[(T,\ell_T,\mathcal{B}_T), (T',\ell_{T'},\mathcal{B}_{T'})] = \mu_k$.

A morphism $[(T,\ell_T,\mathcal{B}_T), (T',\ell_{T'},\mathcal{B}_{T'})]$ of ${\rm SPt}(S)$ is called a \ul{\bf label-change} by $\sigma$ (for a permutation $\sigma$ of $I$) and denoted by $P_\sigma$ if $[T,T']$ is $P_\sigma$ of ${\rm Pt}(S)$, and then we write $[(T,\ell_T,\mathcal{B}_T), (T',\ell_{T'},\mathcal{B}_{T'})] = P_\sigma$.

A morphism $[(T,\ell_T,\mathcal{B}_T), (T, \ol{\ell}_T, \ol{\mathcal{B}}_T)]$ of ${\rm SPt}(S)$ is called a \ul{\bf decomposition-change} and denoted by $F$ or $F_{(\ell_T,\mathcal{B}_T),(\ol{\ell}_T,\ol{\mathcal{B}}_T)}$, and then we write $[(T,\ell_T,\mathcal{B}_T), (T',\ell_{T'},\mathcal{B}_{T'})] = F_{(\ell_T,\mathcal{B}_T),(\ol{\ell}_T,\ol{\mathcal{B}}_T)}$.

These are called the \ul{\bf elementary morphisms} of ${\rm SPt}(S)$.
\end{definition}

\begin{proposition}[the presentation of the symplectic Ptolemy groupoid]
\label{prop:presentation_of_SPt}
One has:
\begin{enumerate}
\item[\rm (1)] The elementary morphisms of ${\rm SPt}(S)$ generate all morphisms of ${\rm SPt}(S)$.

\item[\rm (2)] Any relation among the elementary morphisms of ${\rm SPt}(S)$ are consequences of the ones in Prop.\ref{prop:presentation_of_Pt} together with
\begin{align*}
\begin{array}{rl}
\mbox{\rm (the $FF=F$ identity)} & F_{(\ell_T',\mathcal{B}_T'),(\ell_T'',\mathcal{B}_T'')} \circ F_{(\ell_T,\mathcal{B}_T),(\ell_T',\mathcal{B}_T')} = F_{(\ell_T,\mathcal{B}_T),(\ell_T'',\mathcal{B}_T'')}, \\
& (\mbox{$\forall$three symplectic decompositions $(\ell_T,\mathcal{B}_T),(\ell_T',\mathcal{B}_T'),(\ell_T'',\mathcal{B}_T'')$ for $T$}) \\ 
\mbox{\rm (the compatibility)} & \mbox{the following two diagrams commute} 
\\
& \xymatrix{
(T,\ell_T,\mathcal{B}_T) \ar[r]^-{\mu_k} \ar[d]_{F} & (T',\ell_{T'},\mathcal{B}_{T'}) \ar[d]^{F} \\
(T,\ol{\ell}_T, \ol{\mathcal{B}}_T) \ar[r]^-{\mu_k} & (T',\ol{\ell}_{T'},\ol{\mathcal{B}}_{T'}) \\
}
\qquad
\xymatrix{
(T,\ell_T,\mathcal{B}_T) \ar[r]^-{P_\sigma} \ar[d]_{F} & (T',\ell_{T'},\mathcal{B}_{T'}) \ar[d]^{F} \\
(T,\ol{\ell}_T, \ol{\mathcal{B}}_T) \ar[r]^-{P_\sigma} & (T',\ol{\ell}_{T'},\ol{\mathcal{B}}_{T'}) \\
}
\end{array}
\end{align*}
\end{enumerate}
\end{proposition}
The part (1) is easy to see. Verifying the relations in the part (2) is also easy. To show that these relations generate all relations, one can mimick the proof in \cite[\S4]{K16b}, which uses a lemma of Bakalov-Kirillov \cite{BK}. We leave it as a straightforward exercise to readers.

\begin{remark}
\label{rem:care_on_P}
Some care is need when transplanting the relations of Prop.\ref{prop:presentation_of_Pt} to the groupoid ${\rm SPt}(S)$, due to our definition of the label-change symbol $P_\sigma$. For example, $P_{\rm Id} = {\rm Id}$ has to be modified. If the morphism $[(T,\ell_T,\mathcal{B}_T),(T',\ell_{T'},\mathcal{B}_{T'})]$ of ${\rm Spt}(S)$ is $P_{\rm Id}$, then $T = T'$, hence this morphism is also $F_{(\ell_T,\mathcal{B}_T),(\ell_{T'},\mathcal{B}_{T'})}$, but not necessarily the identity. So, $P_{\rm Id} = {\rm Id}$ should be replaced by  $P_{\rm Id} = F$ in a suitable sense.
\end{remark}

\subsection{The projective representations of the symplectic Ptolemy groupoid}
\label{subsec:projective_representation_of_SPt_S}

Our intertwiners for irreducible self-adjoint representations of the quantum Teichm\"uller space can be viewed as forming a contravariant {\em projective} functor
$$
\pi = \pi^\hbar_\lambda : {\rm SPt}(S) \to {\rm Hilb},
$$
depending on the real parameter $\hbar$ and a `puncture function' $\lambda$ chosen as in eq.\eqref{eq:lambda}, where ${\rm Hilb}$ is the category whose objects are complex Hilbert spaces and whose morphisms are unitary maps. This functor being projective means that the composition of the morphisms goes to the composition of the morphisms up to a multiplicative constant of modulus $1$. Note that the above functor only captures the information of the intertwining operators. To keep track of the actual representations of the quantum Teichm\"uller space that are being intertwined, one has to consider a more complicated category than ${\rm Hilb}$, namely the category of representations of the Chekhov-Fock algebras; we do not do so in the present paper, and refer the readers to \S5.4 of the version 2 of the paper \cite{K16c} (i.e. arXiv:1602.00797v2).

\vs

To each object $(T,\ell_T,\mathcal{B}_T)$ of ${\rm SPt}(S)$, we associate the Hilbert space $\pi(T,\ell_T,\mathcal{B}_T) := \mathscr{H}_{\ell_T,\mathcal{B}_T}$ we constructed in \S\ref{subsec:irreducible_quantum_representation}, keeping in mind its structure as representations of $N_T$, $\mathfrak{n}_T$, and $L^q_T$ which we studied. 

\vs

\ul{\it The flip operator.} To each elementary morphism $[(T,\ell_T,\mathcal{B}_T), (T',\ell_{T'},\mathcal{B}_{T'})]$ of ${\rm SPt}(S)$ that is a flip $\mu_k$, we associate the operator ${\bf K}_k^{\hbar(\epsilon)} : \mathscr{H}_{\ell_{T'},\mathcal{B}_{T'}} \to \mathscr{H}_{\ell_T,\mathcal{B}_T}$. We choose either of the signs $\epsilon \in \{+,-\}$; the ambiguity coming from this choice of a sign is the multiplication by a power of $e^{{\rm i}\pi/4}$, as seen before. 

\vs

\ul{\it The permutation operator.} To each elementary morphism $[(T,\ell_T,\mathcal{B}_T), (T',\ell_{T'},\mathcal{B}_{T'})]$ of ${\rm SPt}(S)$ that is a label-change $P_\sigma$, we associate the operator
$$
{\bf P}_\sigma = {\bf P}_{\sigma;(\ell_T,\mathcal{B}_T),(\ell_{T'},\mathcal{B}_{T'})} : \mathscr{H}_{\ell_{T'},\mathcal{B}_{T'}} \to \mathscr{H}_{\ell_T,\mathcal{B}_T}
$$
defined as follows. First, note that $T$ and $T'$ have the same underlying ideal triangulation; define the bijection $P_\sigma : T \to T'$ induced by the label change $P_\sigma$. That is, $e\in T$ and $P_\sigma(e) \in T'$ are labeled by a same element of the index set $I$, through the respective labelings of $T$ and $T'$ (if $e\in T$ is labeled by $i\in I$ by the labeling of $T$, then $P_\sigma(e) \in T'$ is labeled by $\sigma(i) \in I$ in $T'$). Consider the induced isomorphism
\begin{align}
\label{eq:C_sigma}
C_\sigma : V_{T'} \to V_T,
\end{align}
sending each basis vector $x_{P_\sigma(e)}$ to $x_e$. By construction, $C_\sigma$ respects the forms $B_{T'}$ and $B_T$, thus $C_\sigma(\ell_{T'},\mathcal{B}_{T'})$ is a symplectic decomposition of $(V_T,B_T)$. Following the philosophy of the construction in eq.\eqref{eq:monomial_transformation_part} of the monomial transformation part ${\bf K}'$ of the flip operator, we define
\begin{align}
\label{eq:permutation_operator_definition}
{\bf P}_{\sigma;(\ell_T,\mathcal{B}_T),(\ell_{T'},\mathcal{B}_{T'})} := {\bf F}_{C_\sigma(\ell_{T'},\mathcal{B}_{T'}), \, (\ell_T,\mathcal{B}_T)} \circ {\bf R}_{(\ell_{T'},\mathcal{B}_{T'}), \, C_\sigma(\ell_{T'},\mathcal{B}_{T'})}
\end{align}

\vs

\ul{\it The Weil interwining operator.} To each elementary morphism $[(T,\ell_T,\mathcal{B}_T), (T,\ol{\ell}_T,\ol{\mathcal{B}}_T)]$ of ${\rm SPt}(S)$ that is a decomposition change $F_{(\ell_T,\mathcal{B}_T),(\ol{\ell}_T,\ol{\mathcal{B}}_T)}$, we associate the operator
$$
{\bf F} = {\bf F}_{(\ol{\ell}_T,\ol{\mathcal{B}}_T), (\ell_T,\mathcal{B}_T)} : \mathscr{H}_{\ol{\ell}_T,\ol{\mathcal{B}}_T} \to \mathscr{H}_{\ell_T,\mathcal{B}_T}
$$
i.e. the Weil intertwiner. 

\vs

\ul{\it A general intertwining operator.} Finally, to {\em any} morphism $[(T,\ell_T,\mathcal{B}_T), (T',\ell_{T'},\mathcal{B}_{T'})]$ of ${\rm SPt}(S)$, first write it as (the composition of) a sequence of elementary morphisms, and associate the composition of the corresponding operators ${\bf K}_k^{\hbar(\epsilon)}$, ${\bf P}_\sigma$, and ${\bf F}$. A morphism can be expressed as different sequences of elementary morphisms; to guarantee that the operators associated to these sequences are well-defined up to a constant, it suffices to show that the operators representing the elementary morphisms satisfy the operator identities corresponding to the generating relations of the elementary morphisms of ${\rm SPt}(S)$, written in Prop.\ref{prop:presentation_of_Pt} and Prop.\ref{prop:presentation_of_SPt}, up to constants.

\vs

More precisely, we shall prove the following, which can be viewed as the major part of the \ul{\bf main result} of the present paper.
\begin{theorem}[the operator identities for the operators representing the elementary morphisms]
\label{thm:main}
The operators ${\bf K}_k^{\hbar(\epsilon)}$, ${\bf P}_\sigma$, and ${\bf F}$ representing the elementary morphisms $\mu_k$, $P_\sigma$, and $F$ of the symplectic Ptolemy groupoid ${\rm SPt}(S)$ satisfy the following operator identities corresponding to the relations among the selementary morphisms of ${\rm SPt}(S)$.
\begin{enumerate}
\item[\rm (1)] (the twice-flip identity) Suppose
$$
[(T,\ell_T,\mathcal{B}_T), (T',\ell_{T'},\mathcal{B}_{T'})] = \mu_k
$$
Then $[(T',\ell_{T'},\mathcal{B}_{T'}), (T,\ell_T,\mathcal{B}_T)]=\mu_k$, and
$$
{\bf K}_k^{\hbar(+)} \, {\bf K}_k^{\hbar(-)} = {\bf K}^{\hbar(+)}_{(\ell_T,\mathcal{B}_T),(\ell_{T'},\mathcal{B}_{T'})} \, {\bf K}^{\hbar(-)}_{(\ell_{T'},\mathcal{B}_{T'}),(\ell_T,\mathcal{B}_T)} ~\sim~ {\rm Id}_{\mathscr{H}_{\ell_T,\mathcal{B}_T}}
$$
i.e. the equality holds up to an integer power of $e^{{\rm i} \pi /4}$.

\item[\rm (2)] (the square identity) Let $T^{(0)}$ be an ideal triangulation whose exchange matrix $\varepsilon^{(0)}$ satisfies
$$
\varepsilon^{(0)}_{ij} = 0
$$
for some $i,j \in I$ with $i\neq j$. Let $(T^{(a)}, \ell_{T^{(a)}}, \mathcal{B}_{T^{(a)}})$, $a=0,\ldots,3$, be objects of ${\rm SPt}(S)$ such that
\begin{align*}
& [(T^{(0)},\ell_{T^{(0)}},\mathcal{B}_{T^{(0)}}), (T^{(1)},\ell_{T^{(1)}},\mathcal{B}_{T^{(1)}})] = \mu_i, \quad
[(T^{(1)},\ell_{T^{(1)}},\mathcal{B}_{T^{(1)}}), (T^{(2)},\ell_{T^{(2)}},\mathcal{B}_{T^{(2)}})] = \mu_j, \\
& [(T^{(2)},\ell_{T^{(2)}},\mathcal{B}_{T^{(2)}}), (T^{(3)},\ell_{T^{(3)}},\mathcal{B}_{T^{(3)}})] = \mu_i.
\end{align*}
Then $[(T^{(3)},\ell_{T^{(3)}},\mathcal{B}_{T^{(3)}}), (T^{(0)},\ell_{T^{(0)}},\mathcal{B}_{T^{(0)}})] = \mu_j$, and
\begin{align*}
{\bf K}_i^{\hbar(+)} \, {\bf K}_j^{\hbar(+)} \, {\bf K}_i^{\hbar(-)} \, {\bf K}_i^{\hbar(-)}
& = {\bf K}_{(\ell_{T^{(0)}},\mathcal{B}_{T^{(0)}}),(\ell_{T^{(1)}},\mathcal{B}_{T^{(1)}})}^{\hbar(+)} \, {\bf K}_{(\ell_{T^{(1)}},\mathcal{B}_{T^{(1)}}),(\ell_{T^{(2)}},\mathcal{B}_{T^{(2)}})}^{\hbar(+)}  \\
& \quad \cdot {\bf K}_{(\ell_{T^{(2)}},\mathcal{B}_{T^{(2)}}),(\ell_{T^{(3)}},\mathcal{B}_{T^{(3)}})}^{\hbar(-)} \, {\bf K}_{(\ell_{T^{(3)}},\mathcal{B}_{T^{(3)}}),(\ell_{T^{(0)}},\mathcal{B}_{T^{(0)}})}^{\hbar(-)} \\
& \sim ~ {\rm Id}_{\mathscr{H}_{\ell_{T^{(0)}},\mathcal{B}_{T^{(0)}}}}
\end{align*}
i.e. the equality holds up to an integer power of $e^{{\rm i} \pi /4}$.

\item[\rm (3)] (the pentagon identity) Let $T^{(0)}$ be an ideal triangulation whose exchange matrix $\varepsilon^{(0)}$ satisfies
$$
\varepsilon^{(0)}_{ij} \in \{1,-1\}
$$
for some $i,j \in I$. Let $(T^{(a)}, \ell_{T^{(a)}}, \mathcal{B}_{T^{(a)}})$, $a=0,\ldots,5$, be objects of ${\rm SPt}(S)$ such that
\begin{align*}
& [(T^{(0)},\ell_{T^{(0)}},\mathcal{B}_{T^{(0)}}), (T^{(1)},\ell_{T^{(1)}},\mathcal{B}_{T^{(1)}})] = \mu_i, \quad
[(T^{(1)},\ell_{T^{(1)}},\mathcal{B}_{T^{(1)}}), (T^{(2)},\ell_{T^{(2)}},\mathcal{B}_{T^{(2)}})] = \mu_j, \\
& [(T^{(2)},\ell_{T^{(2)}},\mathcal{B}_{T^{(2)}}), (T^{(3)},\ell_{T^{(3)}},\mathcal{B}_{T^{(3)}})] = \mu_i, \quad
[(T^{(3)},\ell_{T^{(3)}},\mathcal{B}_{T^{(3)}}), (T^{(4)},\ell_{T^{(4)}},\mathcal{B}_{T^{(4)}})] = \mu_j, \\
& [(T^{(4)},\ell_{T^{(4)}},\mathcal{B}_{T^{(4)}}), (T^{(5)},\ell_{T^{(5)}},\mathcal{B}_{T^{(5)}})] = \mu_i.
\end{align*}
Then $[(T^{(5)},\ell_{T^{(5)}},\mathcal{B}_{T^{(5)}}), (T^{(0)},\ell_{T^{(0)}},\mathcal{B}_{T^{(0)}})] = P_{(ij)}$, and
\begin{align*}
& {\bf K}_i^{\hbar(+)} \, {\bf K}_j^{\hbar(+)} \, {\bf K}_i^{\hbar(-)} \, {\bf K}_i^{\hbar(-)} \, {\bf K}_i^{\hbar(-)} \, {\bf P}_{(ij)} \\
& = {\bf K}_{(\ell_{T^{(0)}},\mathcal{B}_{T^{(0)}}),(\ell_{T^{(1)}},\mathcal{B}_{T^{(1)}})}^{\hbar(+)} \, {\bf K}_{(\ell_{T^{(1)}},\mathcal{B}_{T^{(1)}}),(\ell_{T^{(2)}},\mathcal{B}_{T^{(2)}})}^{\hbar(+)} \, {\bf K}_{(\ell_{T^{(2)}},\mathcal{B}_{T^{(2)}}),(\ell_{T^{(3)}},\mathcal{B}_{T^{(3)}})}^{\hbar(-)} \\
& \quad \cdot {\bf K}_{(\ell_{T^{(3)}},\mathcal{B}_{T^{(3)}}),(\ell_{T^{(4)}},\mathcal{B}_{T^{(4)}})}^{\hbar(-)} \, {\bf K}_{(\ell_{T^{(4)}},\mathcal{B}_{T^{(4)}}),(\ell_{T^{(5)}},\mathcal{B}_{T^{(5)}})}^{\hbar(-)} \, {\bf P}_{(ij); (\ell_{T^{(5)}},\mathcal{B}_{T^{(5)}}),(\ell_{T^{(0)}},\mathcal{B}_{T^{(0)}})} \\
& \sim ~ \alpha_\hbar^{-1} \cdot {\rm Id}_{\mathscr{H}_{\ell_{T^{(0)}},\mathcal{B}_{T^{(0)}}}}
\end{align*}
i.e. the equality holds up to an integer power of $e^{{\rm i} \pi/4}$.

\item[\rm (4)] (the `trivial' relations) The operators ${\bf K}_k^{\hbar(\epsilon)}$, ${\bf P}_\sigma$, and ${\bf F}$ also satisfy the remaining relations:
\begin{align*}
& {\bf P}_\sigma {\bf P}_\gamma \sim {\bf P}_{\sigma\circ \gamma}, \qquad
{\bf P}_\sigma {\bf K}_k^{\hbar(\epsilon)} {\bf P}_{\sigma^{-1}} \sim {\bf K}_{\sigma(k)}^{\hbar(\epsilon)}, \qquad {\bf P}_{\rm Id} = {\rm Id}, \\
& {\bf FF} \sim {\bf F}, \qquad\qquad
{\bf F} \, {\bf K}_k^{\hbar(\epsilon)} \sim {\bf K}_k^{\hbar(\epsilon)} \, {\bf F}, \qquad\qquad
{\bf P}_\sigma \, {\bf F} \sim {\bf F} \, {\bf P}_\sigma,
\end{align*}
in appropriate senses.
\end{enumerate}
In particular, these operator induces a projective functor $\pi : {\rm SPt}(S) \to {\rm Hilb}$.
\end{theorem}

\begin{remark}
A similar remark as in Rem.\ref{rem:care_on_P} applies to the part (4) of the above theorem, especially for ${\bf P}_{\rm Id} = {\rm Id}$; it should really read as ${\bf P}_{\rm Id} = {\bf F}$.
\end{remark}

Due to Prop.\ref{prop:equality_of_two_signed_decompositions}, the (tropical) signs for the flips can in fact be chosen arbitrarily, with correspondingly modified phase constants. We will see in the proof of Thm.\ref{thm:main} why we chose the signs like above.

\subsection{A proof of the operator identities and the computation of the phase constants}

Here we prove Thm.\ref{thm:main}. We first deal with the hardest part, namely the part (3).

\vspace{2mm}

\ul{\it The pentagon identity.} Assume the hypothesis of Thm.\ref{thm:main}.(3) and further assume $\varepsilon^{(0)}_{ij}=1$. We first establish certain notations for the five flips, to ease the discussion. For each $a=1,\ldots,5$, consider the $a$-th flip $[(T^{(a-1)}, \ell_{T^{(a-1)}}, \mathcal{B}_{T^{(a-1)}}), (T^{(a)}, \ell_{T^{(a)}}, \mathcal{B}_{T^{(a)}})]$, which is a flip at the edge $k_a \in \{i,j\}$. It is wise to come up with a concise notation to denote the symplectic decompositions:
$$
\mathcal{D}^{(a)} := (\ell_{T^{(a)}},\mathcal{B}_{T^{(a)}}), \quad a=0,1,\ldots,5.
$$
We assigned the (tropical) signs $\epsilon_a$ to each of the five flips for $a=1,\ldots,5$, namely $+,+,-,-,-$. We denote the corresponding linear map
$$
C_{k_a}^{(\epsilon_a)} : V_{T^{(a)}} \to V_{T^{(a-1)}}
$$
in eq.\eqref{eq:C_k_epsilon} by the symbol $C^{(a)}$. Also, let us write
$$
\wh{x}^{(a)} := \wh{x}^\hbar_{k_a;\ell_{T^{(a-1)}}, \mathcal{B}_{T^{(a-1)}}} = \wh{x}^\hbar_{k_a;\mathcal{D}^{(a-1)}},
$$
which is a self-adjoint operator on $\mathscr{H}_{\ell_{T^{(a-1)}},\mathcal{B}_{T^{(a-1)}}} = \mathscr{H}_{\mathcal{D}^{(a-1)}}$. The exchange matrix for $T^{(a-1)}$ is denoted by $\varepsilon^{(a-1)}$. Then $\varepsilon^{(0)}_{ij} = - \varepsilon^{(1)}_{ij} = \varepsilon^{(2)}_{ij} = - \varepsilon^{(3)}_{ij} = \varepsilon^{(4)}_{ij} = 1$.

\vs

Now, in the composition
$$
{\bf K}_{\mathcal{D}^{(0)},\mathcal{D}^{(1)}}^{\hbar(+)} \, {\bf K}_{\mathcal{D}^{(1)},\mathcal{D}^{(2)}}^{\hbar(+)} \, {\bf K}_{\mathcal{D}^{(2)},\mathcal{D}^{(3)}}^{\hbar(-)}\, {\bf K}_{\mathcal{D}^{(3)},\mathcal{D}^{(4)}}^{\hbar(-)} \, {\bf K}_{\mathcal{D}^{(4)},\mathcal{D}^{(5)}}^{\hbar(-)} \, {\bf P}_{(ij); \mathcal{D}^{(5)},\mathcal{D}^{(0)} }
$$
we shall write each ${\bf K}^{\hbar(\epsilon_a)}_{k_a}$ as $\alpha_\hbar^{\epsilon_a} \, \Phi^\hbar(\epsilon_a \, \wh{x}^{(a)})^{\epsilon_a} \, {{\bf K}'}^{(\epsilon_a)}_{\hspace{-1mm} k_a}$, and move all ${{\bf K}'}^{(\epsilon_a)}_{\hspace{-1mm} k_a}$ to the right using eq.\eqref{eq:conjugation_of_primed_part_on_generators} of Lem.\ref{lem:conjugation_of_primed_part}, together with the commutativity of the conjugation by a unitary operator and the functional calculus (i.e. eq.\eqref{eq:conjugation_commutes_with_functional_calculus}). So, moving ${{\bf K}'}^{(\epsilon_a)}_{\hspace{-1mm} k_a}$ to the right results in altering the arguments of $\Phi^\hbar$ correspondingly. To show some steps, observe
\begin{align*}
& {\bf K}_{\mathcal{D}^{(0)},\mathcal{D}^{(1)}}^{\hbar(+)} \, {\bf K}_{\mathcal{D}^{(1)},\mathcal{D}^{(2)}}^{\hbar(+)} \, {\bf K}_{\mathcal{D}^{(2)},\mathcal{D}^{(3)}}^{\hbar(-)}\, {\bf K}_{\mathcal{D}^{(3)},\mathcal{D}^{(4)}}^{\hbar(-)} \, {\bf K}_{\mathcal{D}^{(4)},\mathcal{D}^{(5)}}^{\hbar(-)} \, {\bf P}_{(ij); \mathcal{D}^{(5)},\mathcal{D}^{(0)} } \\
& = \alpha_\hbar^{1+1-1-1-1} \, \Phi^\hbar(\wh{x}^{(1)}) \, {{\bf K}'}^{(+)}_{\hspace{-1mm}k_1} \, 
\Phi^\hbar(\wh{x}^{(2)}) \, {{\bf K}'}^{(+)}_{\hspace{-1mm}k_2} \, 
\Phi^\hbar(-\wh{x}^{(3)})^{-1} \, {{\bf K}'}^{(+)}_{\hspace{-1mm}k_3} \, 
\Phi^\hbar(-\wh{x}^{(4)})^{-1} \\
& \quad \cdot \underbrace{ {{\bf K}'}^{(+)}_{\hspace{-1mm}k_4} }_{\mbox{\tiny move to right}} \, \Phi^\hbar(-\wh{x}^{(5)})^{-1} \, {{\bf K}'}^{(+)}_{\hspace{-1mm}k_5} \, {\bf P}_{(ij); \mathcal{D}^{(5)},\mathcal{D}^{(0)} } \\
& = \alpha_\hbar^{-1} \, \Phi^\hbar(\wh{x}^{(1)}) \, {{\bf K}'}^{(+)}_{\hspace{-1mm}k_1} \, 
\Phi^\hbar(\wh{x}^{(2)}) \, {{\bf K}'}^{(+)}_{\hspace{-1mm}k_2} \, 
\Phi^\hbar(-\wh{x}^{(3)})^{-1} \, \underbrace{ {{\bf K}'}^{(+)}_{\hspace{-1mm}k_3} }_{\mbox{\tiny move to right}} \, 
\Phi^\hbar(-\wh{x}^{(4)})^{-1} \\
& \quad \cdot   \, \Phi^\hbar(C^{(4)}(-\wh{x}^{(5)}))^{-1} \, {{\bf K}'}^{(+)}_{\hspace{-1mm}k_4} \, {{\bf K}'}^{(+)}_{\hspace{-1mm}k_5} \, {\bf P}_{(ij); \mathcal{D}^{(5)},\mathcal{D}^{(0)} } \\
& = \alpha_\hbar^{-1} \, \Phi^\hbar(\wh{x}^{(1)}) \, {{\bf K}'}^{(+)}_{\hspace{-1mm}k_1} \, 
\Phi^\hbar(\wh{x}^{(2)}) \, {{\bf K}'}^{(+)}_{\hspace{-1mm}k_2} \, 
\Phi^\hbar(-\wh{x}^{(3)})^{-1} \, 
\Phi^\hbar(C^{(3)}(-\wh{x}^{(4)}))^{-1} \\
& \quad \cdot   \, \Phi^\hbar(C^{(3)}C^{(4)}(-\wh{x}^{(5)}))^{-1} \, {{\bf K}'}^{(+)}_{\hspace{-1mm}k_3} \, {{\bf K}'}^{(+)}_{\hspace{-1mm}k_4} \, {{\bf K}'}^{(+)}_{\hspace{-1mm}k_5} \, {\bf P}_{(ij); \mathcal{D}^{(5)},\mathcal{D}^{(0)} } \\
& = \cdots \quad (\mbox{do similarly}).
\end{align*}
Here we are using the symbols $C^{(a)}$ also for the induced map on the self-adjoint operators, by a slight abuse of notation; for example, for a self-adjoint operator $\wh{x} = \pi^\hbar_{\mathcal{D}^{(a)}}(x)$ on $\mathscr{H}_{\mathcal{D}^{(a)}}$, the operator $\pi^\hbar_{\mathcal{D}^{(a-1)}}(C^{(a)}(x))$ on $\mathscr{H}_{\mathcal{D}^{(a-1)}}$ is being denoted by $C^{(a)}(\wh{x})$. Anyways, this way we get
\begin{align*}
& {\bf K}_{\mathcal{D}^{(0)},\mathcal{D}^{(1)}}^{\hbar(+)} \, {\bf K}_{\mathcal{D}^{(1)},\mathcal{D}^{(2)}}^{\hbar(+)} \, {\bf K}_{\mathcal{D}^{(2)},\mathcal{D}^{(3)}}^{\hbar(-)}\, {\bf K}_{\mathcal{D}^{(3)},\mathcal{D}^{(4)}}^{\hbar(-)} \, {\bf K}_{\mathcal{D}^{(4)},\mathcal{D}^{(5)}}^{\hbar(-)} \, {\bf P}_{(ij); \mathcal{D}^{(5)},\mathcal{D}^{(0)} } \\
& = \alpha_\hbar^{-1} \, \Phi^\hbar(\wh{x}^{(1)}) \,\Phi^\hbar(C^{(1)}(\wh{x}^{(2)})) \,
\Phi^\hbar(C^{(1)}C^{(2)}(-\wh{x}^{(3)}))^{-1} \,
\Phi^\hbar(C^{(1)}C^{(2)}C^{(3)}(-\wh{x}^{(4)}))^{-1} \\
& \quad \cdot \Phi^\hbar(C^{(1)}C^{(2)}C^{(3)}C^{(4)}(-\wh{x}^{(5)}))^{-1} \, \,
{{\bf K}'}^{(+)}_{\hspace{-1mm} k_1} \,
{{\bf K}'}^{(+)}_{\hspace{-1mm} k_2} \,
{{\bf K}'}^{(-)}_{\hspace{-1mm} k_3} \,
{{\bf K}'}^{(-)}_{\hspace{-1mm} k_4} \,
{{\bf K}'}^{(-)}_{\hspace{-1mm} k_5} \,
{\bf P}_{(ij); \mathcal{D}^{(5)},\mathcal{D}^{(0)} }.
\end{align*}
The arguments of the five $\Phi^\hbar$-factors are all self-adjoint operators on $\mathscr{H}_{\mathcal{D}^{(0)}}$, and they are as follows (we omit $\mathcal{D}^{(0)}$ from the subscripts of the final resulting operators):
\begin{align*}
& \wh{x}^{(1)} = \wh{x}^\hbar_i, \quad 
C^{(1)}(\wh{x}^{(2)}) = C^{(1)}(\wh{x}^\hbar_{j;\mathcal{D}^{(1)}}) = \wh{x}^\hbar_j + \underbrace{ [\varepsilon^{(0)}_{ji}]_+ \wh{x}^\hbar_i }_{\mbox{\tiny $=0$}} = \wh{x}^\hbar_j , \\
& C^{(1)}C^{(2)}(-\wh{x}^{(3)}) = C^{(1)}C^{(2)}(-\wh{x}^\hbar_{i;\mathcal{D}^{(2)}}) = C^{(1)}(-(\wh{x}^\hbar_{i; \mathcal{D}^{(1)}} + \underbrace{ [\varepsilon^{(1)}_{ij}]_+  \wh{x}^\hbar_{j;\mathcal{D}^{(1)}} }_{\mbox{\tiny $=0$}})) = \wh{x}^\hbar_i, \\
& C^{(1)}C^{(2)}C^{(3)}(-\wh{x}^{(4)}) = C^{(1)}C^{(2)}C^{(3)}(-\wh{x}^\hbar_{j;\mathcal{D}^{(3)}}) = C^{(1)}C^{(2)}(-(\wh{x}^\hbar_{j;\mathcal{D}^{(2)}} + \underbrace{ [-\varepsilon^{(2)}_{ji}]_+ }_{\mbox{\tiny $=1$} } \wh{x}^\hbar_{i;\mathcal{D}^{(2)}}  )) \\
& \qquad = C^{(1)}(\wh{x}^\hbar_{j;\mathcal{D}^{(1)}} - (\wh{x}^\hbar_{i;\mathcal{D}^{(1)}} + \underbrace{ [\varepsilon^{(1)}_{ij}]_+ \wh{x}^\hbar_{j;\mathcal{D}^{(1)}} }_{\mbox{\tiny $=0$}} ) )  = (\wh{x}^\hbar_j  + \underbrace{ [\varepsilon^{(0)}_{ji}]_+ \wh{x}^\hbar_i }_{\mbox{\tiny $=0$} } ) + \wh{x}^\hbar_i  = \wh{x}^\hbar_i + \wh{x}^\hbar_j, \\
& C^{(1)}C^{(2)}C^{(3)}C^{(4)}(-\wh{x}^{(5)}) = C^{(1)} C^{(2)} C^{(3)} C^{(4)} (-\wh{x}^\hbar_{i;\mathcal{D}^{(4)}}) \\
& 
\qquad = C^{(1)}C^{(2)}C^{(3)}(-(\wh{x}^\hbar_{i;\mathcal{D}^{(3)}} + \underbrace{ [-\varepsilon^{(3)}_{ij}]_+}_{\mbox{\tiny $=1$}} \, \wh{x}^\hbar_{j;\mathcal{D}^{(3)}})) \\
& \qquad = C^{(1)}C^{(2)} ( \wh{x}^\hbar_{i;\mathcal{D}^{(2)}} - (\wh{x}^\hbar_{j;\mathcal{D}^{(2)}} + \underbrace{ [-\varepsilon^{(2)}_{ji}]_+ }_{\mbox{\tiny $=1$}} \wh{x}^\hbar_{i;\mathcal{D}^{(2)}}  ) )  = C^{(1)} (  \wh{x}^\hbar_{j;\mathcal{D}^{(1)}}) 
=  \wh{x}^\hbar_j + \underbrace{ [\varepsilon^{(0)}_{ji}]_+ \wh{x}^\hbar_i }_{\mbox{\tiny $=0$}}  = \wh{x}^\hbar_j.
\end{align*}
So the composition of the five $\Phi^\hbar$-factors is
$$
\Phi^\hbar(\wh{x}^\hbar_i) \, \Phi^\hbar(\wh{x}^\hbar_j) \, \Phi^\hbar(\wh{x}^\hbar_i)^{-1} \, \Phi^\hbar(\wh{x}^\hbar_i + \wh{x}^\hbar_j)^{-1} \, \Phi^\hbar(\wh{x}^\hbar_j)^{-1}
$$
which equals ${\rm Id}$ by Prop.\ref{prop:Phi_hbar}(4). It remains to deal with the composition of the five ${\bf K}'$-factors and ${\bf P}_{(ij)}$. Observe
\begin{align*}
& {{\bf K}'}^{(+)}_{\hspace{-1mm} k_1} \,
{{\bf K}'}^{(+)}_{\hspace{-1mm} k_2} \,
{{\bf K}'}^{(-)}_{\hspace{-1mm} k_3} \,
{{\bf K}'}^{(-)}_{\hspace{-1mm} k_4} \,
{{\bf K}'}^{(-)}_{\hspace{-1mm} k_5} \\
& = {{\bf K}'}^{(+)}_{\hspace{-1mm} \mathcal{D}^{(0)},\mathcal{D}^{(1)}} \,
{{\bf K}'}^{(+)}_{\hspace{-1mm} \mathcal{D}^{(1)},\mathcal{D}^{(2)}} \,
{{\bf K}'}^{(-)}_{\hspace{-1mm} \mathcal{D}^{(2)},\mathcal{D}^{(3)}} \,
{{\bf K}'}^{(-)}_{\hspace{-1mm} \mathcal{D}^{(3)},\mathcal{D}^{(4)}} \,
{{\bf K}'}^{(-)}_{\hspace{-1mm} \mathcal{D}^{(4)},\mathcal{D}^{(5)}} \\
& = {\bf F}_{C^{(1)}(\mathcal{D}^{(1)}), \, \mathcal{D}^{(0)}} \, 
{\bf R}_{\mathcal{D}^{(1)}, \, C^{(1)}(\mathcal{D}^{(1)})} \,
{\bf F}_{C^{(2)}(\mathcal{D}^{(2)}), \, \mathcal{D}^{(1)}} \, 
{\bf R}_{\mathcal{D}^{(2)}, \, C^{(2)}(\mathcal{D}^{(2)})} \,
{\bf F}_{C^{(3)}(\mathcal{D}^{(3)}), \, \mathcal{D}^{(2)}} \, 
{\bf R}_{\mathcal{D}^{(3)}, \, C^{(3)}(\mathcal{D}^{(3)})} \\
& \quad \cdot {\bf F}_{C^{(4)}(\mathcal{D}^{(4)}), \, \mathcal{D}^{(3)}} \,
{\bf R}_{\mathcal{D}^{(4)}, \, C^{(4)}(\mathcal{D}^{(4)})} \,
\underbrace{ {\bf F}_{C^{(5)}(\mathcal{D}^{(5)}), \, \mathcal{D}^{(4)}} \, 
{\bf R}_{\mathcal{D}^{(5)}, \, C^{(5)}(\mathcal{D}^{(5)})} }
\end{align*}
We shall repeatedly use Lem.\ref{lem:compatibility_between_F_and_R} to rewrite some ${\bf FR}$ as ${\bf RF}$, then combine ${\bf R}$'s and ${\bf F}$'s using Lem.\ref{lem:R_composition_identity} (eq.\eqref{eq:R_composition_identity}) and Cor.\ref{cor:relationship_between_two_compositions_of_F} (eq.\eqref{eq:F_composition_combine}). First, apply Lem.\ref{lem:compatibility_between_F_and_R} to the underbraced part above, for $(\ell_1',\mathcal{B}_1') = C^{(5)}(\mathcal{D}^{(5)})$, $C = (C^{(5)})^{-1}$, $(\ell_2',\mathcal{B}_2') = \mathcal{D}^{(4)}$, to rewrite it as ${\bf R}_{(C^{(5)})^{-1}(\mathcal{D}^{(4)}),\, \mathcal{D}^{(4)}} \, {\bf F}_{\mathcal{D}^{(5)}, \, (C^{(5)})^{-1}(\mathcal{D}^{(4)})}$, to get
\begin{align*}
& {{\bf K}'}^{(+)}_{\hspace{-1mm} k_1} \,
{{\bf K}'}^{(+)}_{\hspace{-1mm} k_2} \,
{{\bf K}'}^{(-)}_{\hspace{-1mm} k_3} \,
{{\bf K}'}^{(-)}_{\hspace{-1mm} k_4} \,
{{\bf K}'}^{(-)}_{\hspace{-1mm} k_5} \\
& = {\bf F}_{C^{(1)}(\mathcal{D}^{(1)}), \, \mathcal{D}^{(0)}} \, 
{\bf R}_{\mathcal{D}^{(1)}, \, C^{(1)}(\mathcal{D}^{(1)})} \,
{\bf F}_{C^{(2)}(\mathcal{D}^{(2)}), \, \mathcal{D}^{(1)}} \, 
{\bf R}_{\mathcal{D}^{(2)}, \, C^{(2)}(\mathcal{D}^{(2)})} \,
{\bf F}_{C^{(3)}(\mathcal{D}^{(3)}), \, \mathcal{D}^{(2)}} \, 
{\bf R}_{\mathcal{D}^{(3)}, \, C^{(3)}(\mathcal{D}^{(3)})} \\
& \quad \cdot {\bf F}_{C^{(4)}(\mathcal{D}^{(4)}), \, \mathcal{D}^{(3)}} \,
\underbrace{ {\bf R}_{\mathcal{D}^{(4)}, \, C^{(4)}(\mathcal{D}^{(4)})} \,
{\bf R}_{(C^{(5)})^{-1}(\mathcal{D}^{(4)}), \, \mathcal{D}^{(4)}} } \, {\bf F}_{\mathcal{D}^{(5)}, \, (C^{(5)})^{-1}(\mathcal{D}^{(4)})}
\end{align*}
Apply Lem.\ref{lem:R_composition_identity} (eq.\eqref{eq:R_composition_identity}) to replace the underbraced ${\bf RR}$ with ${\bf R} = {\bf R}_{(C^{(5)})^{-1}(\mathcal{D}^{(4)}), \, C^{(4)}(\mathcal{D}^{(4)})}$. Then, we shall apply Lem.\ref{lem:compatibility_between_F_and_R} to this new ${\bf R}$ and the ${\bf F}$-factor that is at the immediate left, etc:
\begin{align*}
& {{\bf K}'}^{(+)}_{\hspace{-1mm} k_1} \,
{{\bf K}'}^{(+)}_{\hspace{-1mm} k_2} \,
{{\bf K}'}^{(-)}_{\hspace{-1mm} k_3} \,
{{\bf K}'}^{(-)}_{\hspace{-1mm} k_4} \,
{{\bf K}'}^{(-)}_{\hspace{-1mm} k_5} \\
& = {\bf F}_{C^{(1)}(\mathcal{D}^{(1)}), \, \mathcal{D}^{(0)}} \, 
{\bf R}_{\mathcal{D}^{(1)}, \, C^{(1)}(\mathcal{D}^{(1)})} \,
{\bf F}_{C^{(2)}(\mathcal{D}^{(2)}), \, \mathcal{D}^{(1)}} \, 
{\bf R}_{\mathcal{D}^{(2)}, \, C^{(2)}(\mathcal{D}^{(2)})} \,
{\bf F}_{C^{(3)}(\mathcal{D}^{(3)}), \, \mathcal{D}^{(2)}} \, 
{\bf R}_{\mathcal{D}^{(3)}, \, C^{(3)}(\mathcal{D}^{(3)})} \\
& \quad \cdot \underbrace{ {\bf F}_{C^{(4)}(\mathcal{D}^{(4)}), \, \mathcal{D}^{(3)}} \,
{\bf R}_{(C^{(5)})^{-1}(\mathcal{D}^{(4)}), \, C^{(4)}(\mathcal{D}^{(4)})}  } _{\mbox{\tiny use Lem.\ref{lem:compatibility_between_F_and_R}: $(\ell_1',\mathcal{B}_1')=C^{(4)}(\mathcal{D}^{(4)})$, $C=(C^{(4)}C^{(5)})^{-1}$, $(\ell_2',\mathcal{B}_2')=\mathcal{D}^{(3)}$ } \hspace{-20mm}} \, {\bf F}_{\mathcal{D}^{(5)}, \, (C^{(5)})^{-1}(\mathcal{D}^{(4)})} \\
& = {\bf F}_{C^{(1)}(\mathcal{D}^{(1)}), \, \mathcal{D}^{(0)}} \, 
{\bf R}_{\mathcal{D}^{(1)}, \, C^{(1)}(\mathcal{D}^{(1)})} \,
{\bf F}_{C^{(2)}(\mathcal{D}^{(2)}), \, \mathcal{D}^{(1)}} \, 
{\bf R}_{\mathcal{D}^{(2)}, \, C^{(2)}(\mathcal{D}^{(2)})} \,
{\bf F}_{C^{(3)}(\mathcal{D}^{(3)}), \, \mathcal{D}^{(2)}}  \\
& \quad\cdot \underbrace{ {\bf R}_{\mathcal{D}^{(3)}, \, C^{(3)}(\mathcal{D}^{(3)})} \, {\bf R}_{(C^{(4)} C^{(5)})^{-1}(\mathcal{D}^{(3)}), \, \mathcal{D}^{(3)}}  }_{\mbox{\tiny use Lem.\ref{lem:R_composition_identity} (eq.\eqref{eq:R_composition_identity})}} \, \underbrace{ {\bf F}_{(C^{(5)})^{-1}(\mathcal{D}^{(4)}), \, (C^{(4)} C^{(5)})^{-1}(\mathcal{D}^{(3)})} \, {\bf F}_{\mathcal{D}^{(5)}, \, (C^{(5)})^{-1}(\mathcal{D}^{(4)})} }_{\mbox{\tiny use Cor.\ref{cor:relationship_between_two_compositions_of_F} (eq.\eqref{eq:F_composition_combine})}} \\
& \sim {\bf F}_{C^{(1)}(\mathcal{D}^{(1)}), \, \mathcal{D}^{(0)}} \, 
{\bf R}_{\mathcal{D}^{(1)}, \, C^{(1)}(\mathcal{D}^{(1)})} \,
{\bf F}_{C^{(2)}(\mathcal{D}^{(2)}), \, \mathcal{D}^{(1)}} \, 
{\bf R}_{\mathcal{D}^{(2)}, \, C^{(2)}(\mathcal{D}^{(2)})}  \\
& \quad \cdot \underbrace{ {\bf F}_{C^{(3)}(\mathcal{D}^{(3)}), \, \mathcal{D}^{(2)}}  \,
{\bf R}_{(C^{(4)}C^{(5)})^{-1}(\mathcal{D}^{(3)}), \, C^{(3)}(\mathcal{D}^{(3)})} }_{\mbox{\tiny use Lem.\ref{lem:compatibility_between_F_and_R}: $(\ell_1',\mathcal{B}_1')=C^{(3)}(\mathcal{D}^{(3)})$, $C=(C^{(3)}C^{(4)}C^{(5)})^{-1}$, $(\ell_2',\mathcal{B}_2')=\mathcal{D}^{(2)}$ } \hspace{-20mm}}  \, {\bf F}_{\mathcal{D}^{(5)}, \, (C^{(4)} C^{(5)})^{-1}(\mathcal{D}^{(3)})} \\
& = {\bf F}_{C^{(1)}(\mathcal{D}^{(1)}), \, \mathcal{D}^{(0)}} \, 
{\bf R}_{\mathcal{D}^{(1)}, \, C^{(1)}(\mathcal{D}^{(1)})} \,
{\bf F}_{C^{(2)}(\mathcal{D}^{(2)}), \, \mathcal{D}^{(1)}} \, 
\underbrace{ {\bf R}_{\mathcal{D}^{(2)}, \, C^{(2)}(\mathcal{D}^{(2)})} \,
{\bf R}_{(C^{(3)}C^{(4)}C^{(5)})^{-1}(\mathcal{D}^{(2)}),\, \mathcal{D}^{(2)}}  }_{\mbox{\tiny use Lem.\ref{lem:R_composition_identity} (eq.\eqref{eq:R_composition_identity})}} \\
& \quad \cdot 
\, \underbrace{ {\bf F}_{(C^{(4)}C^{(5)})^{-1}(\mathcal{D}^{(3)}), \, (C^{(3)}C^{(4)}C^{(5)})^{-1}(\mathcal{D}^{(2)})} \, {\bf F}_{\mathcal{D}^{(5)}, \, (C^{(4)} C^{(5)})^{-1}(\mathcal{D}^{(3)})} }_{\mbox{\tiny use Cor.\ref{cor:relationship_between_two_compositions_of_F} (eq.\eqref{eq:F_composition_combine})}} \\
& \sim {\bf F}_{C^{(1)}(\mathcal{D}^{(1)}), \, \mathcal{D}^{(0)}} \, 
{\bf R}_{\mathcal{D}^{(1)}, \, C^{(1)}(\mathcal{D}^{(1)})} \,
{\bf F}_{C^{(2)}(\mathcal{D}^{(2)}), \, \mathcal{D}^{(1)}} \,
{\bf R}_{(C^{(3)}C^{(4)}C^{(5)})^{-1}(\mathcal{D}^{(2)}),\, C^{(2)}(\mathcal{D}^{(2)})} \\
& \quad \cdot {\bf F}_{\mathcal{D}^{(5)},\,(C^{(3)}C^{(4)}C^{(5)})^{-1}(\mathcal{D}^{(2)})} \\
& \sim \cdots (\mbox{do similarly}) \\
& \sim {\bf R}_{(C^{(1)}C^{(2)}C^{(3)}C^{(4)}C^{(5)})^{-1}(\mathcal{D}^{(0)}), \, \mathcal{D}^{(0)}} \, {\bf F}_{\mathcal{D}^{(5)}, \, (C^{(1)}C^{(2)}C^{(3)}C^{(4)}C^{(5)})^{-1}(\mathcal{D}^{(0)})}.
\end{align*}
In fact, such a manipulation works for any sequence of flips, hence shall be adapted also for proofs of the parts (1) and (2) of Thm.\ref{thm:main} shortly. It remains to compute the linear map $C^{(1)}C^{(2)}C^{(3)}C^{(4)}C^{(5)} : V_{T^{(5)}} \to V_{T^{(0)}}$. It is straightforward to compute this composed linear map by hand. Instead, we refer to \cite[\S2.7,\S3.6,\S4.7]{KN} \cite[\S4.5]{K16c}  about the results on the tropical sign sequences, saying that $C^{(1)}C^{(2)}C^{(3)}C^{(4)}C^{(5)}$ equals the label exchange isomorphism $C_{(i j)} = C_{(ij)}^{-1}$. Thus
\begin{align*}
& {{\bf K}'}^{(+)}_{\hspace{-1mm} k_1} \,
{{\bf K}'}^{(+)}_{\hspace{-1mm} k_2} \,
{{\bf K}'}^{(-)}_{\hspace{-1mm} k_3} \,
{{\bf K}'}^{(-)}_{\hspace{-1mm} k_4} \,
{{\bf K}'}^{(-)}_{\hspace{-1mm} k_5} \,
{\bf P}_{(ij);\mathcal{D}^{(5)},\mathcal{D}^{(0)}} \\
& \sim ~ {\bf R}_{C_{(ij)}(\mathcal{D}^{(0)}), \, \mathcal{D}^{(0)}} \, {\bf F}_{\mathcal{D}^{(5)}, \, C_{(ij)}(\mathcal{D}^{(0)})} \, {\bf P}_{(ij);\mathcal{D}^{(5)},\mathcal{D}^{(0)}} \\
& = {\bf R}_{C_{(ij)}(\mathcal{D}^{(0)}), \, \mathcal{D}^{(0)}} \, \underbrace{ {\bf F}_{\mathcal{D}^{(5)}, \, C_{(ij)}(\mathcal{D}^{(0)})} \, {\bf F}_{C_{(ij)}(\mathcal{D}^{(0)}), \, \mathcal{D}^{(5)}} }_{\mbox{\tiny $={\rm Id}$ ($\because$ Prop.\ref{prop:Weil_intertwiner_involutivity})}} \, {\bf R}_{\mathcal{D}^{(0)}, \, C_{(ij)}(\mathcal{D}^{(0)})} \qquad (\because \mbox{eq.\eqref{eq:permutation_operator_definition}}) \\
& = {\bf R}_{C_{(ij)}(\mathcal{D}^{(0)}), \, \mathcal{D}^{(0)}} \, {\bf R}_{\mathcal{D}^{(0)}, \, C_{(ij)}(\mathcal{D}^{(0)})} ~ \underset{\mbox{\tiny Lem.\ref{lem:R_composition_identity} (eq.\eqref{eq:R_composition_identity})}}{=} ~ {\bf R}_{\mathcal{D}^{(0)},\,\mathcal{D}^{(0)}}.
\end{align*}
It is easy to see from the construction of ${\bf R}$ in \S\ref{subsec:projective_representation_of_the_symplectic_group} that
\begin{align}
\label{eq:R_ell_ell}
{\bf R}_{(\ell,\mathcal{B}),(\ell,\mathcal{B})} = {\rm Id}
\end{align}
holds for any symplectic decomposition, hence ${\bf R}_{\mathcal{D}^{(0)},\mathcal{D}^{(0)}} = {\rm Id}$, finishing the proof of the claim of Thm.\ref{thm:main}.(3) in the case $\varepsilon^{(0)}_{ij}=1$. 

\vs

The case $\varepsilon^{(0)}_{ij}=-1$ can be dealt with by a similar method, but it is wiser to use a slightly different choice of the (tropical) signs, namely $+,+,+,-,-$; then, for example, the composition of the five $\Phi^\hbar$-factors is $\Phi^\hbar(\wh{x}^\hbar_i) \Phi^\hbar(\wh{x}^\hbar_i+\wh{x}^\hbar_j) \Phi^\hbar(\wh{x}^\hbar_j) \Phi^\hbar(\wh{x}^\hbar_i)^{-1} \Phi^\hbar(\wh{x}^\hbar_j)^{-1}$, which equals ${\rm Id}$ by Prop.\ref{prop:Phi_hbar}(4). We omit the detailed computation, and refer the interested readers to \cite[\S4.5,\S5.4]{K16c}. \quad \ul{\it end of proof of the pentagon identity.}

\vs

\vs

\ul{\it The twice flip identity.} Assume the hypothesis of Thm.\ref{thm:main}.(1). Observe
\begin{align*}
& {\bf K}^{\hbar(+)}_{(\ell_T,\mathcal{B}_T),(\ell_{T'},\mathcal{B}_{T'})} \, {\bf K}^{\hbar(-)}_{(\ell_{T'},\mathcal{B}_{T'}),(\ell_T,\mathcal{B}_T)} \\
& = \Phi^\hbar(\wh{x}^\hbar_{k; \ell_T,\mathcal{B}_T} ) \,
\underbrace{ {{\bf K}'}^{(+)}_{\hspace{-1mm}(\ell_T,\mathcal{B}_T),(\ell_{T'},\mathcal{B}_{T'})}   }_{\mbox{\tiny move to right}} \,
\Phi^\hbar(- \wh{x}^\hbar_{k; \ell_{T'},\mathcal{B}_{T'}} )^{-1} \,
{{\bf K}'}^{(-)}_{\hspace{-1mm}(\ell_{T'},\mathcal{B}_{T'}),(\ell_T,\mathcal{B}_T)}  \\
& = \underbrace{ \Phi^\hbar(\wh{x}^\hbar_{k; \ell_T,\mathcal{B}_T} ) \,
\Phi^\hbar(\wh{x}^\hbar_{k; \ell_T,\mathcal{B}_T} )^{-1} }_{\mbox{\tiny cancels}}  \,
{{\bf K}'}^{(+)}_{\hspace{-1mm}(\ell_T,\mathcal{B}_T),(\ell_{T'},\mathcal{B}_{T'})}  \,
{{\bf K}'}^{(-)}_{\hspace{-1mm}(\ell_{T'},\mathcal{B}_{T'}),(\ell_T,\mathcal{B}_T)},
\end{align*}
where in the first step we used eq.\eqref{eq:conjugation_of_primed_part_on_generators} of Lem.\ref{lem:conjugation_of_primed_part}, like in the proof of the pentagon identity.  It remains to show that ${{\bf K}'}^{(+)}_{\hspace{-1mm}(\ell_T,\mathcal{B}_T),(\ell_{T'},\mathcal{B}_{T'})}  \,
{{\bf K}'}^{(-)}_{\hspace{-1mm}(\ell_{T'},\mathcal{B}_{T'}),(\ell_T,\mathcal{B}_T)}$ equals identity, up to $e^{{\rm i}\pi/4}$. By a similar process adapted in the proof of the pentagon identity, using Lem.\ref{lem:compatibility_between_F_and_R}, Cor.\ref{cor:relationship_between_two_compositions_of_F} (or eq.\eqref{eq:F_composition_combine}), and eq.\eqref{eq:R_composition_identity}, we obtain
\begin{align*}
& {{\bf K}'}^{(+)}_{\hspace{-1mm}(\ell_T,\mathcal{B}_T),(\ell_{T'},\mathcal{B}_{T'})}  \,
{{\bf K}'}^{(-)}_{\hspace{-1mm}(\ell_{T'},\mathcal{B}_{T'}),(\ell_T,\mathcal{B}_T)} 
~\sim~ {\bf R}_{(C_k^{(+)} \, {C'_k}^{(-)})^{-1}(\ell_{T'},\mathcal{B}_{T'}), \, (\ell_{T'},\mathcal{B}_{T'})} \, {\bf F}_{(\ell_{T'},\mathcal{B}_{T'}), \, (C_k^{(+)} \, {C_k'}^{(-)} )^{-1}(\ell_{T'},\mathcal{B}_{T'})}
\end{align*}
where $C_k^{(+)}:V_{T'} \to V_T$ and ${C'_k}^{(-)}:V_T \to V_{T'}$ are the linear maps in eq.\eqref{eq:C_k_epsilon} corresponding to the first and the second flips, respectively. One can check $C_k^{(+)} \, {C'_k}^{(-)} = {\rm Id}$; indeed, using \eqref{eq:C_k_epsilon} and \eqref{eq:C_k_epsilon_inverse}, one checks $C_k^{(+)} {C'_k}^{(-)} (x_k) = C_k^{(+)}(-x_k) = x_k$ and $C_k^{(+)} {C'_k}^{(-)} (x_e) = C_k^{(+)} (x_e' + [-\varepsilon'_{ek}]_+ x_k') = x_e + [\varepsilon_{ek}]_+ x_k + [-\varepsilon'_{ek}]_+ (-x_k) = x_e + [\varepsilon_{ek}]_+ (x_k - x_k) = x_e$. As a matter of fact, this is why we chose the sign $+$ for the first flip and $-$ for the second flip, being an example of the `tropical sign sequence'.

\vs

Now, note ${\bf F}_{(\ell_{T'},\mathcal{B}_{T'}), (\ell_{T'},\mathcal{B}_{T'})} = {\rm Id}$ as observed in eq.\eqref{eq:F_ell_ell} (plus the definition of ${\bf F}_{(\ell,\mathcal{B}),(\ell,\mathcal{B})}$ in \S\ref{subsec:LSSW_intertwiner}), and note ${\bf R}_{(\ell_{T'},\mathcal{B}_{T'}),(\ell_{T'},\mathcal{B}_{T'})} = {\rm Id}$ as observed in eq.\eqref{eq:R_ell_ell}. This finishes the proof of Thm.\ref{thm:main}.(1). \quad 

\ul{\it end of proof of the twice flip identity.}

\vspace{5mm}

\ul{\it The square identity.} Assume the hypothesis of Thm.\ref{thm:main}.(2). We follow the notation convention as in the proof of the pentagon identity; so we shall investigate the composition
$$
{\bf K}_{\mathcal{D}^{(0)},\mathcal{D}^{(1)}}^{\hbar(+)} \, {\bf K}_{\mathcal{D}^{(1)},\mathcal{D}^{(2)}}^{\hbar(+)} \, {\bf K}_{\mathcal{D}^{(2)},\mathcal{D}^{(3)}}^{\hbar(-)}\, {\bf K}_{\mathcal{D}^{(3)},\mathcal{D}^{(0)}}^{\hbar(-)}
$$
Like before, we write each ${\bf K}^{\hbar(\epsilon_a)}_{k_a}$ as $\alpha_\hbar^{\epsilon_a} \, \Phi^\hbar(\epsilon_a \, \wh{x}^{(a)})^{\epsilon_a} \, {{\bf K}'}^{(\epsilon_a)}_{\hspace{-1mm} k_a}$, and move all ${{\bf K}'}^{(\epsilon_a)}_{\hspace{-1mm} k_a}$ to the right using eq.\eqref{eq:conjugation_of_primed_part_on_generators} of Lem.\ref{lem:conjugation_of_primed_part};
\begin{align*}
& {\bf K}_{\mathcal{D}^{(0)},\mathcal{D}^{(1)}}^{\hbar(+)} \, {\bf K}_{\mathcal{D}^{(1)},\mathcal{D}^{(2)}}^{\hbar(+)} \, {\bf K}_{\mathcal{D}^{(2)},\mathcal{D}^{(3)}}^{\hbar(-)}\, {\bf K}_{\mathcal{D}^{(3)},\mathcal{D}^{(0)}}^{\hbar(-)} \\
& = \alpha_\hbar^{1+1-1-1} \, \Phi^\hbar(\wh{x}^{(1)}) \,\Phi^\hbar(C^{(1)}(\wh{x}^{(2)})) \,
\Phi^\hbar(C^{(1)}C^{(2)}(-\wh{x}^{(3)}))^{-1} \,
\Phi^\hbar(C^{(1)}C^{(2)}C^{(3)}(-\wh{x}^{(4)}))^{-1} \\
& \quad \cdot
{{\bf K}'}^{(+)}_{\hspace{-1mm} k_1} \,
{{\bf K}'}^{(+)}_{\hspace{-1mm} k_2} \,
{{\bf K}'}^{(-)}_{\hspace{-1mm} k_3} \,
{{\bf K}'}^{(-)}_{\hspace{-1mm} k_4}
\end{align*}
We have $\varepsilon^{(0)}_{ij} = - \varepsilon^{(1)}_{ij} = \varepsilon^{(2)}_{ij} = - \varepsilon^{(3)}_{ij} = 0$, so it is easy to verify
$$
\wh{x}^{(1)} = \wh{x}^\hbar_i, \quad
C^{(1)}(\wh{x}^{(2)}) = \wh{x}^\hbar_j, \quad
(C^{(1)}C^{(2)})(-\wh{x}^{(3)}) = \wh{x}^\hbar_i, \quad
(C^{(1)}C^{(2)}C^{(3)})(-\wh{x}^{(4)}) = \wh{x}^\hbar_j,
$$
thus the composition of the four $\Phi^\hbar$-factors is
$$
\Phi^\hbar(\wh{x}^\hbar_i) \, \Phi^\hbar(\wh{x}^\hbar_j) \, \Phi^\hbar(\wh{x}^\hbar_i)^{-1} \, \Phi^\hbar(\wh{x}^\hbar_j)^{-1};
$$
this equals the identity, because since $\wh{x}^\hbar_i$ strongly commutes with $\wh{x}^\hbar_j$ (i.e. $e^{{\rm i} b \wh{x}^\hbar_i} e^{{\rm i} b' \wh{x}^\hbar_j} = e^{{\rm i} b' \wh{x}^\hbar_j} e^{{\rm i} b \wh{x}^\hbar_i} $, $\forall b,b' \in \mathbb{R}$), the results $\Phi^\hbar(\wh{x}^\hbar_i)$ and $\Phi^\hbar(\wh{x}^\hbar_j)$ of their functional calculus also commute, by a basic result in functional analysis. Meanwhile, by a similar process as in the proof of the pentagon identity adapted to the composition of the four ${\bf K}'$-factors, using Lem.\ref{lem:compatibility_between_F_and_R}, Cor.\ref{cor:relationship_between_two_compositions_of_F} (eq.\eqref{eq:F_composition_combine}), and Lem.\ref{lem:R_composition_identity} (eq.\eqref{eq:R_composition_identity}), we obtain
\begin{align*}
{{\bf K}'}^{(+)}_{\hspace{-1mm} k_1} \,
{{\bf K}'}^{(+)}_{\hspace{-1mm} k_2} \,
{{\bf K}'}^{(-)}_{\hspace{-1mm} k_3} \,
{{\bf K}'}^{(-)}_{\hspace{-1mm} k_4}
~\sim~ {\bf R}_{(C^{(1)}C^{(2)}C^{(3)}C^{(4)})^{-1}(\mathcal{D}^{(0)}) , \, \mathcal{D}^{(0)}} \, {\bf F}_{\mathcal{D}^{(0)}, \, (C^{(1)}C^{(2)}C^{(3)}C^{(4)})^{-1}(\mathcal{D}^{(0)}) }
\end{align*}
It remains to compute the linear map $C^{(1)}C^{(2)}C^{(3)}C^{(4)} : V_{T^{(0)}} \to V_{T^{(0)}}$. One can verify that this equals ${\rm Id}$ either by hand, or by using the results on tropical sign sequences \cite[\S4.5]{K16c}. Anyhow, we get
$$
{\bf K}_{\mathcal{D}^{(0)},\mathcal{D}^{(1)}}^{\hbar(+)} \, {\bf K}_{\mathcal{D}^{(1)},\mathcal{D}^{(2)}}^{\hbar(+)} \, {\bf K}_{\mathcal{D}^{(2)},\mathcal{D}^{(3)}}^{\hbar(-)}\, {\bf K}_{\mathcal{D}^{(3)},\mathcal{D}^{(0)}}^{\hbar(-)}
~\sim~ {\bf R}_{\mathcal{D}^{(0)}, \, \mathcal{D}^{(0)}} \, {\bf F}_{\mathcal{D}^{(0)}, \, \mathcal{D}^{(0)}} = {\rm Id},
$$
as desired. \ul{\it end of proof of the square identity.}

\vs

\ul{\it The `trivial' identities.} ${\bf FF} \sim {\bf F}$ and ${\bf F} \, {\bf K}_k^{\hbar(\epsilon)} \sim {\bf K}_k^{\hbar(\epsilon)} \, {\bf F}$ are already proved in Cor.\ref{cor:relationship_between_two_compositions_of_F} and Prop.\ref{prop:compatibility_of_mutation_intertwiner_and_Weil_intertwiners}. Other identities are straightforward exercises about ${\bf F}$'s and ${\bf R}$'s, using Cor.\ref{cor:relationship_between_two_compositions_of_F}, Lem.\ref{lem:R_composition_identity}, and Lem.\ref{lem:compatibility_between_F_and_R}. Here we only present a sketch of proof of the most nontrivial among them, namely ${\bf P}_\sigma \, {\bf K}_{\sigma(k)}^{\hbar(\epsilon)} \, {\bf P}_{\sigma^{-1}} \sim {\bf K}_{k}^{\hbar(\epsilon)}$. As a preparation, we investigate the conjugation action of the permutation operator ${\bf P}_\sigma$ on the generating self-adjoint operators $\wh{x}$. Let $[(T,\mathcal{D}),(T',\mathcal{D}')] = P_\sigma$ in ${\rm SPt}(S)$; so the corresponding operator is ${\bf P}_\sigma = {\bf P}_{\sigma;\mathcal{D},\mathcal{D}'} = {\bf F}_{C_\sigma(\mathcal{D}'),\,\mathcal{D}} \, {\bf R}_{\mathcal{D}', \, C_\sigma(\mathcal{D}')}$. Consider a basis vector $x'_{\sigma(i)} \in V_{T'}$ for any edge index $i \in I$, and its corresponding self-adjoint operator $\wh{x'}^\hbar_{\hspace{-1mm}\sigma(i)} = \wh{x'}^\hbar_{\hspace{-1mm}\sigma(i); \mathcal{D}'}= \pi^\hbar_{\mathcal{D}'}(x'_{\sigma(i)})$. Recall from eq.\eqref{eq:C_sigma} that $C_\sigma:V_{T'} \to V_T$ sends $x'_{\sigma(i)}$ to $x_i$, so from Cor.\ref{cor:R_is_pullback} applied to $C = C_\sigma^{-1}$ we have
\begin{align*}
{\bf R}_{\mathcal{D}', \, C_\sigma(\mathcal{D}')} \, \wh{x'}^\hbar_{\hspace{-1mm}\sigma(i);\mathcal{D}'} \, ({\bf R}_{\mathcal{D}', \, C_\sigma(\mathcal{D}')})^{-1} = \wh{x}_{i; C_\sigma(\mathcal{D}')}
\end{align*}
and from eq.\eqref{eq:conjugation_by_F} we have
$$
{\bf F}_{C_\sigma(\mathcal{D}'),\,\mathcal{D}} \, \wh{x}_{i; C_\sigma(\mathcal{D}')} \, ({\bf F}_{C_\sigma(\mathcal{D}'),\,\mathcal{D}})^{-1} = \wh{x}_{i;\mathcal{D}}.
$$
Combining, we get the sought-for conjugation action of ${\bf P}_\sigma$:
\begin{align}
\label{eq:conjugation_action_of_P_sigma}
{\bf P}_{\sigma;\mathcal{D},\mathcal{D}'} \, \wh{x'}^\hbar_{\hspace{-1mm}\sigma(i);\mathcal{D}'} \, ({\bf P}_{\sigma;\mathcal{D},\mathcal{D}'})^{-1} = \wh{x}^\hbar_{i;\mathcal{D}}, \qquad \forall i \in I,
\end{align}
which is exactly what we expect from the label permutation operator ${\bf P}_\sigma$.

\vs

Let's show ${\bf K}_{k}^{\hbar(\epsilon)}  \, {\bf P}_\sigma \sim {\bf P}_\sigma \, {\bf K}_{\sigma(k)}^{\hbar(\epsilon)}$. Let $[T^{(a)},\mathcal{D}^{(a)}]$, $a=0,1,2,3,4$, be objects of ${\rm Spt}(S)$ with
\begin{align*}
\hspace{-2mm}
[(T^{(0)},\mathcal{D}^{(0)}),(T^{(1)},\mathcal{D}^{(1)})] = \mu_k, \quad
[(T^{(1)},\mathcal{D}^{(1)}),(T^{(2)},\mathcal{D}^{(2)})] = P_\sigma, \quad
[(T^{(0)},\mathcal{D}^{(0)}),(T^{(3)},\mathcal{D}^{(3)})] = P_\sigma.
\end{align*}
Then $[(T^{(3)},\mathcal{D}^{(3)}),(T^{(2)},\mathcal{D}^{(2)})] = \mu_{\sigma(k)}$. Note
\begin{align*}
{\bf K}_{k}^{\hbar(\epsilon)}  \, {\bf P}_\sigma  & = {\bf K}_{\mathcal{D}^{(0)},\mathcal{D}^{(1)}}^{\hbar(\epsilon)} \,
{\bf P}_{\sigma; \mathcal{D}^{(1)},\mathcal{D}^{(2)} } \\
& = \alpha_\hbar^{\epsilon} \, \Phi^\hbar(\epsilon \, \wh{x}^\hbar_{k;\mathcal{D}^{(0)}})^\epsilon \, {\bf F}_{C_k^{(\epsilon)}(\mathcal{D}^{(1)}), \, \mathcal{D}^{(0)}} \, {\bf R}_{\mathcal{D}^{(1)}, \, C_k^{(\epsilon)}(\mathcal{D}^{(1)})} \, \underbrace{ {\bf F}_{C_\sigma(\mathcal{D}^{(2)}),\,\mathcal{D}^{(1)}} \, {\bf R}_{\mathcal{D}^{(2)}, \, C_\sigma(\mathcal{D}^{(2)})} }_{\mbox{\tiny use Lem.\ref{lem:compatibility_between_F_and_R}}} \\
& = \alpha_\hbar^{\epsilon} \, \Phi^\hbar(\epsilon \, \wh{x}^\hbar_{k;\mathcal{D}^{(0)}})^\epsilon \, {\bf F}_{C_k^{(\epsilon)}(\mathcal{D}^{(1)}), \, \mathcal{D}^{(0)}} \, 
\underbrace{ {\bf R}_{\mathcal{D}^{(1)}, \, C_k^{(\epsilon)}(\mathcal{D}^{(1)})} \, 
{\bf R}_{C_\sigma^{-1}(\mathcal{D}^{(1)}), \, \mathcal{D}^{(1)}} }_{\mbox{\tiny use Lem.\ref{lem:R_composition_identity} (eq.\eqref{eq:R_composition_identity})} } \,
{\bf F}_{\mathcal{D}^{(2)}, \, C_\sigma^{-1}(\mathcal{D}^{(1)})} \\
& = \alpha_\hbar^{\epsilon} \, \Phi^\hbar(\epsilon \, \wh{x}^\hbar_{k;\mathcal{D}^{(0)}})^\epsilon \, 
\underbrace{ {\bf F}_{C_k^{(\epsilon)}(\mathcal{D}^{(1)}), \, \mathcal{D}^{(0)}} \, 
{\bf R}_{C_\sigma^{-1}(\mathcal{D}^{(1)}), \, C_k^{(\epsilon)}(\mathcal{D}^{(1)})} }_{\mbox{\tiny use Lem.\ref{lem:compatibility_between_F_and_R} }}  \,
{\bf F}_{\mathcal{D}^{(2)}, \, C_\sigma^{-1}(\mathcal{D}^{(1)})} \\
& = \alpha_\hbar^{\epsilon} \, \Phi^\hbar(\epsilon \, \wh{x}^\hbar_{k;\mathcal{D}^{(0)}})^\epsilon \,
{\bf R}_{(C_k^{(\epsilon)} C_\sigma )^{-1}(\mathcal{D}^{(0)}), \, \mathcal{D}^{(0)}}  \,
\underbrace{ {\bf F}_{C_\sigma^{-1}(\mathcal{D}^{(1)}), \, (C_k^{(\epsilon)} C_\sigma )^{-1}(\mathcal{D}^{(0)})} \,
{\bf F}_{\mathcal{D}^{(2)}, \, C_\sigma^{-1}(\mathcal{D}^{(1)})} }_{\mbox{\tiny use Cor.\ref{cor:relationship_between_two_compositions_of_F} (eq.\eqref{eq:F_composition_combine}) }} \\
& \sim \alpha_\hbar^{\epsilon} \, \Phi^\hbar(\epsilon \, \wh{x}^\hbar_{k;\mathcal{D}^{(0)}})^\epsilon \,
{\bf R}_{(C_k^{(\epsilon)} C_\sigma )^{-1}(\mathcal{D}^{(0)}), \, \mathcal{D}^{(0)}}  \,
{\bf F}_{D^{(2)}, \, (C_k^{(\epsilon)}C_\sigma)^{-1}(\mathcal{D}^{(0)})}
\end{align*}
while
\begin{align*}
{\bf P}_\sigma \, {\bf K}_{\sigma(k)}^{\hbar(\epsilon)} \, & = {\bf P}_{\sigma; \mathcal{D}^{(0)},\mathcal{D}^{(3)} } \, {\bf K}_{\mathcal{D}^{(3)},\mathcal{D}^{(2)}}^{\hbar(\epsilon)}
 \\
 & = \alpha_\hbar^{\epsilon} \, {\bf P}_{\sigma; \mathcal{D}^{(0)},\mathcal{D}^{(3)} } \, \Phi^\hbar(\epsilon \, \wh{x}^\hbar_{\sigma(k);\mathcal{D}^{(3)}})^\epsilon \, {\bf F}_{C_{\sigma(k)}^{(\epsilon)}(\mathcal{D}^{(2)}), \, \mathcal{D}^{(3)}} \, {\bf R}_{\mathcal{D}^{(2)}, \, C_{\sigma(k)}^{(\epsilon)}(\mathcal{D}^{(2)})} \\
& = \alpha_\hbar^{\epsilon} \, \Phi^\hbar(\epsilon \, \wh{x}^\hbar_{k;\mathcal{D}^{(0)}})^\epsilon \, {\bf P}_{\sigma; \mathcal{D}^{(0)},\mathcal{D}^{(3)} } \, {\bf F}_{C_{\sigma(k)}^{(\epsilon)}(\mathcal{D}^{(2)}), \, \mathcal{D}^{(3)}} \, {\bf R}_{\mathcal{D}^{(2)}, \, C_{\sigma(k)}^{(\epsilon)}(\mathcal{D}^{(2)})} \quad (\because\mbox{eq.\eqref{eq:conjugation_action_of_P_sigma}}) \\
& = \alpha_\hbar^{\epsilon} \, \Phi^\hbar(\epsilon \, \wh{x}^\hbar_{k;\mathcal{D}^{(0)}})^\epsilon \, {\bf F}_{C_\sigma(\mathcal{D}^{(3)}),\,\mathcal{D}^{(0)}} \, {\bf R}_{\mathcal{D}^{(3)}, \, C_\sigma(\mathcal{D}^{(3)})} \, 
\underbrace{ 
{\bf F}_{C_{\sigma(k)}^{(\epsilon)}(\mathcal{D}^{(2)}), \, \mathcal{D}^{(3)}} \, {\bf R}_{\mathcal{D}^{(2)}, \, C_{\sigma(k)}^{(\epsilon)}(\mathcal{D}^{(2)})} 
}_{\mbox{\tiny use Lem.\ref{lem:compatibility_between_F_and_R}}}
 \\
& = \cdots (\mbox{do similarly using Cor.\ref{cor:relationship_between_two_compositions_of_F}, Lem.\ref{lem:R_composition_identity}, and Lem.\ref{lem:compatibility_between_F_and_R}}) \\
& \sim  \alpha_\hbar^{\epsilon} \, \Phi^\hbar(\epsilon \, \wh{x}^\hbar_{k;\mathcal{D}^{(0)}})^\epsilon \,
{\bf R}_{(C_\sigma C_{\sigma(k)}^{(\epsilon)})^{-1}(\mathcal{D}^{(0)}), \, \mathcal{D}^{(0)}} \,
{\bf F}_{\mathcal{D}^{(2)}, \, (C_\sigma C_{\sigma(k)}^{(\epsilon)})^{-1}(\mathcal{D}^{(0)})}
\end{align*}
One can finish proving ${\bf K}_{k}^{\hbar(\epsilon)}  \, {\bf P}_\sigma \sim {\bf P}_\sigma \, {\bf K}_{\sigma(k)}^{\hbar(\epsilon)}$ by observing $C_k^{(\epsilon)} C_\sigma = C_\sigma C_{\sigma(k)}^{(\epsilon)}$. \qed

\section{The conclusion and beyond}

\subsection{The projective representations of the mapping class group}

Consider the mapping class group ${\rm MCG}(S) = {\rm Homeo}^+(S)/{\rm Homeo}(S)_0$ of $S$, i.e. the group of all orientation-preserving homeomorphisms $S \to S$ up to isotopy. Then ${\rm MCG}(S)$ naturally acts on the set of all labeled ideal triangulations of $S$; this action can be defined on the Ptolemy groupoid ${\rm Pt}(S)$, if we define the action on the morphisms as $g.[T,T'] := [gT, gT']$, for each $g\in {\rm MCG}(S)$. In fact, this action also extends naturally to the symplectic Ptolemy groupoid ${\rm SPt}(S)$. Indeed, for any $g\in {\rm MCG}(S)$ and any object $(T,\ell_T,\mathcal{B}_T)$ of ${\rm SPt}(S)$, we have a natural symplectic isomorphism $C_g : V_T \to V_{gT}$ sending each basis vector $x_e$ to $x_{g(e)}$, and this isomorphism induces a symplectic decomposition $(\ell_{gT},\mathcal{B}_{gT}) := C_g(\ell_T,\mathcal{B}_T)$ of $(V_{gT}, B_{gT})$; then, define the action on objects as $g.(T,\ell_T,\mathcal{B}_T) := (gT, \ell_{gT},\mathcal{B}_{gT})$, and that on morphisms as $g[(T,\ell_T,\mathcal{B}_T),(T',\ell_{T'},\mathcal{B}_{T'})] := [(gT,\ell_{gT},\mathcal{B}_{gT}),(gT',\ell_{gT'},\mathcal{B}_{gT'})]$. Then we can easily observe:
\begin{proposition}
The projective functor $\pi:{\rm SPt}(S) \to {\rm Hilb}$ constructed in \S\ref{subsec:projective_representation_of_SPt_S} is invariant under the ${\rm MCG}(S)$-action on ${\rm SPt}(S)$. \qed
\end{proposition}
Indeed, the Hilbert spaces and the operators in the present paper used only the combinatorics of the ideal triangulations $T$, but did not really use the full information of $T$ as an ideal triangulation. So, we may formulate our result as the projective functor 
$$
\ol{\pi } = \ol{\pi}^\hbar_\lambda : \ol{\rm SPt}(S)/{\rm MCG}(S) \to {\rm Hilb}
$$
from the quotient category
$$
\ol{\rm SPt}(S) := {\rm SPt}(S)/{\rm MCG}(S)
$$
In fact, this is what had really been aimed for in the quantum Teichm\"uller theory, so we might just have formulated the problem this way from the beginning. By construction, for any object $O$ of $\ol{\rm SPt}(S)$, the automorphism group ${\rm Aut}(O)$ is isomorphic to ${\rm MCG}(S)$, so $\ol{\pi}$ yields a projective unitary representation
$$
\ol{\pi}^\hbar_O = \ol{\pi}^\hbar_{O;\lambda} : {\rm MCG}(S) \to {\rm U}(\mathscr{H}_O)
$$
of the mapping class group ${\rm MCG}(S)$ on the Hilbert space $\mathscr{H}_O \cong L^2(\mathbb{R}^r)$. It is also clear that the projective representations $\ol{\pi}^\hbar_O$ for different $O$'s are unitarily equivalent. The main result of the present paper is then an explicit construction of these representations, as well as the computation of the phase constants coming from the projective-ness of the representations. Namely, for each $g_1,g_2 \in {\rm MCG}(S)$, we have
$$
\ol{\pi}^\hbar_O(g_1g_2) = c^\hbar_{g_1, g_2} \, \ol{\pi}^\hbar_O(g_1) \, \ol{\pi}^\hbar_O(g_2) 
$$
for some constants $c^\hbar_{g_1,g_2} \in {\rm U}(1)$, and we computed these constants.

\vs

These constants ${\rm MCG}(S) \times {\rm MCG}(S) \to {\rm U}(1)$, $(g_1,g_2) \mapsto c_{g_1,g_2}^\hbar$, satisfy the group $2$-cocycle condition, hence represents a class in the second group cohomology $H^2({\rm MCG}(S);{\rm U}(1))$, which is in one-to-one correspondence with the central extensions of ${\rm MCG}(S)$ by ${\rm U}(1)$. This cohomology class of ${\rm MCG}(S)$ coming from various versions of quantum Teichm\"uller theories has been interesting objects of study, as studied in \cite{FS} \cite{FuKa} \cite{K16a} \cite{K16c}  \cite{Xu}. The non-triviality of this cohomology class for the Chekhov-Fock-Goncharov type quantum Teichm\"uller space has only been hinted and conjectured, and used in \cite{FS} \cite{Xu} without proof. The present paper can be used to finally fill in this gap.

\vs

However, as seen in the present paper, for the phase constants, not only the special number $\alpha_\hbar$ defined in eq.\eqref{eq:alpha_hbar} is involved, but also powers of $e^{{\rm i}\pi/4}$ appear, from the theory of the Shale-Weil intertwiners; this phenomenon has not been mentioned at all in the quantum Teichm\"uller theory literature. So, to verify the asserted results of the previous papers \cite{FS} \cite{Xu}, one must really construct the representations of ${\rm MCG}(S)$ in each of the asserted cases and check the relations, this time not ignoring but keeping precise track of $e^{{\rm i}\pi/4}$ too. Here is how to do it. Pick any object $(T,\ell_T,\mathcal{B}_T)$ of ${\rm SPt}(S)$, representing an object $O$ of $\ol{\rm SPt}(S)$. It is known that ${\rm MCG}(S)$ is finitely presented; pick generators $g_1,\ldots,g_N$. For each generator $g_i$, consider the morphism $[(T,\ell_T,\mathcal{B}_T), g_i(T,\ell_T,\mathcal{B}_T)]$ of ${\rm SPt}(S)$, express it as a composition of the elementary morphisms, and form the composition of the corresponding operators ${\bf K}, {\bf F}, {\bf P}$, to get a unitary operator $\wh{g}_i = \ol{\pi}^\hbar_O(g_i)$ on $\mathscr{H}_{\ell_T,\mathcal{B}_T}$. Each relation among $g_1,\ldots,g_N$ is satisfied by the operators $\wh{g}_1,\ldots,\wh{g}_N$ up to a phase constant in ${\rm U}(1)$; these constants yield a (presentation of a) central extension of ${\rm MCG}(S)$, representing the sought-for second group cohomology class. Note that these phase constants will be product of integer powers of $\alpha_\hbar$ and those of $e^{{\rm i}\pi/4}$. We recall the well-known fact in the theory of the Shale-Weil projective representations of the symplectic groups that the intertwiners ${\bf F}$ can be scaled suitably so that the phase constants are $\pm 1$, and also that these constants cannot be all $1$; see original papers on the Shale-Weil theory \cite{Shale} \cite{Segal} \cite{Weil}, as well as \cite{Rao}. In our case of the projective representations of ${\rm MCG}(S)$ just mentioned, it may be possible to fine-tune the scalars so that the phase constants are powers of just $\alpha_\hbar$; we leave this for a future research, which should be related to comparing $H^2({\rm MCG}(S);\mathbb{Z})$ and $H^2({\rm MCG}(S);\mathbb{Z} \oplus \mathbb{Z}/2\mathbb{Z})$. Another hurdle left is that one should think of how to apply the present work to the cases of non-finite-type surfaces as studied in \cite{FS} \cite{K16a}.

\vs

Besides this story of the phase constants, another important question is that the ${\rm MCG}(S)$ representation constructed above using the results of the present paper is irreducible as a representation of the group ${\rm MCG}(S)$. If it were true, we expect that a similar argument as in the proofs of Stone-von Neumann theorem (e.g. \cite[Thm.14.8]{Hall}) is needed.

\subsection{A generalization to the cluster $\mathscr{X}$-varieties}

The present paper is written in terms of the Teichm\"uller space of a punctured surface, using the ideal triangulations of the surface. In fact, the whole construction could be applied to the skew-symmetric cluster $\mathscr{X}$-varieties of Fock-Goncharov. A labeled ideal triangulation $T$ corresponds to a labeled $\mathscr{X}$-seed $(\varepsilon, \{X_1,\ldots,X_n\})$, consisting of an $n \times n$ skew-symmetric integer matrix $\varepsilon = (\varepsilon_{ij})_{i,j=1}^n$, with the cluster $\mathscr{X}$-variables $X_1,\ldots,X_N$, which form a transcendental basis of some ambient field $\mathcal{F}$. By the {\em mutation} $\mu_k$ along $k\in I = \{1,\ldots,n\}$ one can transform a seed $T = (\varepsilon,\{X_1,\ldots,X_n\})$ into another seed $\mu_k(T) = T'= (\varepsilon',\{X_1',\ldots,X_n'\})$, related to the previous one by the formulas in eq.\eqref{eq:mutation_formula_for_X}. By the {\em seed automorphism} $P_\sigma$ for a permutation $\sigma$ of $I$, one can transform a seed $T$ into $P_\sigma(T)=T'$ using the formulas $\varepsilon'_{\sigma(i)\sigma(j)} = \varepsilon_{ij}$, $X'_{\sigma(i)} = X_i$. Start from one seed, and collect all seeds obtained by sequences of the mutations and the seed automorphisms. By gluing the split algebraic tori ${\rm Spec}(\mathbb{Z}[X_1^{\pm 1},\ldots,X_n^{\pm 1}])$ for these seeds along the birational maps corresponding to $\mu_k$ and $P_\sigma$, one obtains the cluster $\mathscr{X}$-variety. Replacing the commutative rings $\mathbb{Z}[X_1^{\pm 1},\ldots,X_n^{\pm 1}]$ by the non-commutative algebras $\mathcal{X}^q_T$ (Def.\ref{def:CF_algebra}), related with one another by the quantum mutations $\Phi^q_{TT'}$ and the quantum seed automorphisms $P_\sigma$, yields a quantum cluster $\mathscr{X}$-variety \cite{FG09}.

\vs

One can consider the groupoid of seeds, generated by the elementary morphisms corresponding to the mutations and the seed automorphisms, analogous to ${\rm Pt}(S)$. To apply the present paper's construction, for each seed $T$ one must choose a symplectic decomposition $(\ell_T,\mathcal{B}_T)$ of the vector space $(V_T,B_T)$ equipped with a skew-symmetric bilinear form. Hence one can define a groupoid of symplectic decompositions, analogous to ${\rm SPt}(S)$. Using the present paper's construction, one would assign the unitary operators ${\bf K}_k$, ${\bf F}$, and ${\bf P}_\sigma$, to the elementary morphisms of this groupoid, which would yield a family of {\em irreducible} representations of the quantum cluster $\mathscr{X}$-variety. 

\vs

However, at this moment, such a construction of irreducible representations of a general skew-symmetric quantum cluster $\mathscr{X}$-variety is only conjectural, for the following reasons. First, for the final projective representations to be well-defined, the operators ${\bf K}_k$, ${\bf F}$, and ${\bf P}_\sigma$ must satisfy the operator identities, up to multiplicative constants, corresponding to the (classical) identities satisfied by the elementary morphisms. One problem is that the complete list of identities satisfied by the (classical) elementary morphisms is not settled yet. There are identities satisfied by mutations that are not consequences of the ones mentioned in the present paper; see \cite[Examples 4.1, 4.2]{KY} \cite[Rem.9.19]{FST} e.g. for the relation consisting of $32$ mutations. In principle, one must first find all classical relations, and then check these relations for the operators. If one is interested in just proving the well-definedness, but not so much in computing the phase constants, one might try mimicking an indirect argument of Fock-Goncharov \cite{FG09} using the (strong) irreducibility of the constructed representations. Second, for the case of ideal triangulations of a punctured surface, there was a natural preferred basis of the radical $V_T^\perp$ of $(V_T,\mathcal{B}_T)$ enumerated by the punctures of $S$, which is compatible under the `cluster' linear maps $C_k^{(\epsilon)}$ (Lem.\ref{lem:C_k_epsilon_is_symplectic}(2)); in particular, the functions $\til{f} : \mathfrak{h}_T \to \mathbb{R}$ for all possible ideal triangulations $T$ were conveniently controlled by a single function $f:\mathcal{P} \to \mathbb{R}$ on the set of punctures. For a general cluster $\mathscr{X}$-variety, there is no preferred basis of $V_T^\perp$ (or the kernel of the matrix $\varepsilon$), so one must come up with a basis of $V_T^\perp$ compatible under the maps $C_k^{(\epsilon)}$. If these two obstacles can be overcome, we expect that it is also possible to incorporate the skew-symmetrizable cluster $\mathscr{X}$-varieties too, in the style of \cite{FG09} \cite{K16c}.

\end{document}